\documentclass[onefignum,onetabnum]{siamart190516}

\usepackage{algpseudocode}
\usepackage{amssymb}
\usepackage{amsmath}
\usepackage{bbm}
\usepackage{tikz}
\usepackage{pgfplots}
\usepackage{subfig}
\usepackage{todonotes}	
\usepackage{hyperref}

\usepgfplotslibrary{fillbetween}

\usetikzlibrary{decorations.text}
\usetikzlibrary{shapes}
\usetikzlibrary{trees}
\usetikzlibrary{calc}

\newcommand{\R}{\mathbb{R}}
\newcommand{\Oad}{\mathcal{O}^{\text{ad}}}

\newcommand{\kk}[1]{\textcolor{black}{#1}}
\newcommand{\hg}[1]{\textcolor{black}{#1}}
\newcommand{\otd}[1]{\textcolor{black}{#1}}

\newsiamremark{remark}{Remark}

\title{
Non-convex shape optimization by dissipative Hamiltonian flows}

\author{M.\ Bolten\thanks{IMACM, School of Mathematics and Natural Science, University of Wuppertal, D-42119 Wuppertal (\email{bolten@math.uni-wuppertal.de},\email{doganay@math.uni-wuppertal.de},\email{hgottsch@uni-wuppertal.de},\email{klamroth@math.uni-wuppertal.de})} \and O.\ T.\ Doganay\footnotemark[1] \and H.\ Gottschalk\footnotemark[1] \and K. Klamroth\footnotemark[1]}

\begin{document}

\maketitle

\begin{abstract}
Shape optimization with constraints given by partial differential equations (PDE) is a highly developed field of optimization theory. The elegant adjoint formalism allows to compute shape gradients at the computational  cost of a further PDE solve. Thus, gradient descent methods can be applied to shape optimization problems. However, gradient descent methods that can be understood as approximation to gradient flows get stuck in local minima, if the optimization problem is non-convex. In machine learning, the optimization in high dimensional non-convex energy landscapes has been successfully tackled by momentum methods, which can be understood as passing from gradient flow to dissipative Hamiltonian flows. In this paper, we adopt this strategy for non-convex shape optimization. In particular, we provide a mechanical shape optimization problem that is motivated by optimal reliability considering also material cost and the necessity to avoid certain obstructions in installation space. We then show how this problem can be solved effectively by port Hamiltonian shape flows.  
\end{abstract}



\section{Introduction}
Shape optimization is an active and interdisciplinary field in engineering and mathematics \cite{bucur,delfour,haslinger,sokolowski}. In many applications,  shapes occur as the domain of a partial differential equation (PDE) that models physical phenomena and the fitness of a shape depends on the solution of the PDE, which is also called the state equation \cite{troelsch}. As the numerical solution to the PDE is often compute-intensive, a straightforward computation of shape sensitivities by finite difference methods often comes with prohibitive computational cost. However, the elegant adjoint formalism, in its continuous \cite{bucur,delfour,haslinger,sokolowski} or discrete \cite{Frey2009,giles2003algorithm,gottschalk2019shape,gottschalk2,gottschalk3} variants \kk{(see} \cite{nadarajah2000comparison} \kk{for a comparison),} permits the computation of shape gradients with one additional PDE solve, only. This has been exploited in numerous works to optimize the shape of mechanical components \cite{Hahn2019,bolten2015,bucur,chenais}, see also \cite{bolten2021tracing,Doganay2019,gottschalk2021analytical,haslinger} for first steps in the direction of multi-criteria shape optimization.

Looking at the optimization strategies applied, the gradient information is either used in the gradient descent algorithm \cite{wright1999numerical} \kk{or} (pseudo) Newton methods. As an alternative, one can harness global surrogate models with gradient information and then perform surrogate based optimization, as it is done with gradient enhanced Kriging (GEK),  and then apply the EGO search heuristics \cite{backhaus2012gradient,sobester2008engineering}. However, both approaches are beset with  certain limitations: Gradient descent or Newton methods are likely to get stuck in local minima if the shape optimization problem has non-convex characteristics, whereas GEK scales badly in high dimensional search spaces, which are typical for shape optimization.

In machine learning, non-convex optimization problems in extremely high dimension and with complex energy landscapes are solved during \kk{neural network} training \cite{goodfellow2016deep}. Solving here has to be understood not as \kk{necessarily aiming at} the convergence to a global optimum, but rather \kk{aiming at the} convergence to a local minimum of the loss function, at which the model performs sufficiently well. However, this typically is not the \kk{``nearest''} local minimum. 
To overcome \kk{unfavorable local minima that do not achieve satisfactory objective values}, momentum based methods are used \cite{goh2017momentum,kova:cont:2021}.

In this context, it has been proposed to understand momentum as a physical momentum as in classical Hamiltonian mechanics where objective or loss functions assume the role of potential energy, \kk{hence interpreting the trajectory of the solution during optimization as a heavy ball with friction (HBF) \cite{anti:seco:1993,atto:theh:2000,poly:some:1964}. See also \cite{ochs:loca:2018,ochs:unif:2019,ochs:ipia:2014} for variants tailored for non-convex and non-smooth problems, and \cite{atto:mult:2015,sonn:mult:2022} for multiobjective versions. 
Convergence properties have been discussed, among others, in \cite{alva:onth:2000,hara:conv:1998}.
We note that the HBF dynamic can be interpreted \cite{ochs:adap:2019,su:adif:2015}
as a continuous version of the ``fast iterative shrinkage-thresholding algorithm'' (FISTA) \cite{beck:afas:2009}. 
It has been observed \cite{cabo:onth:2009} that an asymptotically vanishing damping effect makes the momentum more effective asymptotically, see also \cite{atto:fast:2018,cham:onth:2015} for corresponding convergence results. A detailed analysis of the interrelation between the continuous (Hamiltonian) systems and their discretized versions can be found in \cite{alv:asym:2010,alv:auni:2011,peyp:evol:2010}.
}
This has been recently cast \cite{kova:cont:2021,massaroli2019port,poli2020port} in the port Hamiltonian language \cite{van2014port}. It has been observed in these works that the stationary points of such Hamiltonian systems are in one-to-one correspondence with the critical points of the original optimization problem. Therefore, if the Hamiltonian system is made dissipative by introduction of Newtonian friction, it should ultimately settle to one of these stationary points and thereby efficiently solve the optimization problem in a less local way as pure gradient descent flow. The actual optimization algorithms here are understood as discretized gradient or Hamiltonian flows.   

In this work, we adopt this strategy and apply it to shape optimization problems. In order to do so, we first propose new shape optimization test cases, which are manifestly non-convex in the sense that they contain non-optimal local minima. While essentially no-one would believe that shape optimization is convex in general, there seem to be very little explicit and well understood \kk{non-convex} example problems, see however \cite{kawohl2000some,kovtunenko2022shape} for a collection of historical problems and a theoretical investigation of shape differentiability in a non convex setting, respectively. 

\kk{To fill this gap,} we modify a mechanical shape optimization problem where one strives to maximize the reliability of a simple mechanical component while keeping the material consumption \kk{bounded.} 
This problem has been extensively studied by ourselves and co-authors and has exposed astonishingly stable convergence properties.

To render this problem manifestly non-convex, we introduce obstacles which, when penetrated by the shape, result in an additional penalty proportional to the penetration area.  This can be seen as a relaxed version of a common problem, where the installation space is \kk{partially} occupied by several components which mutually have to avoid each other. For some recent studies on shape optimization under installation space constraints, that however do not specifically refer to non convexity, we refer to \cite{kodiyalam2001multidisciplinary,muller2022scalable,werner2021package}. If now the component's \kk{initial guess, i.e., the} starting shape \kk{for the optimization process,} is on the wrong side of \kk{such a} barrier and the gradient of the penetration cost does not outweigh the gradient of the \kk{original} objective function, a gradient flow is trapped on the wrong side of the obstacle while the Hamiltonian flow might overcome it with the aid of momentum and thereby reach the better configuration on the other side. 

That this actually happens for adequate settings of our (discretized) dissipative Hamiltonian flow is shown in this work. \hg{We also observe that the shapes obtained after penetration of the obstacle are even superior to shapes we obtained in a previous work by gradient descent methods \cite{bolten2021tracing,Doganay2019}. Taking the solutions obtained by the dissipative Hamiltonian flow as an initial point for a bi-criteria tracing of a local Pareto front, we find that this improvement is consistently achieved over large parts of local Pareto fronts, which again emphasizes the importance of non-local optimization methods in (multi-objective) shape optimization.}

Our paper is organized as follows:
in Section~\ref{sec:dHS} we introduce the Hamiltonian approach to optimization and discuss the main properties of this method. We also recall the convergence of Hamiltonian flows to critical points which further motivates our approach. Our non-convex shape optimization problem based on avoidance \kk{of certain areas} in the installation space is introduced in Section~\ref{sec:NCSO}.  Numerical experiments are documented in Section~\ref{sec:Numerics}\hg{, where we also use the improved starting points from the dissipative Hamiltonian flow for tracing a local Pareto front which is consistently improving previous results based on gradient descent.} We give our conclusions and recommendations for future work in the final Section~\ref{sec:CO}

\section{\label{sec:dHS}Dissipative Hamiltonian flows and optimization}

We consider an objective function $f:\mathbb{R}^n\to \mathbb{R}$, $n\in \mathbb{N}$ \kk{and an unconstrained minimization problem $\min_{q\in\R^n} f(q)$.} Assuming $\R^n\ni q\mapsto f(q)$,  to be a lower bounded differentiable map with compact level sets and a locally Lipschitz first derivative, it is then easy to see that the gradient flow
\begin{equation}
    \label{eq:gradientFlow}
    \dot q(t)=-\nabla f(q(t)),~~q(0)=q_0,
\end{equation}
has a global solution for $t\in [0,\infty)$. In non-convex optimization, the first goal is to find critical points which fulfill the first order optimality conditions \cite{wright1999numerical}. From the perspective of gradient flows, this is equivalent to find stationary points \kk{$q_c$} of the dynamical system \eqref{eq:gradientFlow} fulfilling $\nabla f(q_c)=0$.  It is well-known that if $f$ also is a Morse function, i.e.\ is second order differentiable with isolated  critical points, then $\lim_{t\to\infty} q(t)=q_c$ for some critical point $q_c\in\R^d$ holds for all starting points $\kk{q_0}\in \R^n$, see \cite[Lemma 8.4.7]{jost2008riemannian} for a slightly stronger result. 

Discretizing \eqref{eq:gradientFlow} with the Euler scheme with stepsize $\alpha$ then leads to
\begin{equation}
    \label{eq:EulerScheme}
\frac{q(t+\alpha)-q(t)}{\alpha}\approx \dot q(t)= -\nabla f(q(t))   ~~\Leftrightarrow~~q(t+\alpha)\approx q(t)-\alpha \nabla f(q(t)),
\end{equation}
    where for $t=\alpha k$, $k\in\mathbb{N}$, on the right hand side we get the update rule for the iterate $q_{k+1}=q_k- \alpha \nabla f(q_k) \approx q(\alpha (k+1))$ for the gradient descent algorithm with stepsize $\alpha>0$ and initial parameter $q_0$. 
    As the Euler scheme converges for $\alpha\to 0$, we see that the iterates of the gradient descent algorithm in this limit follow the gradient flow.   

    The nice asymptotic convergence properties of the gradient flow therefore also shed light on the convergence of gradient descent methods, which are of course well understood, see e.g.\  \cite{wright1999numerical}. This advantage however leads to the disadvantage that gradient descent algorithms with small step size -- like the gradient flow -- tend to get stuck in the first local minimum it encounters.

    Also, in complex energy landscapes as e.g.\ encountered in machine learning, it is common to re-define the update scheme via a time series soothing approach \cite{metcalfe2009introductory}  with $\bar{\alpha}>0$, and $0\leq \beta<1$
     \begin{equation}
       \label{eq:discreteHamiltonian}  
\begin{array}{c}q_{k+1}=q_k-\bar{\alpha}\nabla_q f(q_{k})+\beta (q_k-q_{k-1})\\
\Leftrightarrow \\\left\{ \begin{array}{rl}
   p_{k+1}&=\displaystyle p_k -\alpha \nabla_q f(q_k)-\alpha \frac{\gamma}{m} p_k\\
    &\\
      q_{k+1}&= \displaystyle q_k+\frac{\alpha}{m} p_{k+1}
\end{array}    
\right.,
\end{array}
     \end{equation}
    where we used $p_k=\frac{m}{\alpha}(q_k-q_{k-1})$, $m=\frac{\alpha^2}{\bar{\alpha}}$ and $\gamma=\frac{\alpha(1-\beta)}{\bar{\alpha}}$ and initialized at $q_0,q_{-1}\in\mathbb{R}^n$ or $q_0, p_0=\frac{m}{\alpha}(q_0-q_{-1})$, respectively. $p_k$ is referred to as momentum and $m>0$ is called the mass and $\gamma>0$ the coefficient of Newtonian friction. $\alpha>0$ is a parameter that can be freely chosen, setting a 'time' scale. It has been observed in a number of papers that \eqref{eq:discreteHamiltonian} can be interpreted as a first order discretization of the dissipative Hamiltonian system 
    \begin{equation}
        \label{eq:Hamiltonian_system_seperate}
        \left\{ \begin{array}{rl}
   \displaystyle \dot p(t)&=\displaystyle - \nabla_q f(q(t))- \frac{\gamma}{m} p(t)\\
    &\\
      \displaystyle \dot q(t)&= \frac{1}{m} p(t)
\end{array}\right.,~~q(0)=q_0,~~p(0)=p_0.
    \end{equation}
    As usual in Hamiltonian dynamics, this system can now be brought in a compact, energy based form
    \begin{equation}
    \label{eq:Hamiltonan_system}
    \dot x(t)=
    \left(\begin{array}{c}
         \dot q(t)  \\
         \dot p(t)
    \end{array}\right)
    = (J-R)\nabla_x\mathcal{H}(x), ~~x(0)=x_0=\left(\begin{array}{c} q_0\\ p_0 \end{array}\right),
    \end{equation}
    where we used $x(t)=\left(\begin{array}{c} q(t)\\ p(t) \end{array}\right)$, $J=\left(\begin{array}{cc} 0&\mathbbm{1}\\ -\mathbbm{1} &0\end{array}\right)$ and $R=\left(\begin{array}{cc} \kappa&0\\ 0 &\frac{\gamma}{m}\mathbbm{1}\end{array}\right)$, where $\mathbbm{1}$ is the $n\times n$ unit matrix. Here we introduced an additional parameter $\kappa\geq 0$ which for $\kappa=0$ reproduces \eqref{eq:Hamiltonian_system_seperate}. If $\kappa>0$, the matrix $R$ becomes strictly positive definite, which makes the analysis of dissipativity more simple. Note that $\kappa>0$ creates a term $-\kappa \nabla_qf(q(t))$ on the right hand side of the lower equation in \eqref{eq:Hamiltonian_system_seperate}, which combines dissipative Hamiltonian mechanics with the gradient flow \eqref{eq:gradientFlow}, see \cite{massaroli2019port,poli2020port}  for applications in machine learning. 
    
    The Hamiltonian function $\mathcal{H}(x)$ is defined as energy via
    \begin{equation}
        \label{eq:Hamiltonian function}
        \mathcal{H}(x)=\mathcal{H}\left(\left(\begin{array}{c} q\\ p \end{array}\right)\right)=E_\text{kin.}(p)+E_\text{pot.}(q)=\frac{\|p\|^2}{2m}+f(q).
    \end{equation}
    Hence, our objective function $f(q)$ plays the role of potential energy $E_\text{pot}$, whereas the term $\frac{\|p\|^2}{2m}$ is the kinetic energy $E_\text{kin.}$ with $\|p\|^2=p^\top p$ the squared Euclidean norm on $\mathbb{R}^n$.    

    We now compile some well-known facts about dissipative Hamiltonian systems. Again, we assume that $f(q)$ has the properties given above. By $\bar f:=\sup_{q\in\R^n}f(q)\in \R\cup\{\infty\}$ we denote the supremum of $f(q)$. First, by the Picard-Lindelöf theorem, $x(t)$ has local solutions for $t\in [0,T]$, for some $T\in(0,\infty)$. Second, for any such solution, we obtain the dissipativity inequality for $0\leq s<t\leq T$, i.e.,
    \begin{align}
    \label{eq:Dissipativity}
    \begin{split}
    \mathcal{H}(x(t))-\mathcal{H}(x(s))&=\int_s^t \left\langle \nabla H(x(\tau)),(J-R)\nabla H(x(\tau))\right\rangle\, \mathrm{d}\tau \\
    &=-\int_s^t\left\langle \nabla H(x(\tau)),R\nabla H(x(\tau))\right\rangle\, \mathrm{d}\tau\leq 0
    \end{split}
    \end{align}
    holds as $R$ is positive semi-definite and $J^\top=-J$ is skew symmetric and hence $\langle v,Jv\rangle=0$ for all $v\in\mathbb{R}^n$. Here we used the chain rule $\dot {\mathcal{H}}(x(\tau))=\langle \nabla_x\mathcal{H}(x(\tau)),\dot x(\tau)\rangle$ along with \eqref{eq:Hamiltonan_system}.  Note that this inequality remains valid for abitrary Lipschitz differentiable Hamiltonian functions $\mathcal{H}:\mathbb{R}^{2n}\to\mathbb{R}$, skew symmetric $J\in \text{Mat}_{2n\times 2n}(\mathbb{R})$ and positive semidefinite $R\in \text{Mat}_{2n\times 2n}(\mathbb{R})$.

    From the dissipation inequalities one now easily obtains the existence of global solutions \kk{if} the starting point $x_0=\left(\begin{array}{c} q_0\\ p_0 \end{array}\right)$ \kk{satisfies} $\frac{\|p_0\|^2}{2m}< \bar f-f(q_0)$. \kk{Recall that $\mathcal{H}(x_0)=f(q_0)+\frac{\|p_0\|^2}{2m}$ and} let $\varphi\in (\mathcal{H}(x_0),\bar f)$. 
    By assumption, \kk{the trajectory} $q(t)$, the $q$-coordinate of the \kk{trajectory of } $x(t)$, starts in \kk{the level set} $\{q\in\mathbb{R}^n:f(q)\leq\varphi\}$ and the boundary of this set, the level \kk{curve} $\{q\in\mathbb{R}^n:f(q)=\varphi\}$  can never be crossed due to \eqref{eq:Dissipativity}. Thus the trajectory $x(t)$ can not produce runaway solutions and never leaves a compact set, on which the Lipschitz constant of $\nabla_qf(q)$ can be chosen uniformly. Under these conditions, the solution is known to exist for all times \cite{MR2439721}.

    Like in the case of gradient flows, it thus makes sense to ask for the asymptotic behavior of $x(t)$ and study its relation to the stationary points of the dynamical system \eqref{eq:Hamiltonan_system}. 

    We first analyze the matrix $J-R$ and show that it is invertible. In fact $(J-R)\left(\begin{array}{c} q\\ p \end{array}\right)=\left(\begin{array}{c} p-\kappa q\\ -q-\gamma p \end{array}\right)=0$ is only possible if $p=0$ and $q=0$ as $\kappa\geq 0$ and $\gamma,m>0$, from which $\mathrm{ker}(J-R)=\{0\}$ follows. Thus, stationary points $x_s=\left(\begin{array}{c} q_s\\ p_s \end{array}\right)$ of \eqref{eq:Hamiltonan_system}, where $(J-R)\nabla_x\mathcal{H}(x_s)=0$ fulfill $\nabla_x \mathcal{H}(x_s)=0$. 

    Let us analyze this condition further. First, since $\nabla_x=\left(\begin{array}{c} \nabla_q\\ \nabla_p \end{array}\right)$ and $\nabla_x\mathcal{H}\left(x\right)=\left(\begin{array}{c} \nabla_q f(q)\\ \frac{p}{m} \end{array}\right)=0$, $\nabla_x \mathcal{H}(x_s)=0$ is equivalent to $\nabla_q f(q_s)=0$ and $p_s=0$. Thus, the coordinates $q_s$ of a stationary point $x_s$ are critical points of the original optimization problem \kk{$\min_{q\in\R^n} f(q)$}, i.e.\ they satisfy the first order optimality condition $\nabla_q f(q_s)=0$. 

    In the next step, one  uses the dissipativity inequality \eqref{eq:Dissipativity} to prove under the given hypotheses that
        \begin{equation}
        \label{eq:vanishingGrad}
        \lim_{t\to\infty}\|\nabla_x\mathcal{H}(x(t))\|=0 
        \end{equation}
   holds. More precisely, one assumes that for some $\varepsilon>0$, $\|\nabla_x\mathcal{H}(x(t_n))\|>\varepsilon$ holds for some sequence $t_n\to \infty$ and derives a contradiction to \eqref{eq:Dissipativity}, see e.g.\ \cite{atto:theh:2000} for a detailed proof for the (more involved) situation with pure Newtonian friction, i.e.\ with $\kappa=0$.

If $f$ is a Morse function, i.e.\ if all critical points of $f$ are isolated, this immediately implies the convergence of $q(t)$ to a critical point $q_s$ with $\nabla_qf(q_s)=0$. The details of the proof, which once more is obtained by contradiction, again can be found in \cite{atto:theh:2000}.

We have thus seen that the momentum method, in the continuum time formulation, leads to guaranteed convergence to critical points \kk{under appropriate assumptions.} These results in continuum time are of particular interest if the time steps are small. In the context of shape optimization with PDE constraints this generally holds true, as the avoidance of strong mesh distortion in the morphing steps of shapes interdict large step sizes in order to guarantee a numerically clean solution to the state equation.

\section{\label{sec:NCSO}Optimizing reliability  under spatial and cost constraints} 

In the following, we extend a biobjective shape optimization problem for ceramic components under tensile load investigated in \cite{bolten2021tracing,Doganay2019}, \kk{that considers reliability and volume (cost) as optimization criteria,} by incorporating an additional objective functional penalizing the penetration area of a shape with a given obstacle. 
Hence, we consider a multiobjective shape optimization problem \kk{where we combine three objective functionals into a weighted sum.} 
\kk{While} the results of the numerical experiments in \cite{bolten2021tracing,Doganay2019} \kk{suggest that the corresponding Pareto fronts are (at least locally) convex,} 
we \kk{show in this paper that} introducing a \kk{circular} obstacle $\varpi=\varpi(x_{mp},r)$, where $x_{\text{mp}}$ is the midpoint and $r$ the radius of the circle, \kk{leads to non-convex optimization problems in general.} Note that the circle is not a \kk{forbidden} area for the \kk{shape,} but an \kk{intersection with it is} expensive due to the penalization of the penetration area. We consider the reliability of the ceramic component, its volume, and its penalized penetration area with a given circle as individual objectives \kk{that are combined into one weighted sum objective.} 
Following \cite{bolten2021tracing,Doganay2019}, the reliability of the component is evaluated as its probability of failure as discussed in \cite{bolten2015} and implemented in \cite{Hahn2019} for 2D shapes. We only give a brief summary of the model and refer to \cite{Doganay2019} for further details.

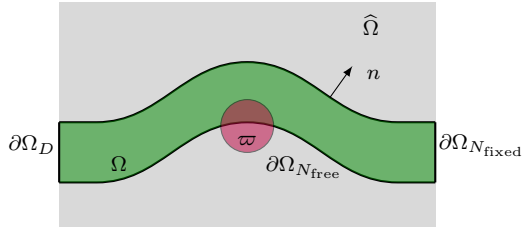
\begin{figure}[htb]
  \centering
  \begin{tikzpicture}[font={\footnotesize}]
    \coordinate (bb_ll) at (0,0.4);
    \coordinate (bb_lr) at (5,0.4);
    \coordinate (bb_tl) at (0,3.4);
    \coordinate (bb_tr) at (5,3.4);
    \coordinate (bar_ll) at (0,1.);
    \coordinate (bar_lr) at (5,1.);
    \coordinate (bar_tl) at (0,1.8);
    \coordinate (bar_tr) at (5,1.8);
    \fill [fill=black!15!white] (bb_ll) rectangle (bb_tr);
    \draw [thick,name path=A] (bar_ll) to[out=0,in=180] (0.5,1) to[out=0,in=180] (2.5,1.8) to[out=0,in=180] (4.5,1) to[out=0,in=180] (bar_lr);
    \draw [thick,name path=B] (bar_tl) to[out=0,in=180] (0.5,1.8) to[out=0,in=180] (2.5,2.6) to[out=0,in=180] (4.5,1.8) to[out=0,in=180] (bar_tr);
    \tikzfillbetween[of=A and B] {black!45!green, opacity=0.5};
    \node[label=left:\(\widehat{\Omega}\)] at (4.5,3.1) {};
    \node[label=left:\(\Omega\)] at (1.15,1.25) {};
    \node[label=right:\(\partial\Omega_{N_{\text{fixed}}}\)] at (4.8,1.5) {};
    \node[label=right:\(\partial\Omega_{N_{\text{free}}}\)] at (2.5,1.2) {};
    \node[label=left:\(\partial\Omega_{D}\)] at (0.2,1.5) {};
    \draw [thick] (bar_ll) -- (bar_tl) ;
    \draw [thick] (bar_lr) -- (bar_tr) ;
    \draw [-latex] (3.6,2.13) -- (3.9,2.55) node [label={[xshift=8pt,yshift=-12pt]\({n}\)}] {} ;    
    \filldraw [fill=purple, opacity=0.5](2.5,1.75) circle (10pt);    
    \node[label=\(\varpi\)] at (2.5,1.25) {};
  \end{tikzpicture}
  \caption{Illustration of an exemplary admissible shape \(\Omega\in \Oad\) and a circle $\varpi$ as an obstacle, compare with Figure~3 in \cite{bolten2021tracing}.\label{fig:adm_shape}}
\end{figure}

We consider a compact body $\Omega\subset \R^2$ that is filled with ceramic material and that has a piecewise Lipschitz boundary. Moreover, we assume that the boundary $\partial\Omega$ of $\Omega$ is \kk{subdivided} into three parts
\begin{align*}
	\partial\Omega = \text{cl}({\partial\Omega}_D) \cup \text{cl}({\partial\Omega}_{N_{\text{fixed}}}) \cup \text{cl}({\partial\Omega}_{N_{\text{free}}}),
\end{align*}
where $\text{cl}(\cdot)$ denotes the closure, \kk{that} the Dirichlet boundary condition holds on $\partial\Omega_D$, \kk{that} $\partial\Omega_{N_{\text{fixed}}}$ describes the part on which surface forces may act on, and \kk{that} $\partial\Omega_{N_{\text{free}}}$ is free to be modified during the optimization. Furthermore, we assume that all feasible shapes are contained in a bounded open set $\widehat{\Omega} \subset \R^2$ satisfying the \emph{cone property}, see, e.g., \cite{bolten}. In Figure~\ref{fig:adm_shape}, an example also containing an obstacle in form of a circle $\varpi\subset \widehat{\Omega} $ is illustrated.

The set of admissible shapes can then be defined as
\[\Oad:=\{\Omega \subset \widehat{\Omega}:\; {\partial\Omega}_D \subset {\partial\Omega},\; {\partial\Omega}_{N_{\text{free}}} \subset \partial\Omega,\; \widehat{\Omega} \text{ and }\Omega\text{ satisfy the cone property}\}.\]
Following \cite{braess,munz}, ceramics behave according to linear elasticity theory and therefore the state equation describing the behavior of the ceramic component under tensile load is given as the following partial differential equation:
\begin{equation} \label{stateequation}
\begin{array}{rcll}    
  -\text{div}(\sigma(u(z))) & = & \bar{f}(z) &  \text{for} \; z\in\Omega \\    
  u(z) & = & 0 &  \text{for} \; z\in\partial\Omega_D    \\
 \sigma(u(z)){n}(z) & = & \bar{g}(z) &  \text{for} \; z\in\partial\Omega_{N_{\text{fixed}}} \\
  \sigma(u(z)){n}(z) & = & 0 & \text{for} \; z\in\partial\Omega_{N_{\text{free}}} .
\end{array}
\end{equation}
Here, the volume forces are given by $\bar{f}\in L^2(\Omega , \R^2)$ and the forces acting on the surface $\partial \Omega_{N_{\text{fixed}}}$, e.g.\ the tensile load, \kk{are given} by $\bar{g} \in L^2(\partial\Omega_{N_{\text{fixed}}} , \R^2)$. The outward pointing normal \kk{is assumed to be defined} nearly everywhere on $\partial \Omega$ and is denoted by ${n}(z)$ at $z\in\partial\Omega$. The displacement of the component is represented by $u\in H^1(\Omega, \R^2)$ and the Jacobian of $u$ by $\nabla u$. Hence, the linear strain tensor $\varepsilon\in L^2(\Omega, \R^{2\times 2})$ is given by $\varepsilon(u(z)):=\frac{1}{2}(\nabla u(z) + (\nabla u(z))^{\top})$. Furthermore, for the stress tensor $\sigma\in L^2(\Omega, \R^{2\times 2})$ we have that $\sigma(u(z))=\hat{\lambda}\, \text{tr}(\varepsilon(u(z)))\,I+ 2\,\hat{\mu}\,\varepsilon(u(z))$, where the Lam\'e constants $\hat{\lambda}=\frac{\nu E}{(1+\nu)(1-2\nu)}$ and $\hat{\mu}=\frac{E}{2(1+\nu)}$ are derived from Young's modulus $E$ and Poisson's ratio $\nu$.

The 
reliability of the component is then modelled by a Poisson point process following \cite{bolten2021tracing,bolten2015,Doganay2019}. More precisely, \kk{we use an} intensity measure  that counts the \kk{potential} cracks in the component which may initiate ruptures under tensile load. We then obtain the following Weibull-type  objective functional \kk{representing the probability of failure of the shape:}
\begin{align*}
J_1(\Omega):=J_1(\Omega, \nabla u):=\frac{1}{2\pi}\int\limits_{\Omega}\int\limits_{S^1}\left(\frac{\Bigl(n^\top \sigma(\nabla u(z))\,n \Bigr)^+}{\sigma_0}\right)^{\bar m} \text{d}n\,\text{d}z .
\end{align*}
Here, $S^1$ denotes the unit sphere in $\R^2$ and $\bar m$ the \emph{Weibull module} which typically assumes values between $5$ and $25$. Furthermore, $\sigma_0$ is a positive constant and $(\cdot)^+:=\max(\cdot, 0)$. For further details we refer to \cite{bolten2015}.
The second objective functional corresponds to the volume of the shape \kk{(representing its material consumption and hence its cost) and} is given by $J_2(\Omega):=\int_{\Omega}\text{d}z$. For a given circle $\varpi=\varpi(x_{\text{mp}},r)$ \kk{that represents an area that should be avoided by the shape,} the penalizing objective functional has the form $J_3(\Omega):=J_3(\Omega, \varpi):=c_\text{P}\int_{\Omega\cap\varpi}\text{d}z$, where $c_\text{P}>0$ is a penalization parameter.

Now we can formulate a non-convex multiobjective shape optimization problem as
\begin{equation}
\begin{split}
\min_{\Omega \in \Oad} & ~J(\Omega):=(J_1(\Omega),J_2(\Omega),J_3(\Omega))\\
\text{s.t. } & u \in H^1(\Omega, \R^2) \text{ solves the state equation } (\ref{stateequation}).
\end{split}\label{ceramicMOP}
\end{equation}
We are interested in finding Pareto optimal shapes $\Omega\in\Oad$ for which the improvement in one objective always leads to a deterioration in at least one other objective. More formally, a shape $\Omega\in\Oad$ is called \emph{Pareto optimal} when its image $J(\Omega)$ is \emph{non-dominated}, i.e., when there is no other shape $\Omega'\in\Oad$ such that $J_i(\Omega')\leq J_i(\Omega)$ for all $i=1,2,3$ 
and $J(\Omega')\neq J(\Omega)$. We refer to \cite{ehrg:mult:2005,miet:nonl:1998} for further details on multiobjective optimization in general, and on scalarization techniques in particular. 
In the following, we utilize a \emph{weighted sum scalarization} of the three objectives $J_1,J_2$ and $J_3$:  For a weight vector $\lambda=(\lambda_1,\lambda_2, \lambda_3)\in\R^3$ with \kk{$\lambda_1,\lambda_2,\lambda_3>0$ and} $\lambda_1+\lambda_2+\lambda_3 =1$, \kk{the weighted sum scalarization} is given by
\[J_\lambda:=\lambda_1J_1+\lambda_2J_2+\lambda_3J_3: \Oad \rightarrow \R.        \]
A shape $\Omega\in \Oad$ is then called (locally) optimal with respect to $J_\lambda$ if $J_\lambda(\Omega)\leq J_\lambda(\Omega')$ for all $\Omega'\in \Oad$ \kk{($\Omega'$ in some neighborhood of $\Omega$), respectively.} 
Moreover, $\Omega$ is called critical (or Pareto critical) for $J_\lambda$ if $\nabla J_\lambda(\Omega)= \lambda_1\nabla J_1(\Omega)+ \lambda_2\nabla J_2(\Omega)+\lambda_3\nabla J_3(\Omega) =0$. \hg{Here $\nabla$ has to be understood as shape gradient \cite{schulz,sokolowski}, however if there exists a (surjective) finite dimensional parametrization of the admissible shapes $\mathbb{R}^n\supset \mathcal{U}\ni q\mapsto \Omega(q)\in\Oad$ for some open set $\mathcal{U}$, $\nabla$ can be replaced by the gradient $\nabla=\nabla_q$ and we write $J_\lambda(q)=J(\Omega(q))$.} 
Note that while an optimal solution of a weighted sum scalarization is always Pareto optimal, 
\kk{for non-convex problems} the weighted sum scalarization can not recover \kk{the complete Pareto set in general.} Since we do not follow a multiobjective approach in this work and are mainly interested in non-convex \kk{single-objective} shape optimization problems, we consider the weighted sum scalarization for \kk{a fixed} weighting vector $\lambda$ \kk{that represents some particular preferences.}

\hg{However, specific Pareto critical solutions found for a particular choice of parameters $\lambda$ can be used as starting points for further tracing the (local) Pareto front. In fact, one Pareto critical solution $q(\lambda)$ with respect to $J_\lambda$, satisfying $\nabla_q J_\lambda (q(\lambda))=0$, under adequate conditions \cite{bolten2021tracing} on the non degeneracy on the Hessian $\nabla_q^2 J_\lambda(q)$, by the implicit function theorem leads to the existence of a local manifold of Pareto critical solutions $\Omega(q(\lambda'))$ for $\lambda'$ in an open neighborhood of $\lambda$.}

\hg{In our experiments, we show that such local Pareto fronts can be efficiently traced on the basis of the ordinary differential equation that results from the formula of the implicit derivative in the implicit function theorem, see \cite{bolten2021tracing} for the details. In this way we obtain consistently improved local Pareto fronts from the improved specific solutions obtained by port Hamiltonian flows.}

\section{\label{sec:Numerics}Numerical experiments}

The momentum method described in Section~\ref{sec:dHS} is now tested and compared to a classical gradient descent approach on two particular instances of problem \eqref{ceramicMOP}. Towards  this end, the Hamiltonian flows (i.e., the trajectories of $q(t)$) are discretized and the occuring ODEs 
\otd{are solved using the symplectic Euler method, see, e.g.\ \cite{Hairer2006}.}

To evaluate the objective functions and gradients of $J_1$ and $J_2$ the implementation  of \cite{Hahn2019} is used. There, triangular Lagrangian finite elements are used to discretize two-dimensional shapes $\Omega \in \Oad $ by an $n_x\times n_y$ finite element mesh \(Z:=(Z^{\Omega}_{ij})_{n_x\times n_y}\). All integrals are computed via numerical quadrature. For the computation of the gradient of the intensity measure $J_1$ the computationally efficient adjoint approach is applied. Following again \cite{bolten2021tracing,Doganay2019}, we utilize a geometry definition that effectively reduces the number of variables by taking advantage of the geometry of the considered shapes. In a first step, we fix all $x$-components of the grid points and represent the discretized shape $Z$ via its \emph{mean line} and \emph{thickness} values  $\varrho^{\text{ml}}\in\R^{ n_x}$ and $\varrho^{\text{th}}\in\R^{ n_x}_+$. In a second step, we use B-splines with a prespecified number of $n_B$ basis functions $\vartheta_j,\ j=1,\hdots, n_B$ to fit these meanline and thickness values (see, e.g., \cite{nurbs}) to achieve smoothed meanline and thickness values via
\begin{equation*}
\hat{\varrho}^{\text{ml}}(z):=\sum_{j=1}^{n_B} q^{\text{ml}}_j \, \vartheta_j(z) \quad \text{and}\quad \hat{\varrho}^{\text{th}}(z):=\sum_{j=1}^{n_B} q^{\text{th}}_j \, \vartheta_j(z), \qquad z\in \R.
\end{equation*}
We then consider the B-spline coefficients $q=(q^{\text{ml}}, q^{\text{th}}) \in \R^{ n_B}\times \R^{ n_B}_+$ as our optimization variables, replacing $J_i(\Omega)$ by $J_i(Z)\approx J_i(q)$, $i\in\{1,2,3,\lambda\}$. 
To evaluate the objective function $J_3$, the area of intersection of the triangular finite elements of the discretized shapes and a given circular obstacle $\varpi$ is computed using the \emph{R} package \emph{'sf'}. Furthermore, the gradient $\nabla J_3$ is approximated with the finite difference approach. Since $J_2(\kk{q}),J_3(\kk{q})\in C^\infty$ and as shown in \cite{Doganay2019} $J_1(\kk{q})\in C^{\bar m}$ we have $J_\lambda(\kk{q})\in C^{\bar m}$ for all $\lambda\in\R^3_+$. For our numerical experiments we set $\bar m=5$ as in \cite{bolten2021tracing,Doganay2019}. Note that the optimality conditions of Section~\ref{sec:NCSO} still hold for \kk{$f=J_{\lambda}$ and $q$ as defined here.} 

\subsubsection*{Test Cases}
We modify the two 2D test cases that were investigated in \cite{bolten2021tracing,Doganay2019} by introducing a circle $\varpi$ as an obstacle for the shapes. \kk{Other than that,} the same set of model parameters and boundaries are used, \kk{i.e., we consider ceramic shapes made from} beryllium oxide (BeO). Following \cite{munz,crcmaterials}, the material parameters of BeO \kk{are chosen as follows:} Young's modulus $\texttt{E}=320\,\text{GPa}$, Poisson's ratio $\nu=0.25$, and  ultimate tensile strength $\text{uts}=140\,\text{MPa}$. We set $\bar m=5$ for the Weibull module. For both test cases we fix the length at $1.0\,\text{m}$ and the height of the left and right boundaries at $0.2\,\text{m}$. Here, the Dirichlet boundary $\partial\Omega_{D}$ is located on the left, i.e., it is fixed and force free, while the Neumann boundary $\partial\Omega_{N_{\text{fixed}}}$ is located on the right, i.e., it \kk{the shape is} fixed \kk{at the left side} and surface forces $\bar{g}$ may act on it \kk{on the right.} The remaining parts of the boundary are force free and can be adapted in an optimization scheme, i.e., \kk{they correspond to} $\partial\Omega_{N_{\text{free}}}$. \kk{As in} \cite{bolten2021tracing,Doganay2019} we set the tensile load to $\bar{g} = 10^{7}\,\text{Pa}$ and the gravity forces to $\bar{f} = 0\,\text{Pa}$. We use a triangular $41\times 7$ mesh, i.e.,  $n_x=41$ and $n_y=7$, for the discretization of the shapes. B-splines with $n_B=5$ basis functions are used to fit the meanline and thickness values, tallying ten B-spline coefficients. Omitting the fixed coefficients corresponding to the fixed boundaries yields in total six optimization variables, i.e., $q\in\R^6$, see also \cite{bolten2021tracing,Doganay2019}.
For both test cases the circles $\varpi=\varpi(x_{\text{mp}},r)$ are placed in such a way that the local optimal solutions computed in \cite{bolten2021tracing,Doganay2019} are beneath them and do not intersect them. The starting solutions differ from the ones utilized in \cite{Doganay2019}, as they are constructed such that they lie above the circles without intersecting them. For this 
setup there exists a local minimum of $J_\lambda$ above the circle, and another one beneath the circle. We compare the solutions of the presented \kk{momentum method} with the results of \kk{simple} gradient descent approach with Armijo step lengths, starting from the same initial solutions. 
All numerical experiments \kk{are based on an implementation} in \emph{R} (version 3.6.3), where an implementation of \cite{Hahn2019} is used to evaluate the objective values and the (adjoint) gradients \kk{of $J_2$} on the finite element grid.

\otd{We solve the occurring ODEs with the symplectic Euler method, see, e.g.\ \cite{Hairer2006}, and compare with \eqref{eq:discreteHamiltonian} and \eqref{eq:Hamiltonian_system_seperate}. The update scheme reads than as }

\[\begin{array}{c}
\otd{\left\{ \begin{array}{rl}
   p_{k+1}&=\displaystyle p_k -\alpha \nabla_q f(q_k)-\alpha \frac{\gamma}{m^2} p_k\\
    &\\
      q_{k+1}&= \displaystyle q_k-\alpha \kappa\nabla_q f(q_k)+\frac{\alpha}{m} p_{k+1}
\end{array}  
\right. ,}
\end{array}
\]
\otd{where $\alpha>0$ denotes a step in time.} 
\hg{The symplectic Euler method is known for a much improved energy balance, as compared to the standard Euler update scheme.  }

For the visualization of the circles we used the \emph{R} package \emph{`plotrix'}. \kk{While this is a straight-forward approach to visualize our results, plotrix does not always show the exact sizes of the respective circles since it does not account for the scaling of the $y$-axis in the plots. This may lead to the impression that the final shapes have a non-empty intersection with the obstacle even though this is actually not the case.} 
The plots of the starting shapes are an exception: In these plots the depicted circles are in full correspondence with the actual circles.

\subsubsection*{Test Case 1: A Straight Joint} 
\label{subsubsec:StraightJoint}

\begin{figure}[h]
	\begin{center}
		\subfloat[Exemplary optimal solution of \cite{bolten2021tracing} with added circle $\varpi$ (purple) placed above it.\label{subfig:TC1_PT_circ}]{\includegraphics[width=0.45\textwidth]{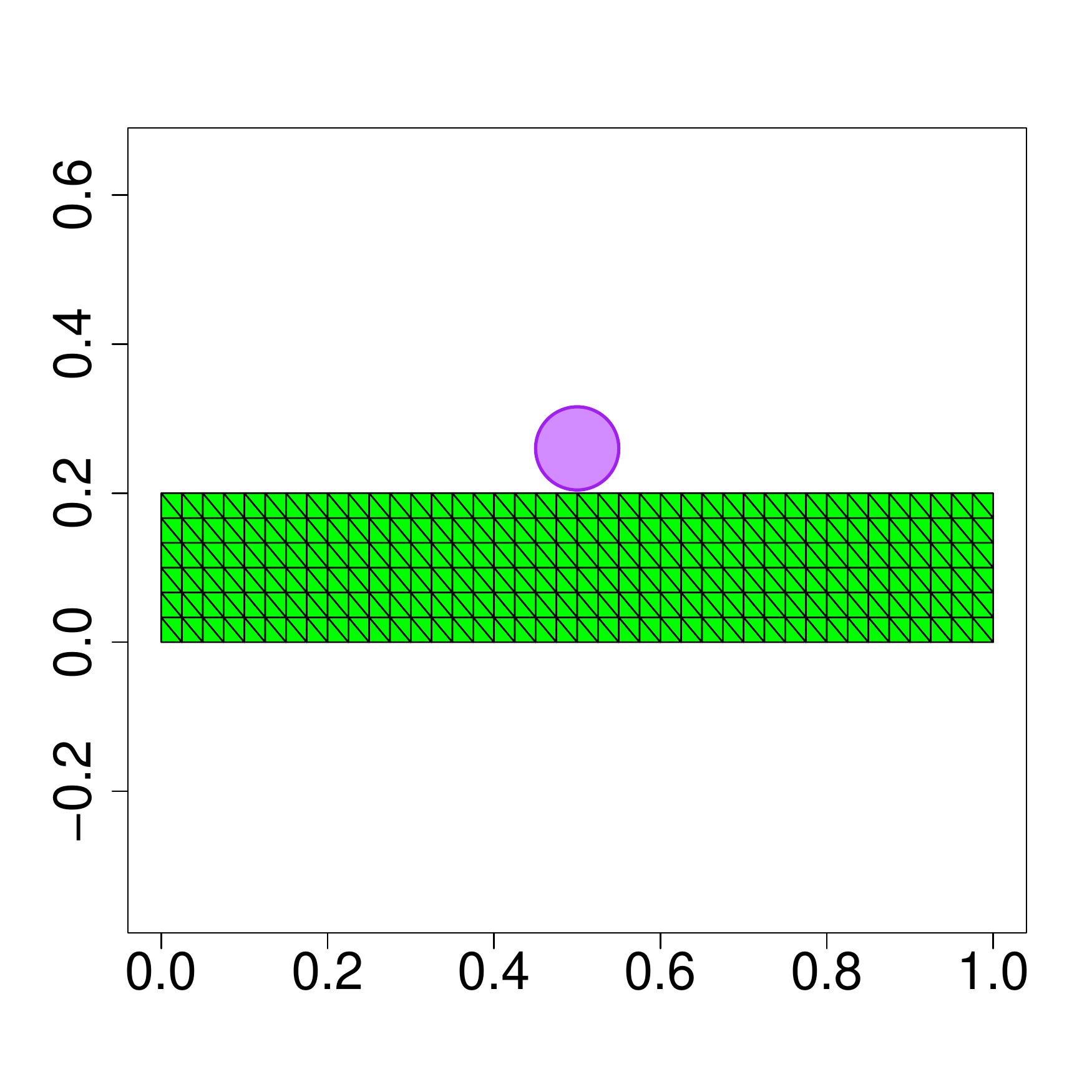}}
		\hspace{\fill}
		\subfloat[Initial shape $q^0$ 
            located above the obstacle $\varpi$. \label{subfig:TC1_initShape}
		]{\includegraphics[width=0.45\textwidth]{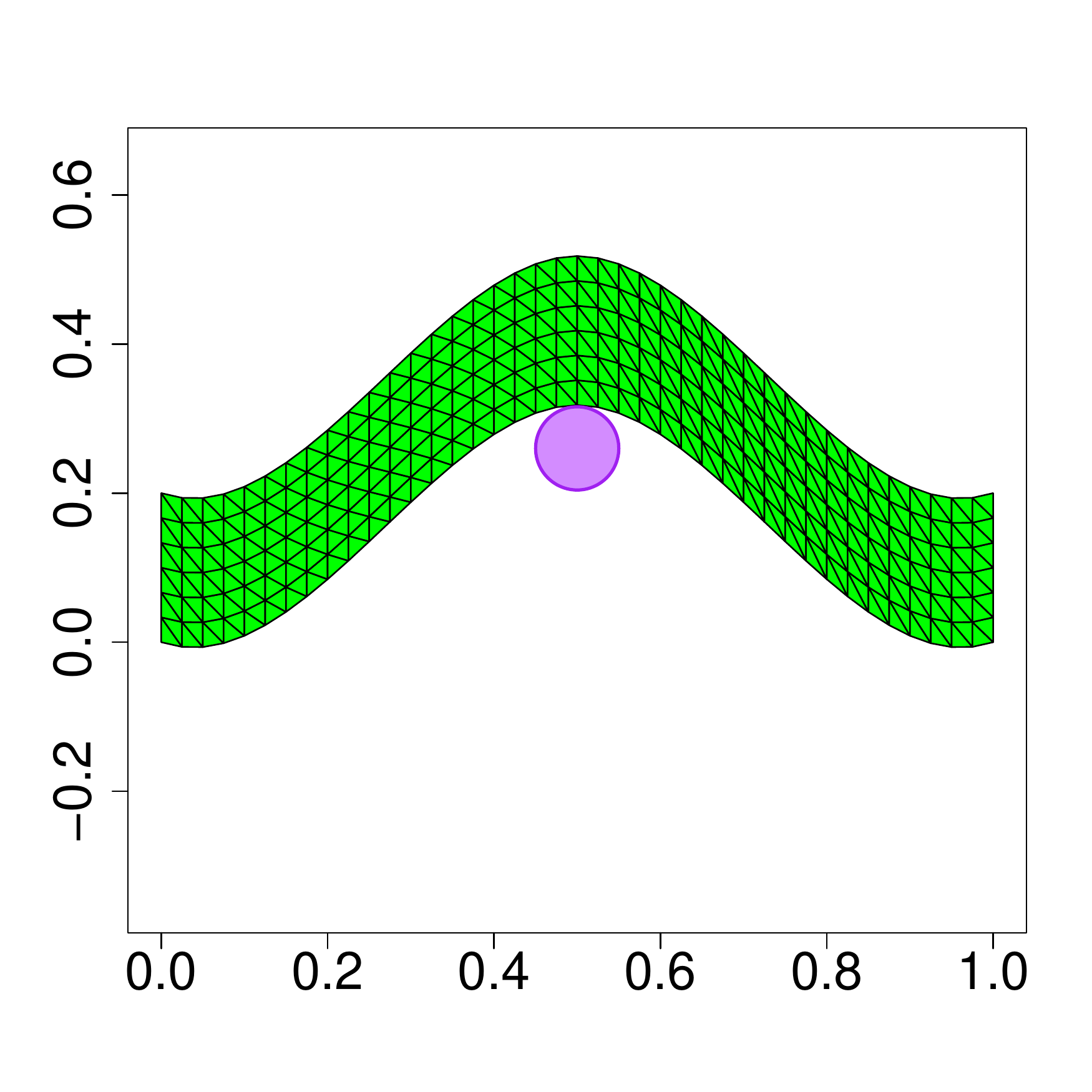}}
 	\end{center}			
	\caption{Test case 1 -- placement of circular obstacle $\varpi$ and initial shape.\label{fig:TC1_start}}
\end{figure}

As in \cite{bolten2021tracing,Doganay2019} we fix the left and right boundaries at the same height for the first test case. The unpenalized biobjective numerical studies of \cite{bolten2021tracing,Doganay2019} yielded straight rods with varying thickness as \kk{Pareto optimal solutions, which is consistent with our intuition.} Now we place the circular obstacle  $\varpi=\varpi(x_{\text{mp}},r)$ above the known solutions from the unconstrained problem. For this test case we choose a circle with midpoint $x_{\text{mp}}=(0.5, 0.26)$ and radius $r=0.05$, i.e., $\varpi=\varpi((0.5, 0.26), 0.05)$, see Figure~\ref{subfig:TC1_PT_circ} for an illustration. As stated before, we modify the starting solution for the numerical approaches in such a way that \kk{the initial shape} is located above the circle $\varpi$. The starting solution $q^0$ for this test case is illustrated in Figure~\ref{subfig:TC1_initShape}.

 For our numerical studies we choose the weight vector $\lambda$ as $\lambda=(0.4, 0.3, 0.3)$. \kk{Since the initial shape has a volume of $J_2=0.2$ and since the values of $J_1$ become very small for straight rods, this weight selection slightly favors the volume ($J_2$) over the reliability ($J_1$) at optimality. To ensure that the optimized shapes do not intersect the circular obstacle, we} set the penalization parameter $c_\text{P}$ of $J_3$ to $c_\text{P}=100$ \kk{and hence strongly penalizing non-empty intersections of the shape with the circular obstacle.} 
 
\begin{figure}[]
	\begin{center}
		\subfloat[Final solution $q^{\text{GD}}$ of the gradient descent method starting in $q^0$. Note that the actual area of intersection with $\varpi$ is smaller than the depicted one.\label{subfig:TC1_solGD}]{\includegraphics[width=0.45\textwidth]{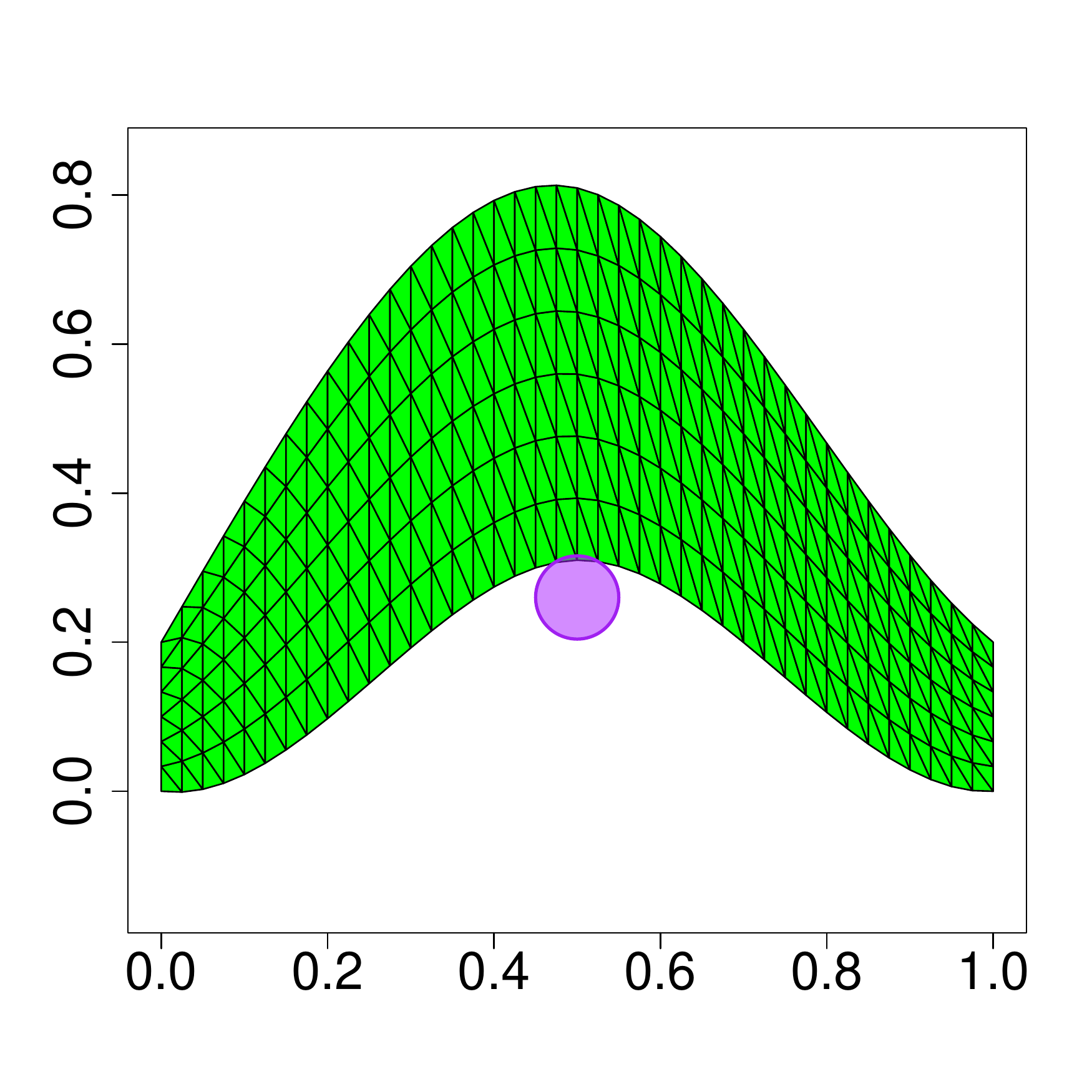}}
		\hspace{\fill}
		\subfloat[Final solution $q^{\text{HF}}$ of the \kk{momentum method after \otd{$250$} iterations} starting in $q^0$. 
        \label{subfig:TC1_soldHF}
		]{\includegraphics[width=0.45\textwidth]{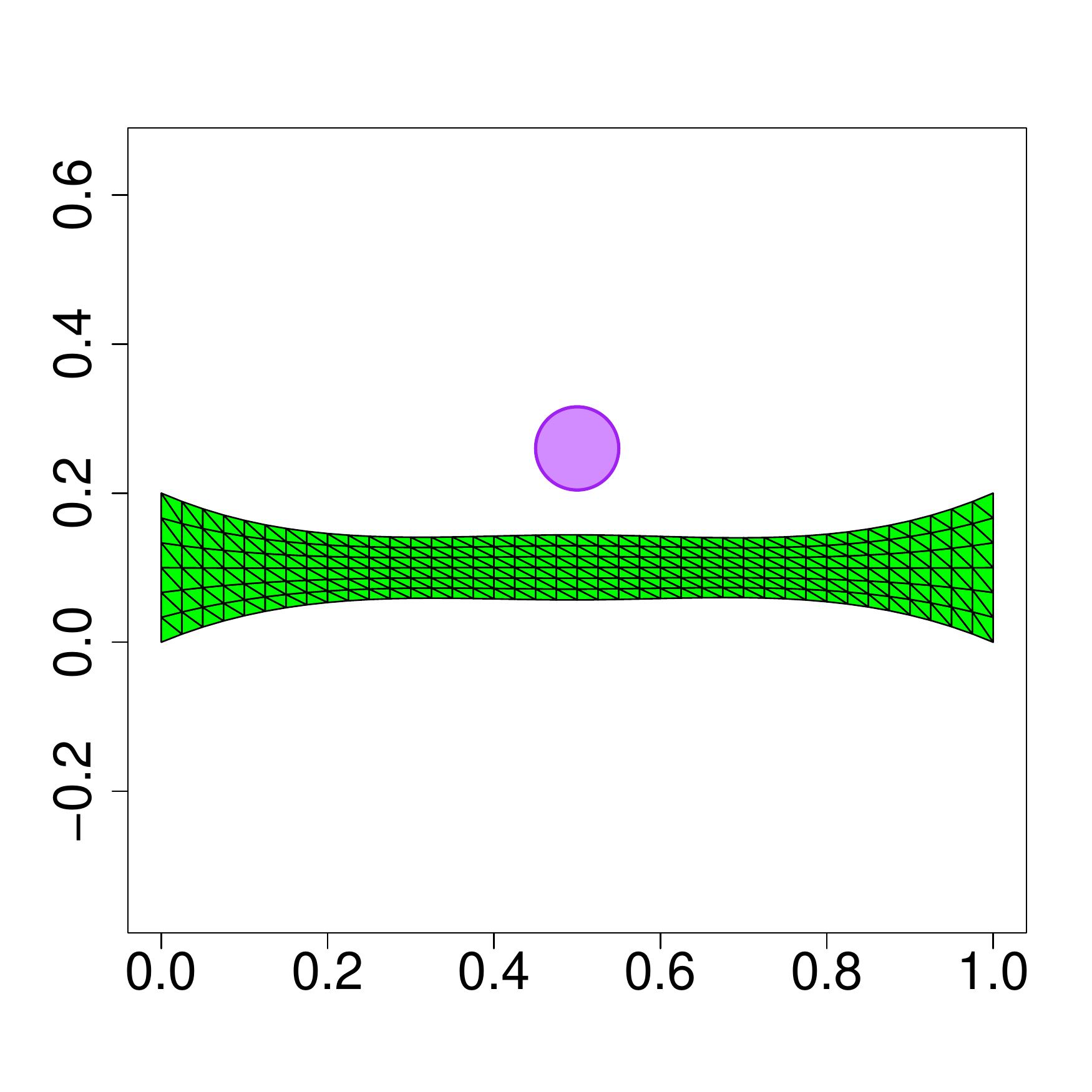}}
 	\end{center}		
	\caption{Test case 1 -- comparison of the results of gradient descent \kk{and momentum method}.\label{fig:TC1_comp}}
\end{figure}

 Starting in $q^0$ we first apply a gradient descent method with Armijo step lengths, a maximum iteration number of $200$ and the stopping condition $\|\nabla J_\lambda(q)\|\leq 10^{-5}$. The solution $q^{\text{GD}}$ \kk{obtained after} $152$ iterations \kk{satisfies this stopping condition and} has the objective value $J_\lambda(q^{\text{GD}})\approx 0.1584$. \kk{As can be seen in} Figure~\ref{subfig:TC1_solGD}, \kk{the gradient descent approach} did not \kk{pass over} the obstacle $\varpi$, i.e., \kk{its trajectory} got stuck in a \kk{local minimum} located above $\varpi$. 

 Starting from the same initial shape $q^0$, we next apply the momentum method and \kk{compute discretized Hamiltonian flows for} $\kappa=10^{-3}$, mass $m=10$ and dissipation parameter \otd{$\gamma=100$}. The initial momentum $p^0$ is set to a vector of zeros of the same dimension as $q^0$, i.e., $p^0=0\in\R^6$. Furthermore, the maximum time $T$ of the dynamics is set to \otd{$T=1$, with $250$ time steps, i.e., $\alpha=1/250$}. This approach yields a solution $q^{\text{HF}}$ with objective value \otd{$J_\lambda(q^{\text{HF}})\approx 0.0365$} that lies beneath the circular obstacle $\varpi$ and that resembles a straight rod corresponding to an established (local) minimum of \cite{bolten2021tracing}, see Figure~\ref{subfig:TC1_soldHF}.  
 
 The histories of the potential energy, i.e., $J_\lambda(q)$, the kinetic energy, \otd{and the total energy} during the \otd{$250$} iterations \kk{approximating} the dissipative Hamiltonian flow, starting in $q^0$ with $p_0=0$, \otd{are illustrated in Figure~\ref{fig:TC1_Ehist}.} 
  
\hg{In the first iteration, we observe a small increase of the total energy, as also the symplectic Euler scheme is not exact with respect to the energy balance. For a standard Euler scheme, this violation is much more pronounced, as we have observed in numerical tests.  }

The computed discretized Hamiltonian flow is in accordance with the modelling assumptions. Indeed, we observe that the potential energy of the initial shape is $3.9354$, which then drops to \otd{$0.3354$} after the \otd{early phase of the approach} and increases again to \otd{$0.3778$} in the following iterations, to subsequently decrease again while converging towards an optimal solution. \otd{Note that around the time $t=0.7$ there is a drop in the potential energy which corresponds to the shape overcoming the circular obstacle $\varpi$.} We observe that, as the potential energy decreases in \otd{the early phase of the approach}, the kinetic energy increases from $0$ to \otd{$2.6300$}. The kinetic energy then decreases while the potential energy increases, nicely capturing the interplay between these two energies.
Note that during the optimization process the potential energy may increase, however, without \otd{the total energy} surpassing the initial \otd{total} energy, which is in accordance with the fact that without an external energy supply the total energy in a dissipative Hamiltonian system can only decrease due to dissipation. \otd{In our numerical studies we experienced a small increase of the total energy at the start of the approach and when overcoming the circular obstacle $\varpi$ due to the precision of the symplectic Euler scheme.}
 
\begin{figure}[]
	\begin{center}
		\subfloat{\includegraphics[width=0.2\textwidth]{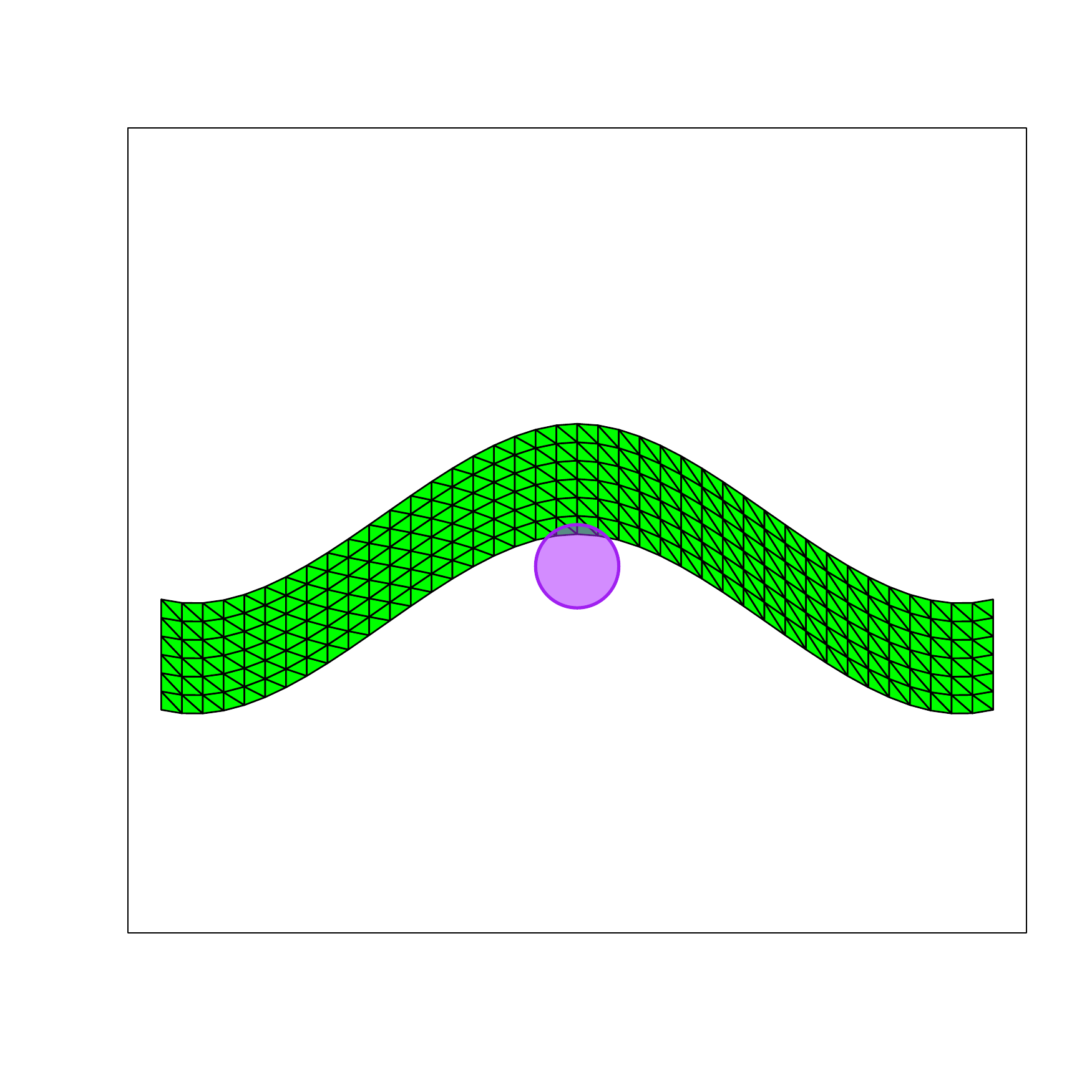}}
		\subfloat{\includegraphics[width=0.2\textwidth]{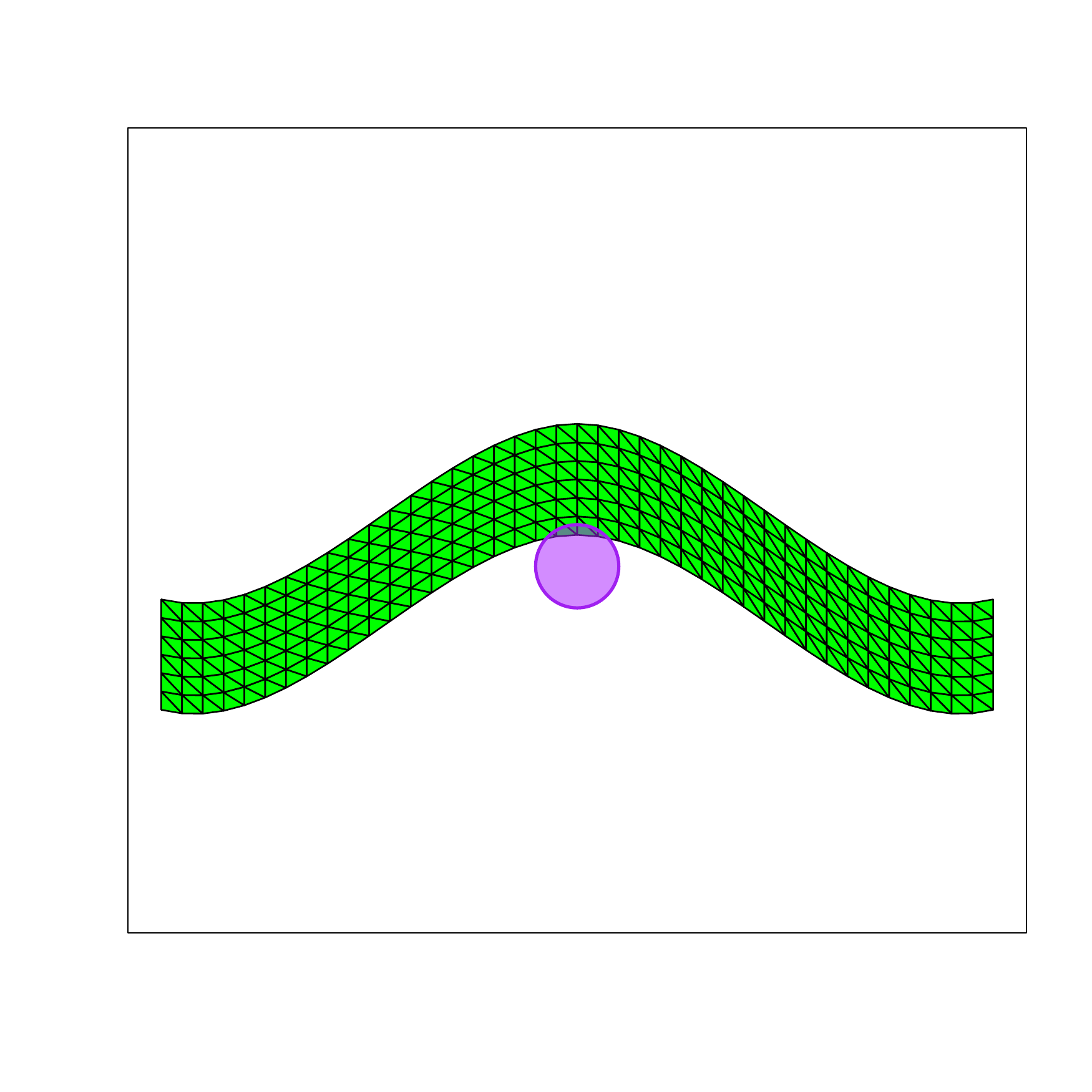}}
		\subfloat{\includegraphics[width=0.2\textwidth]{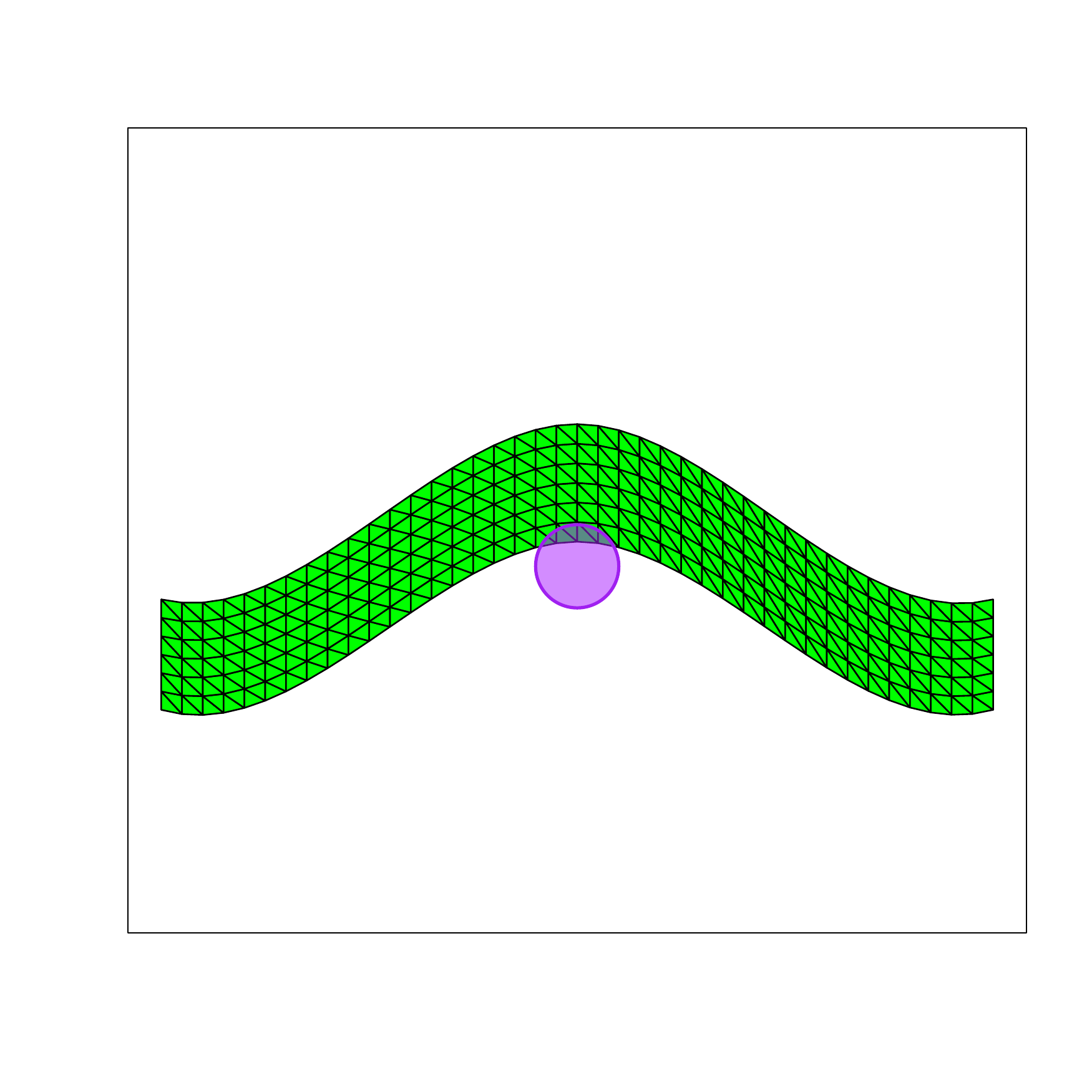}}
		\subfloat{\includegraphics[width=0.2\textwidth]{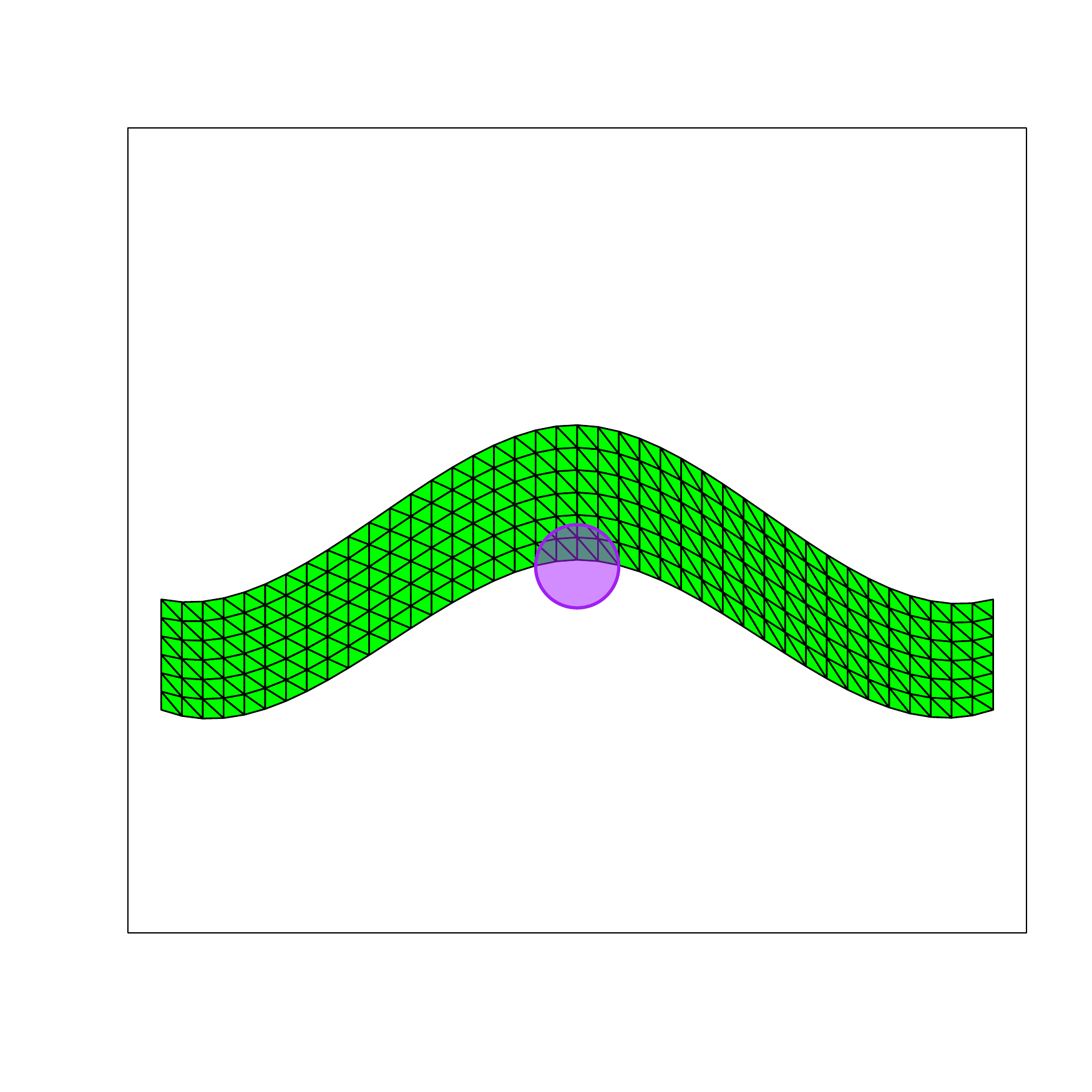}}
		\subfloat{\includegraphics[width=0.2\textwidth]{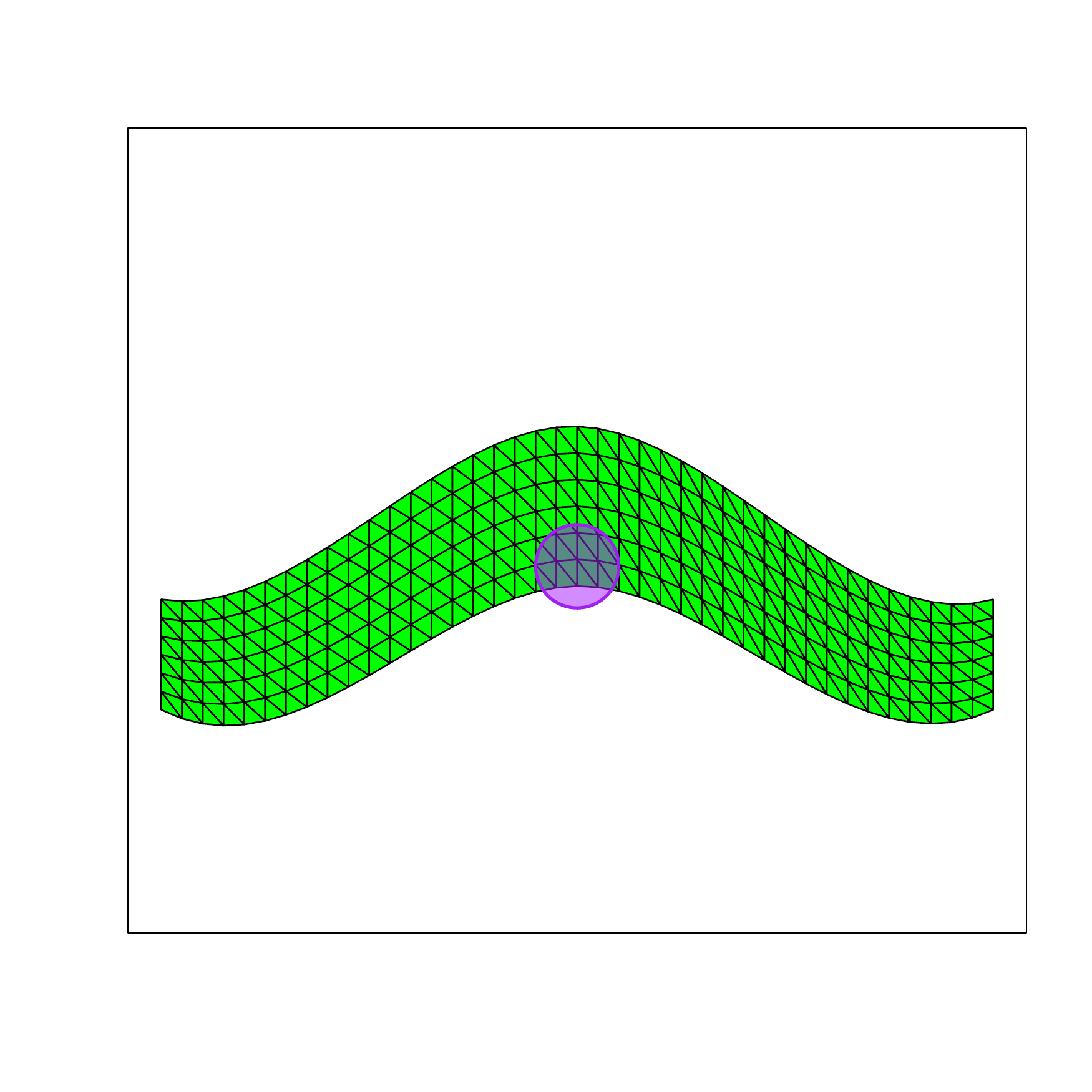}}

		\subfloat{\includegraphics[width=0.2\textwidth]{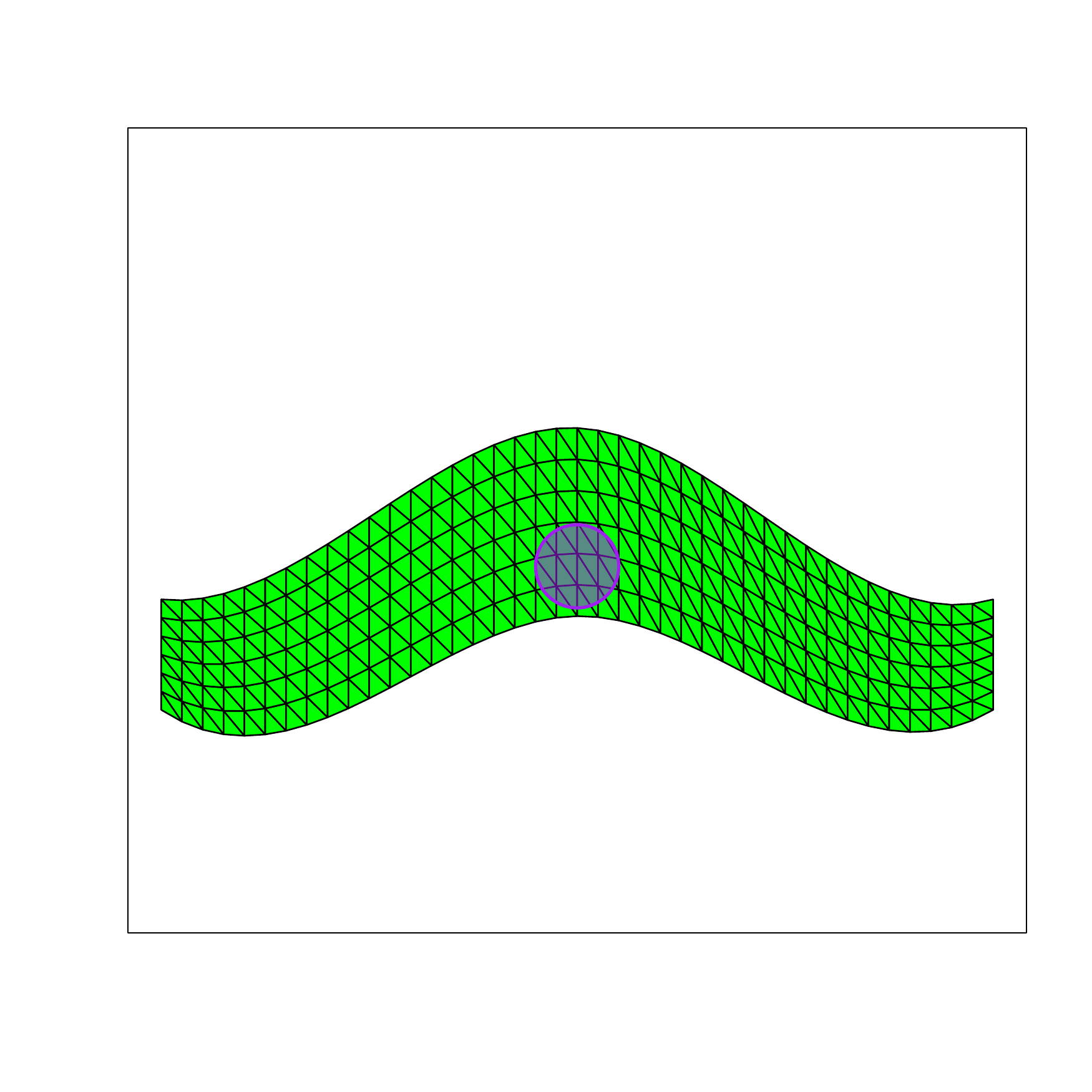}}
		\subfloat{\includegraphics[width=0.2\textwidth]{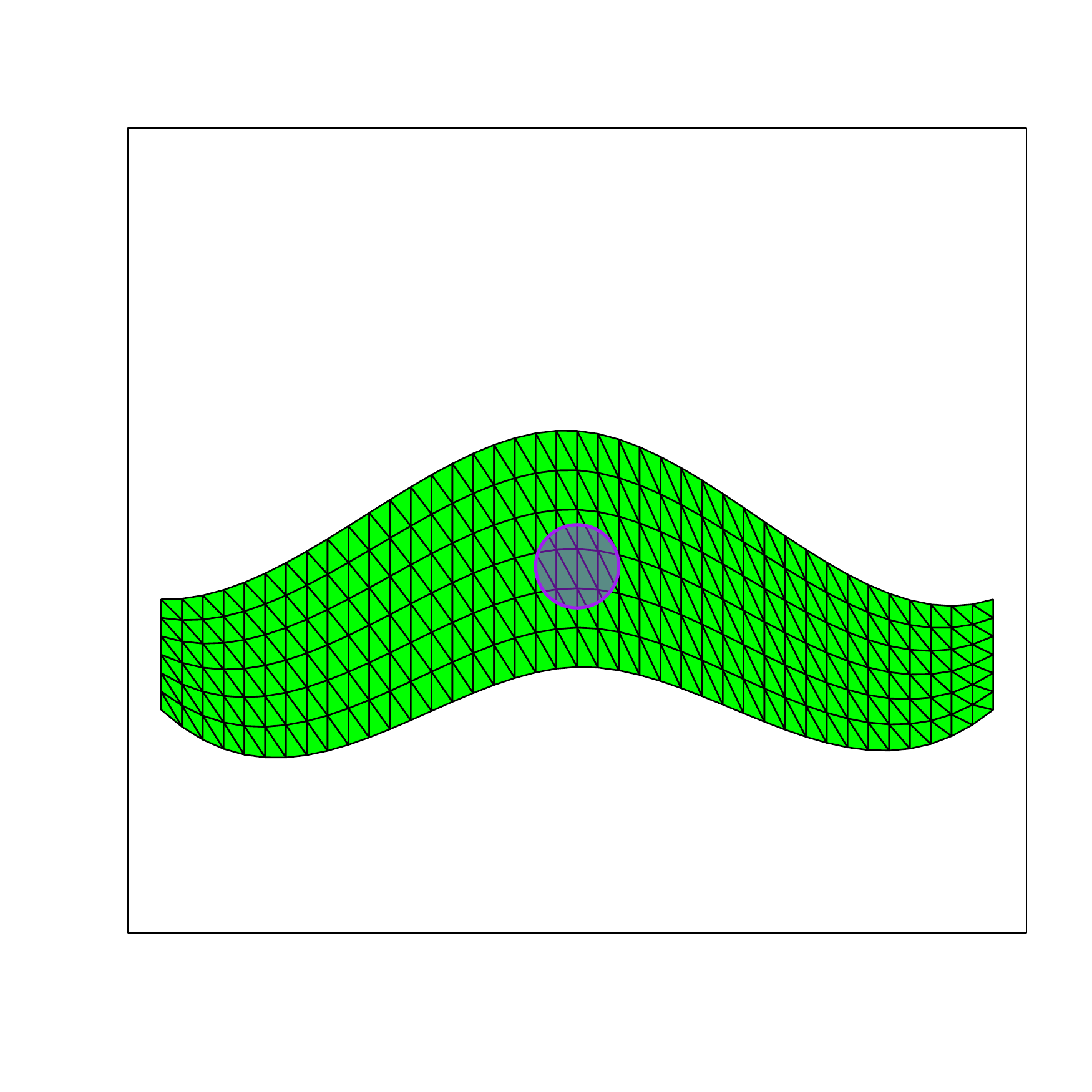}}
		\subfloat{\includegraphics[width=0.2\textwidth]{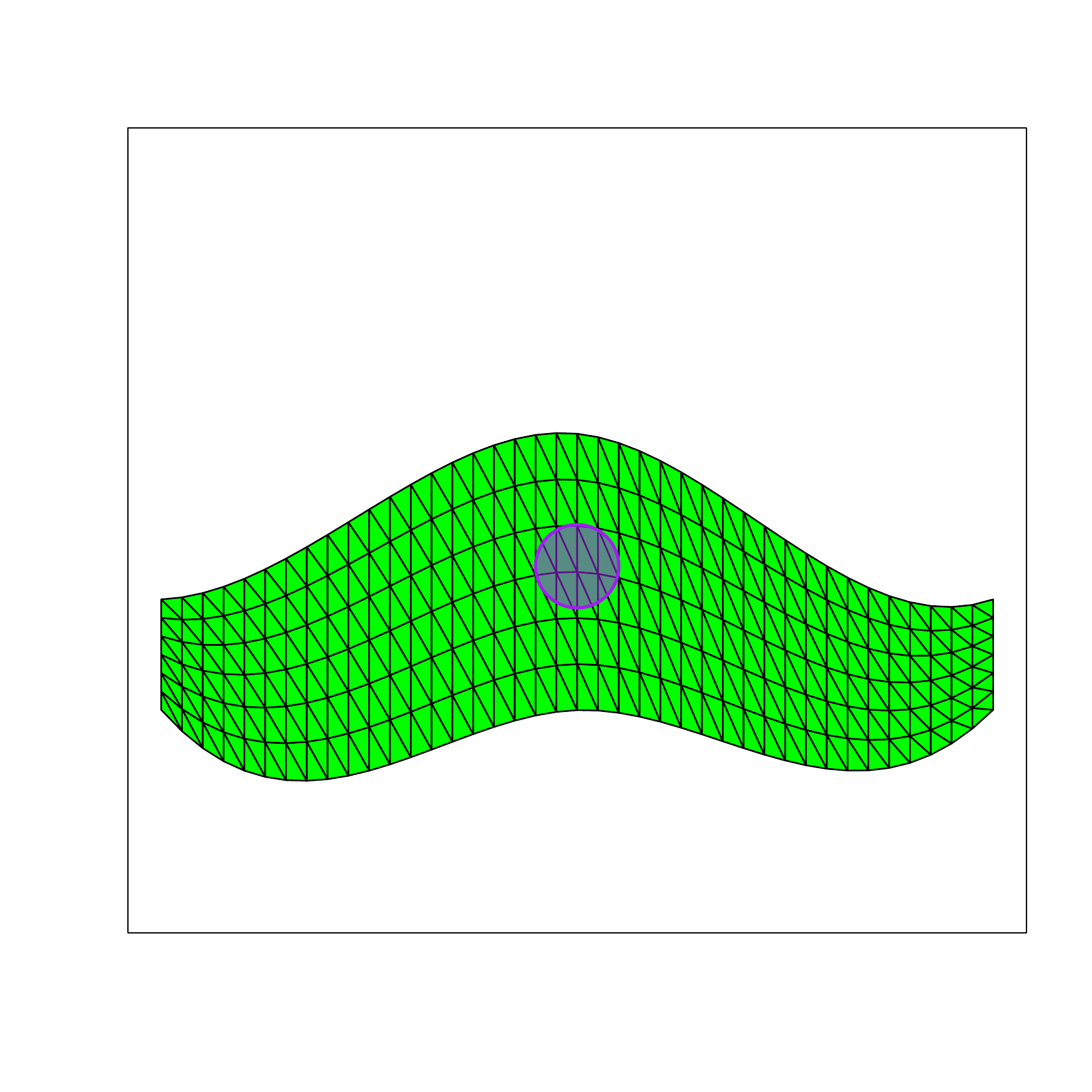}}
		\subfloat{\includegraphics[width=0.2\textwidth]{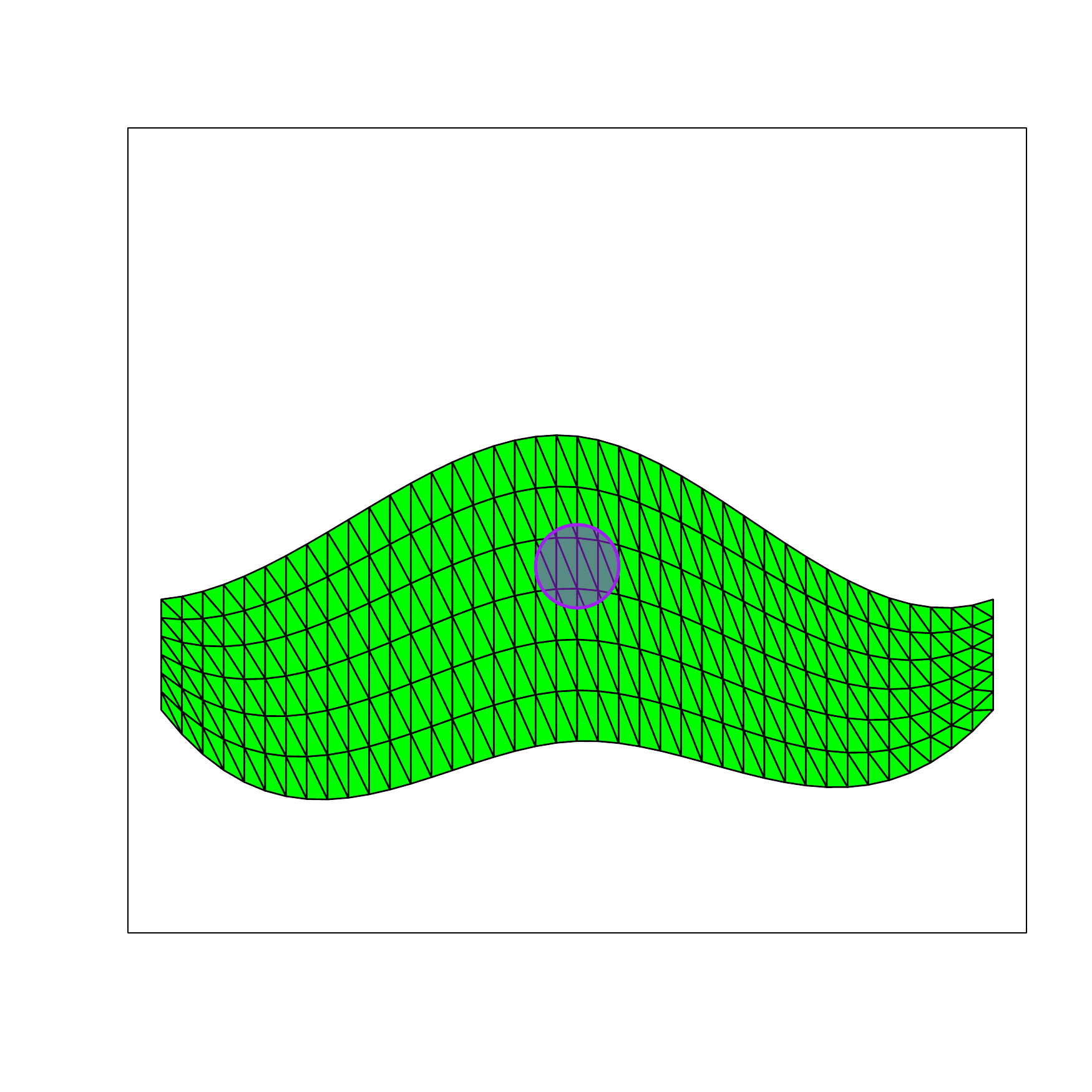}}
		\subfloat{\includegraphics[width=0.2\textwidth]{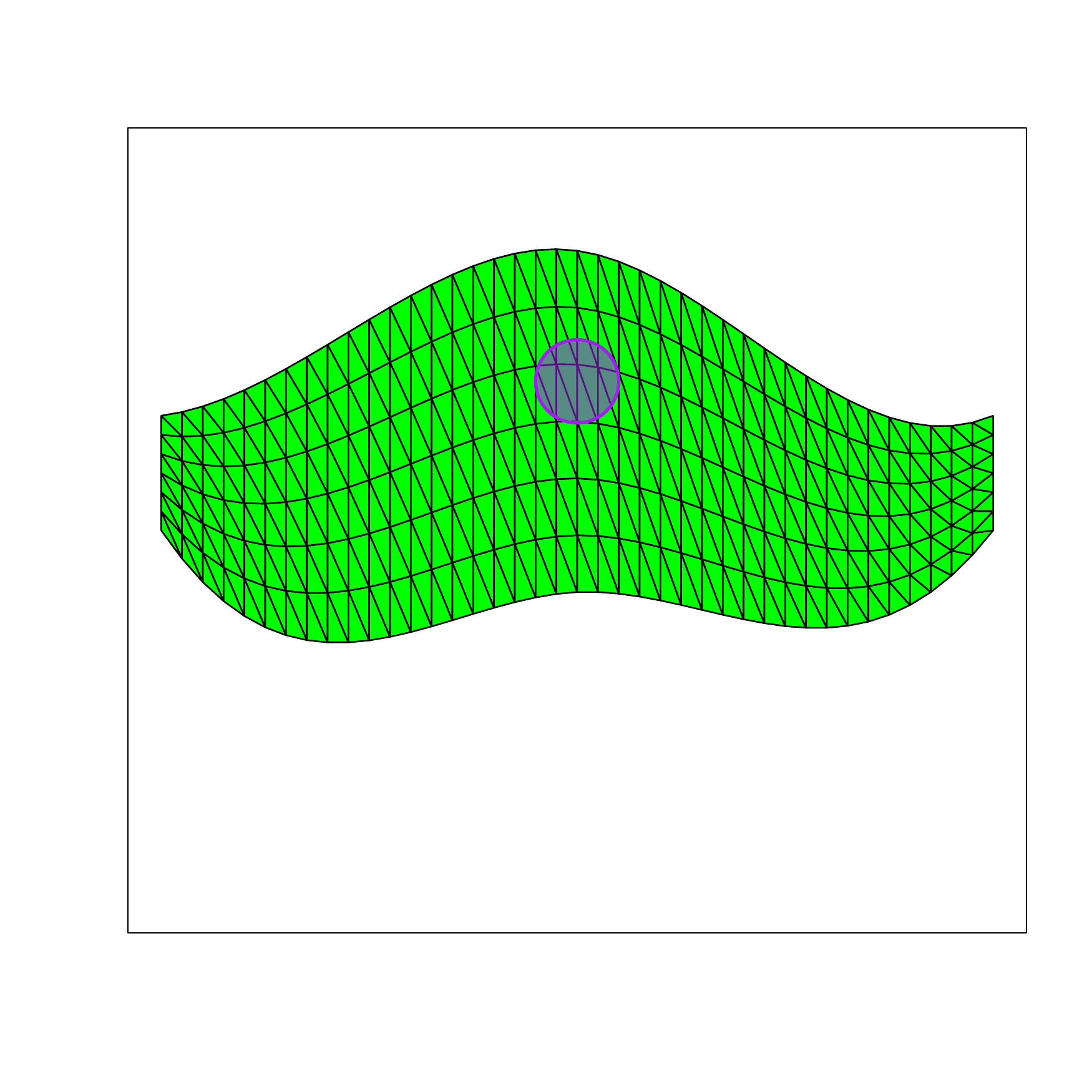}}

		\subfloat{\includegraphics[width=0.2\textwidth]{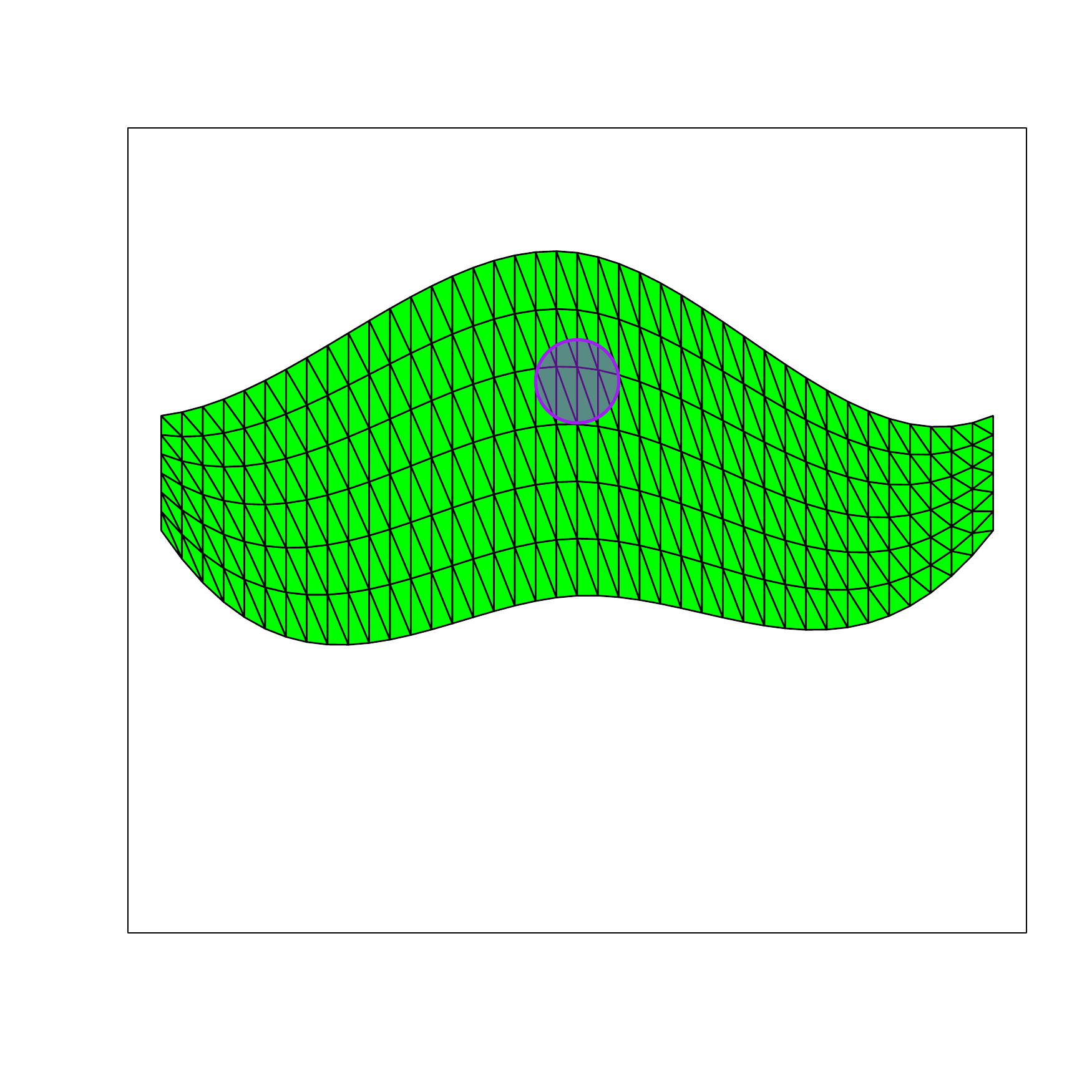}}
		\subfloat{\includegraphics[width=0.2\textwidth]{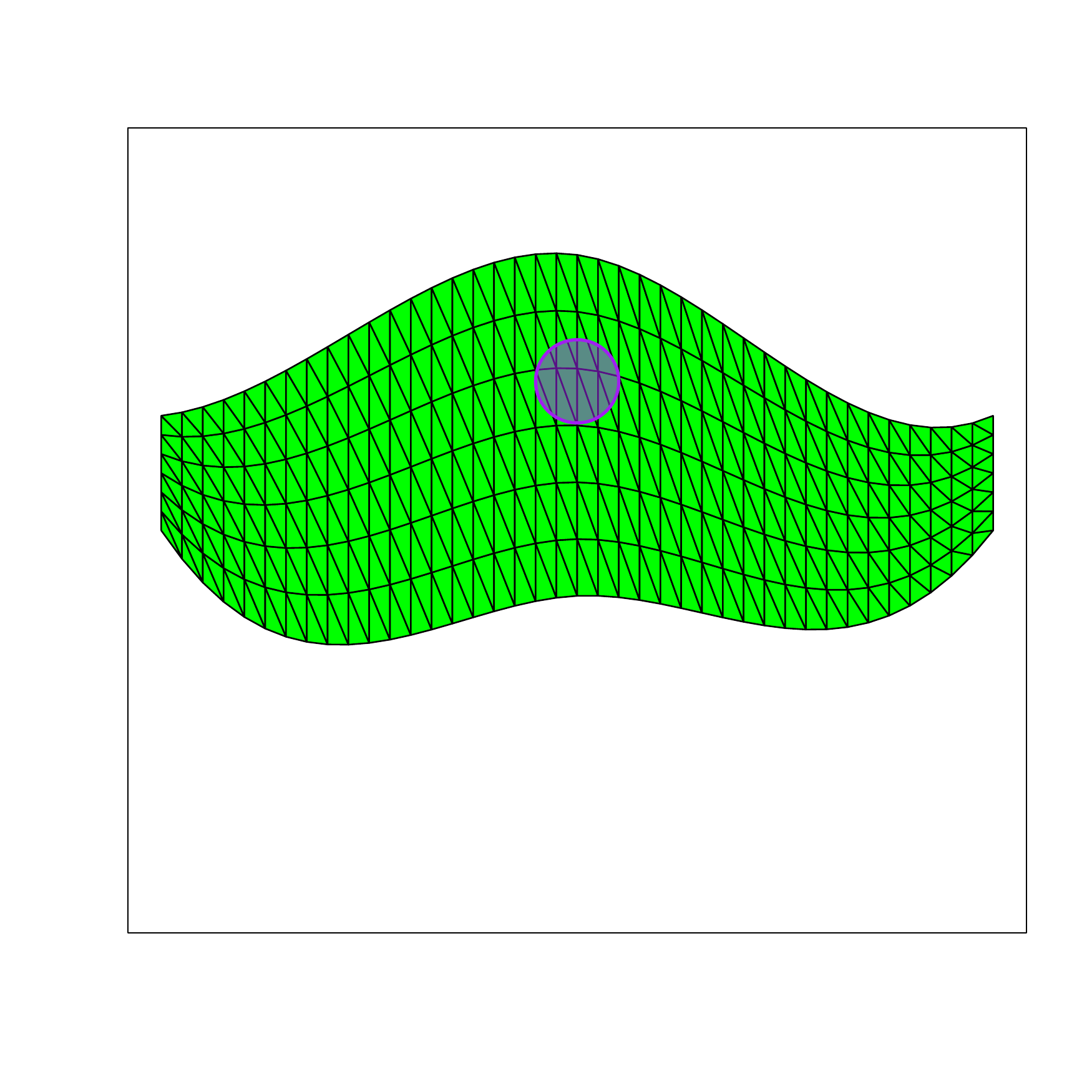}}
		\subfloat{\includegraphics[width=0.2\textwidth]{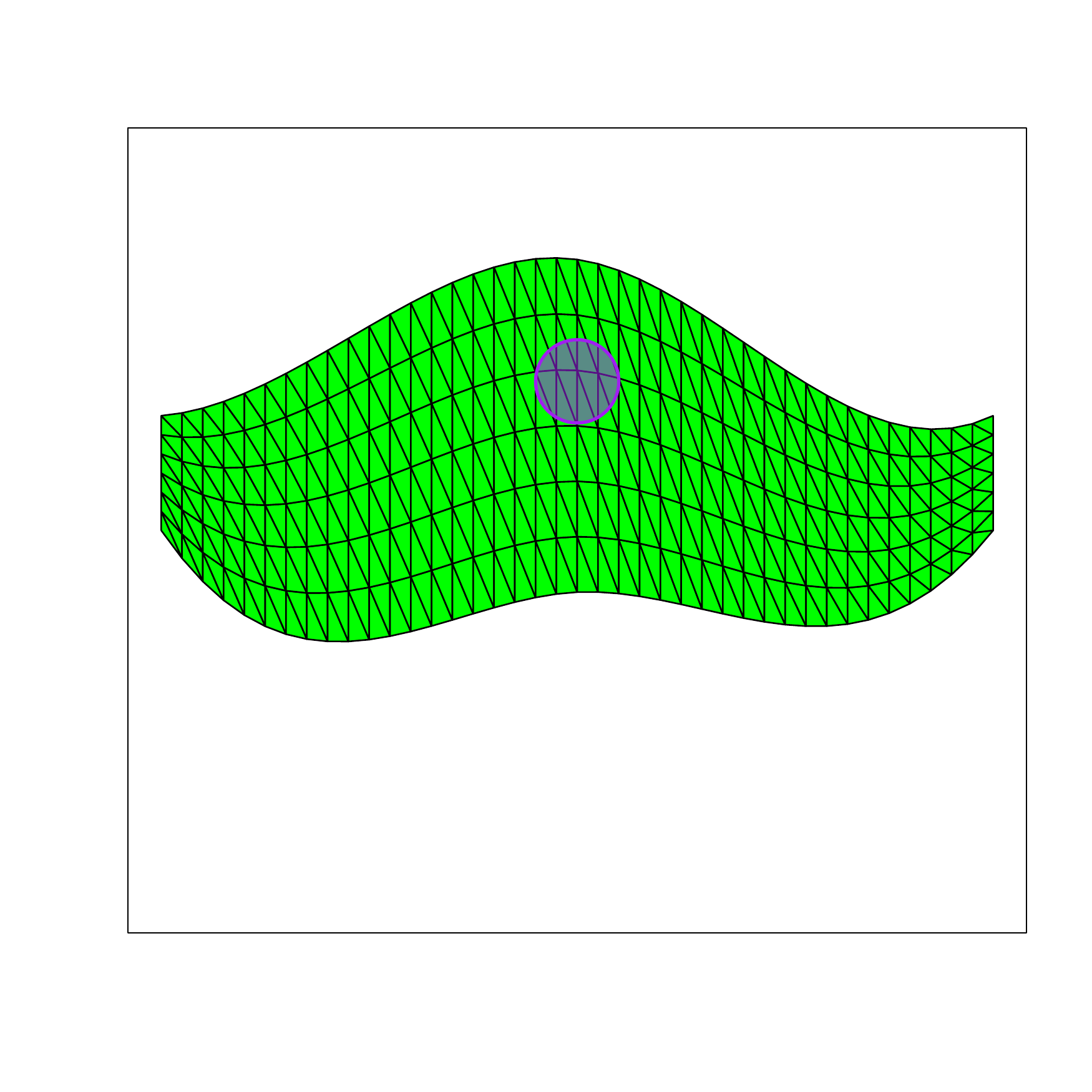}}
		\subfloat{\includegraphics[width=0.2\textwidth]{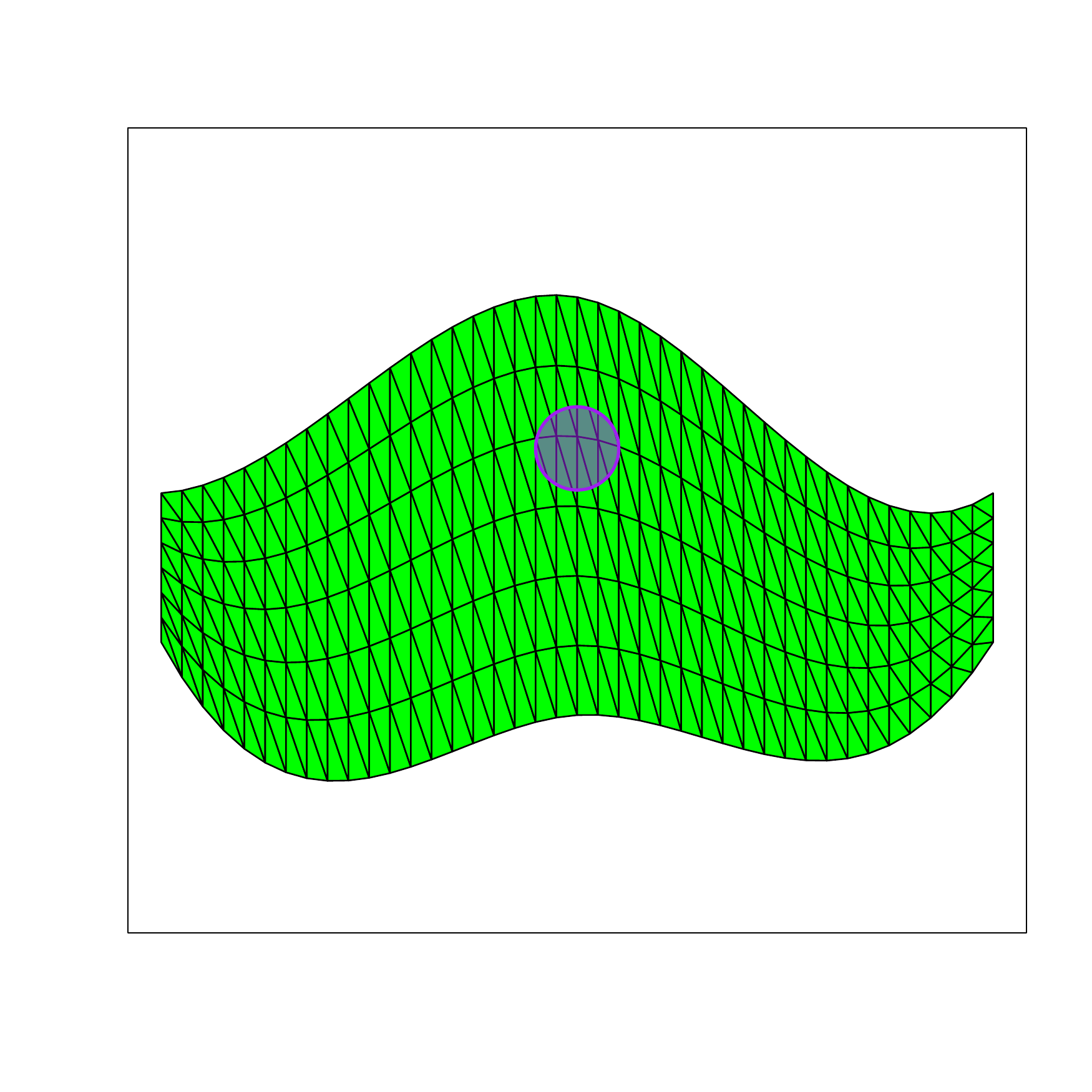}}
		\subfloat{\includegraphics[width=0.2\textwidth]{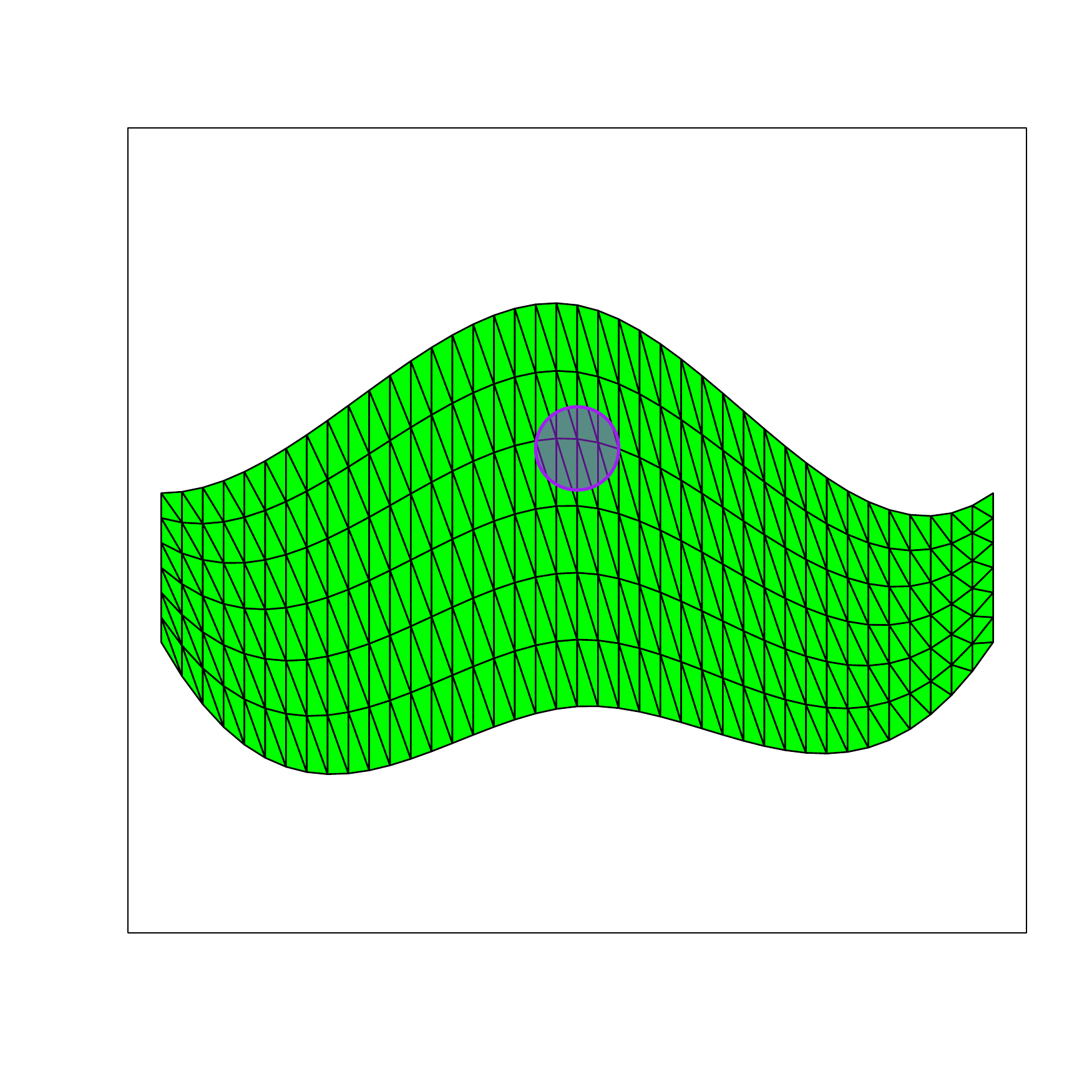}}

		\subfloat{\includegraphics[width=0.2\textwidth]{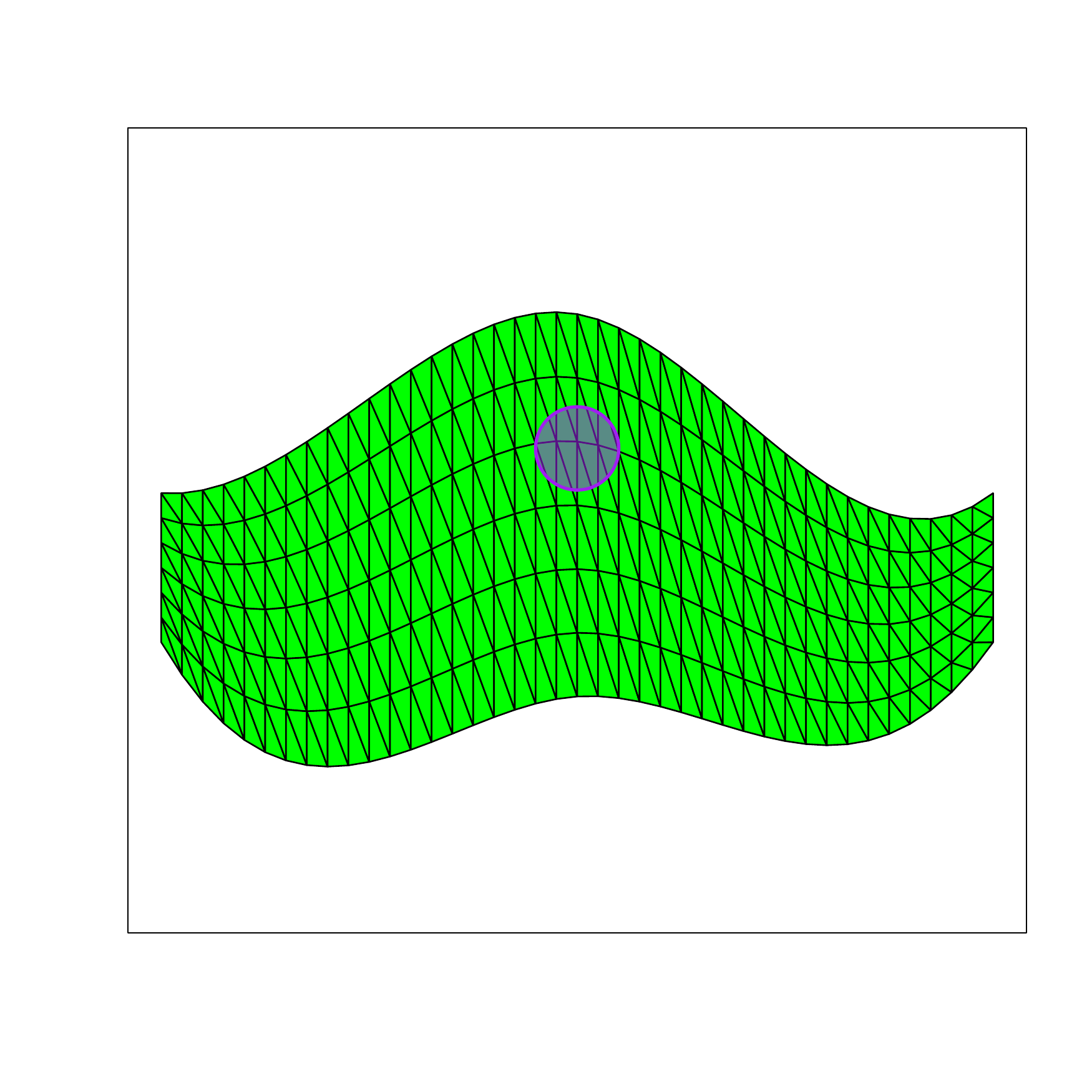}}
		\subfloat{\includegraphics[width=0.2\textwidth]{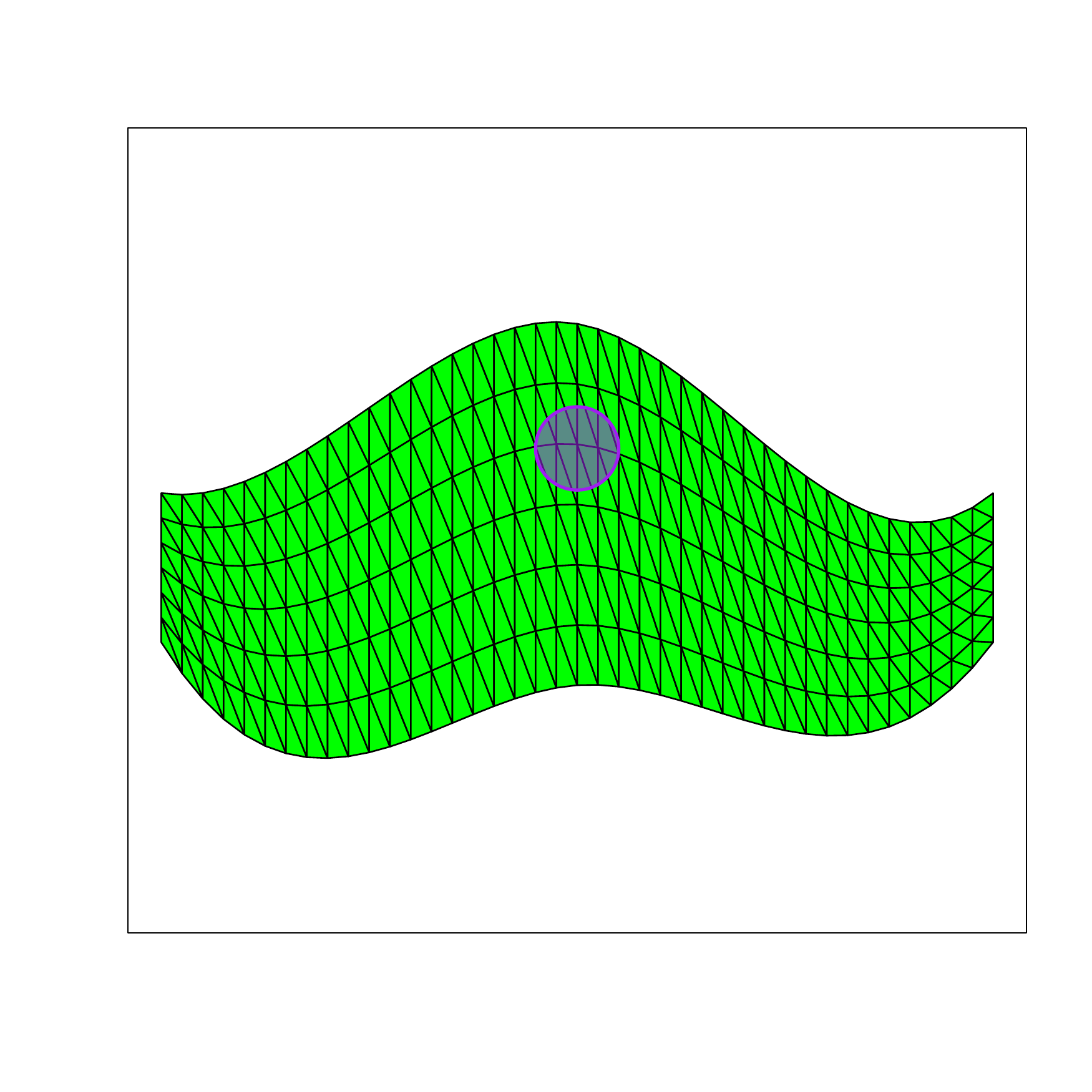}}
		\subfloat{\includegraphics[width=0.2\textwidth]{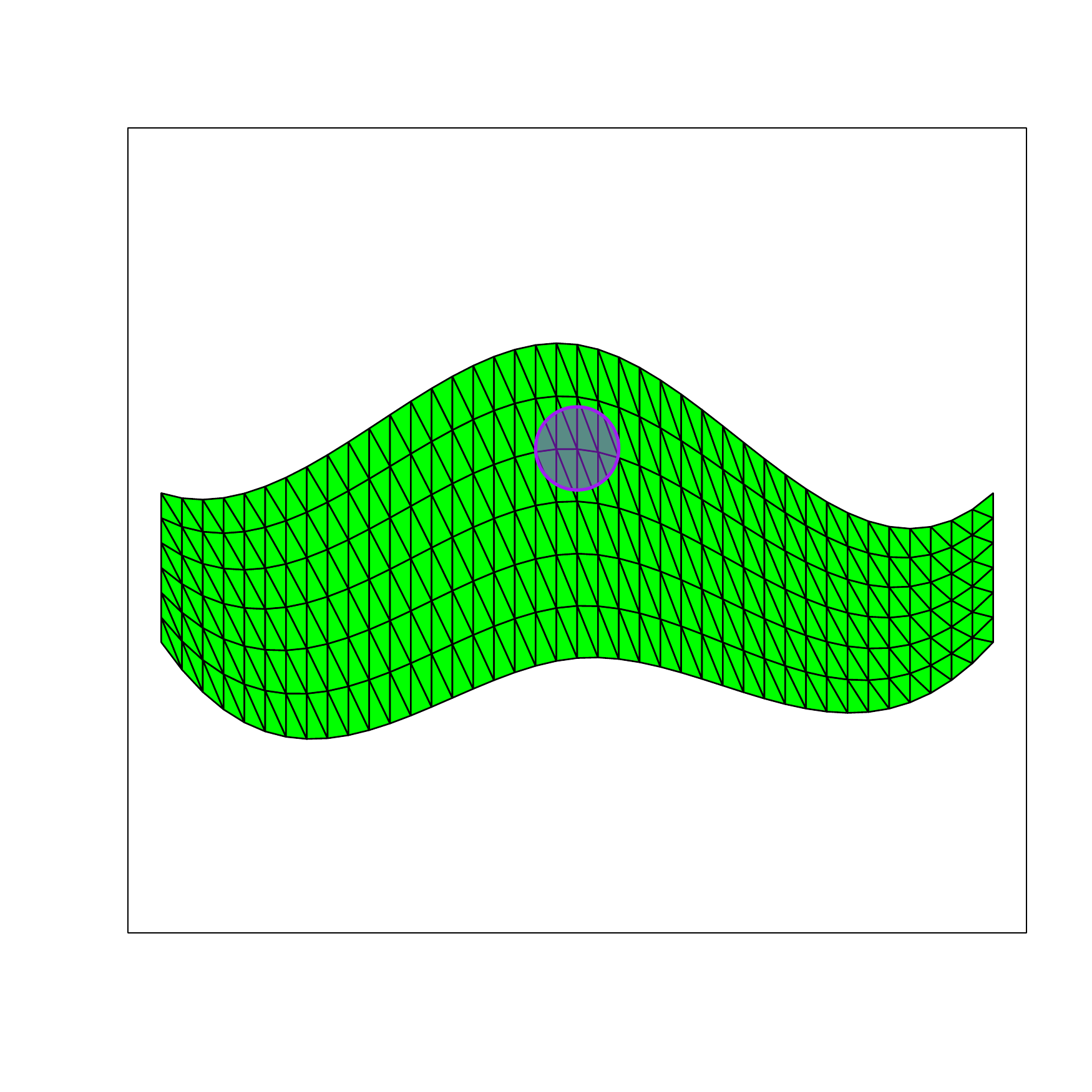}}
		\subfloat{\includegraphics[width=0.2\textwidth]{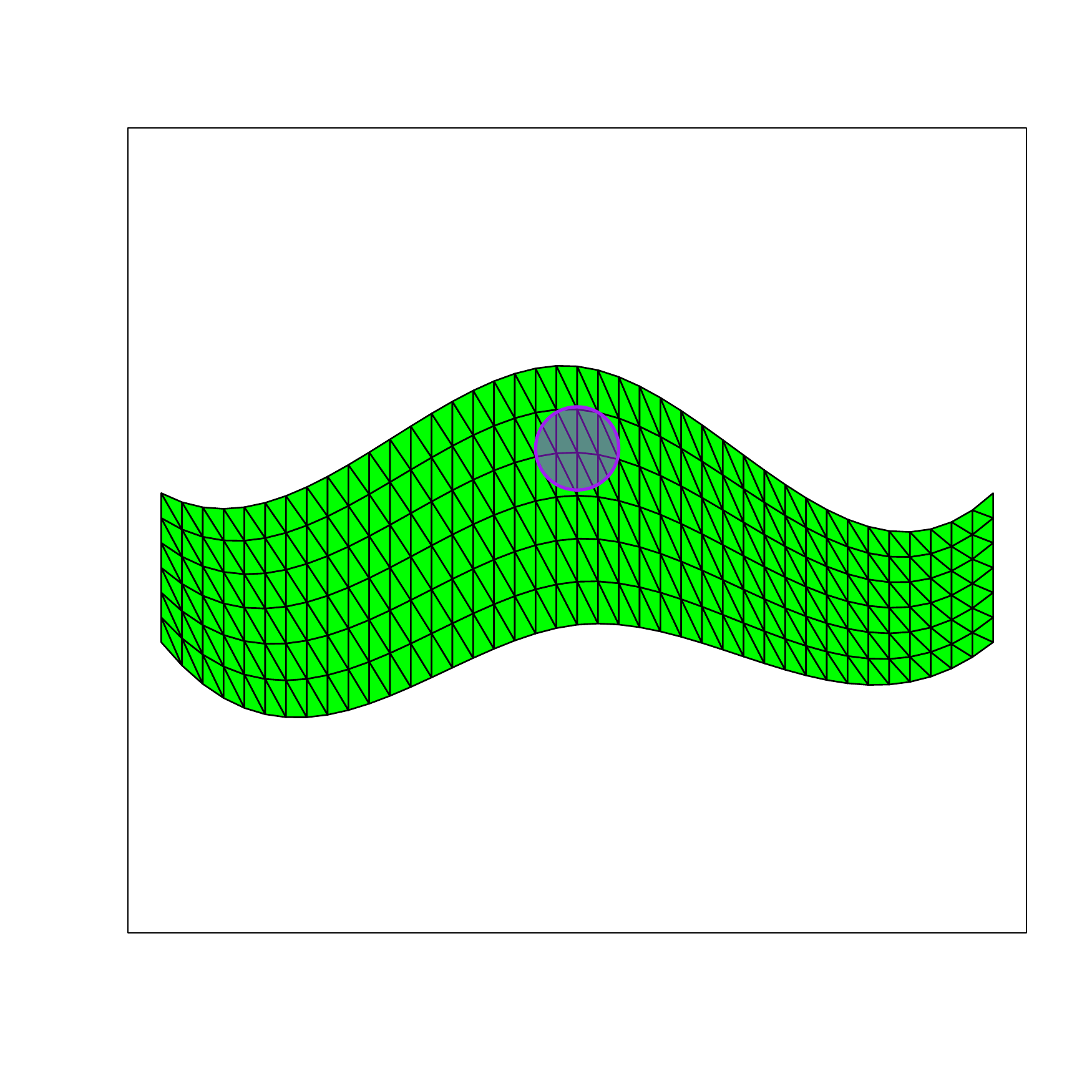}}
		\subfloat{\includegraphics[width=0.2\textwidth]{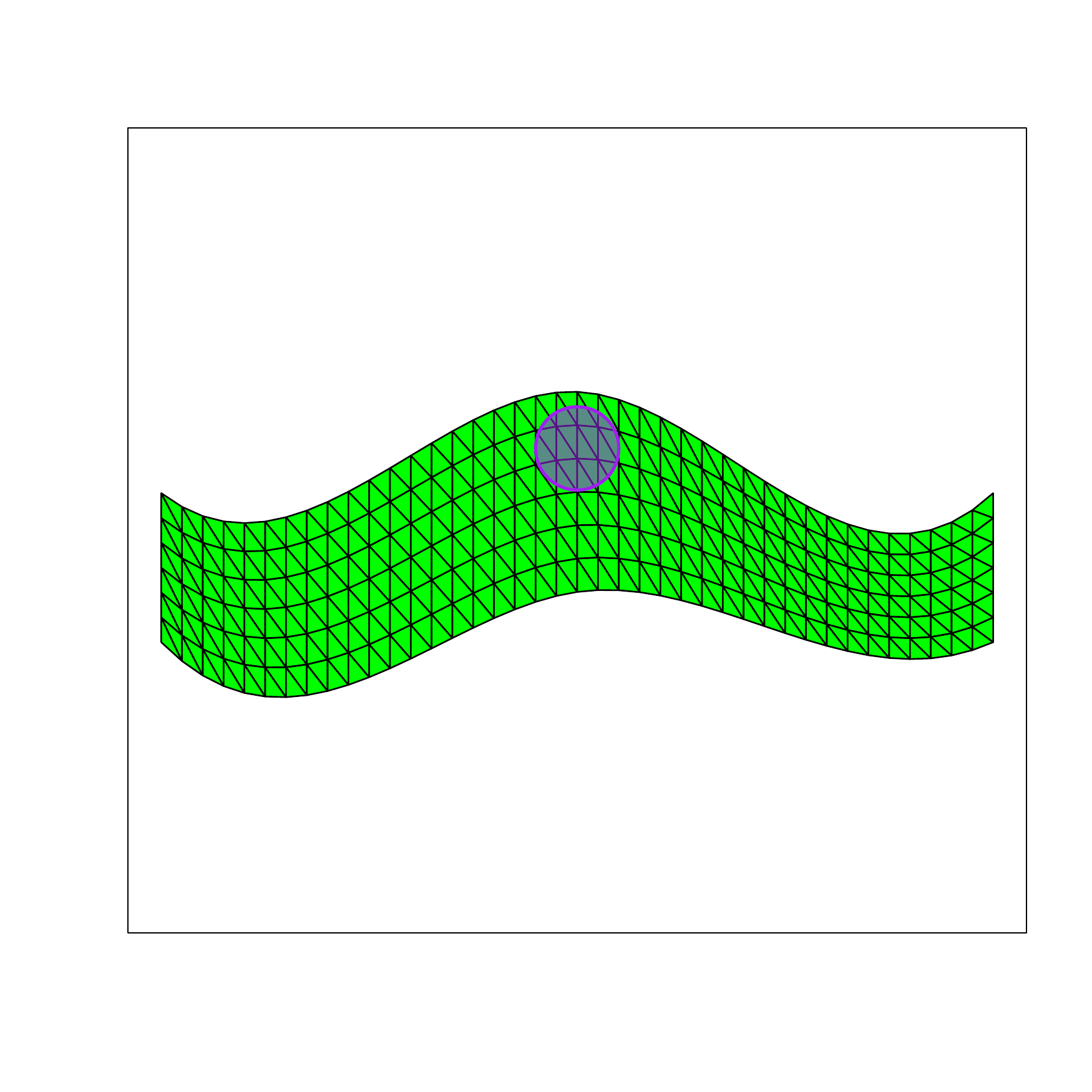}}

		\subfloat{\includegraphics[width=0.2\textwidth]{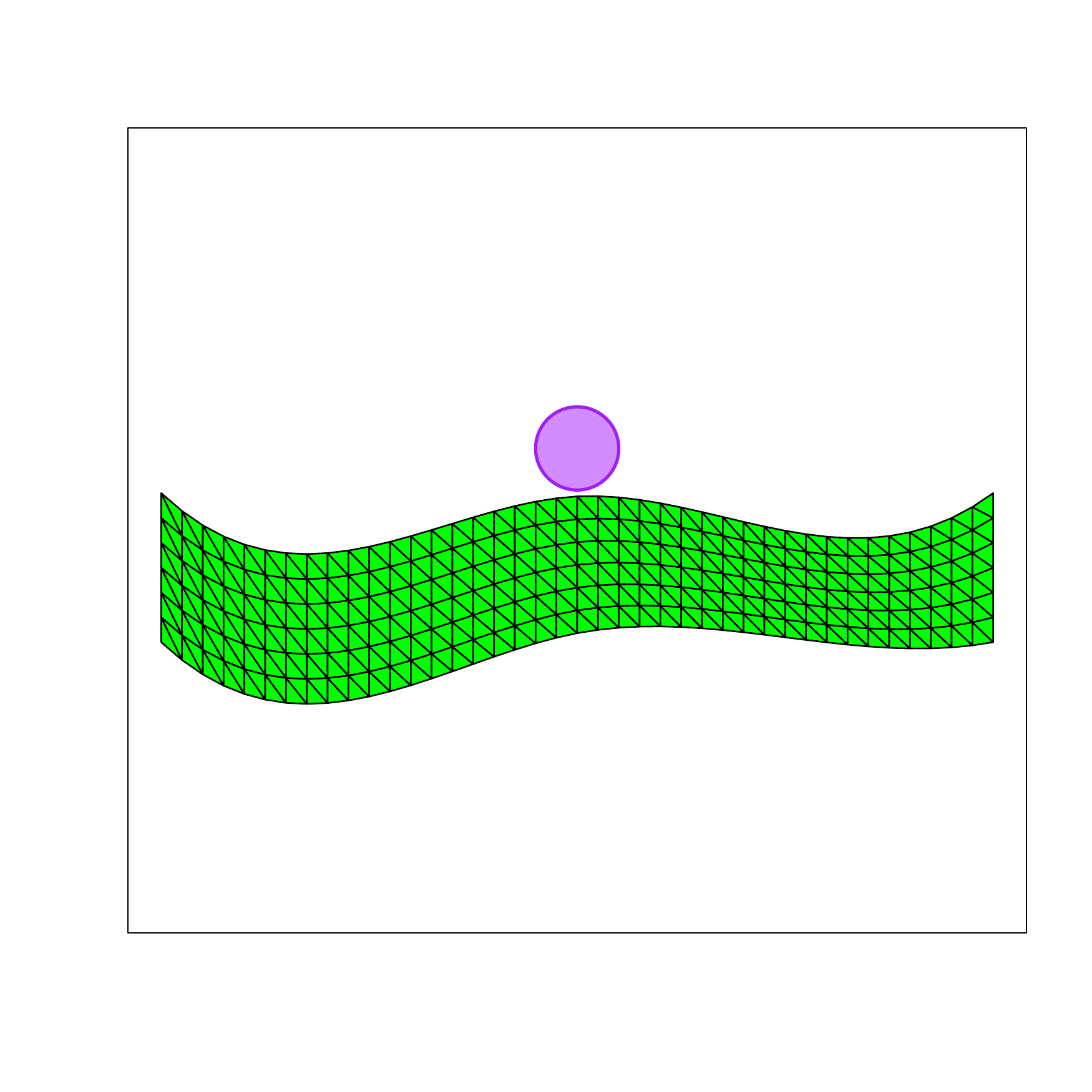}}
		\subfloat{\includegraphics[width=0.2\textwidth]{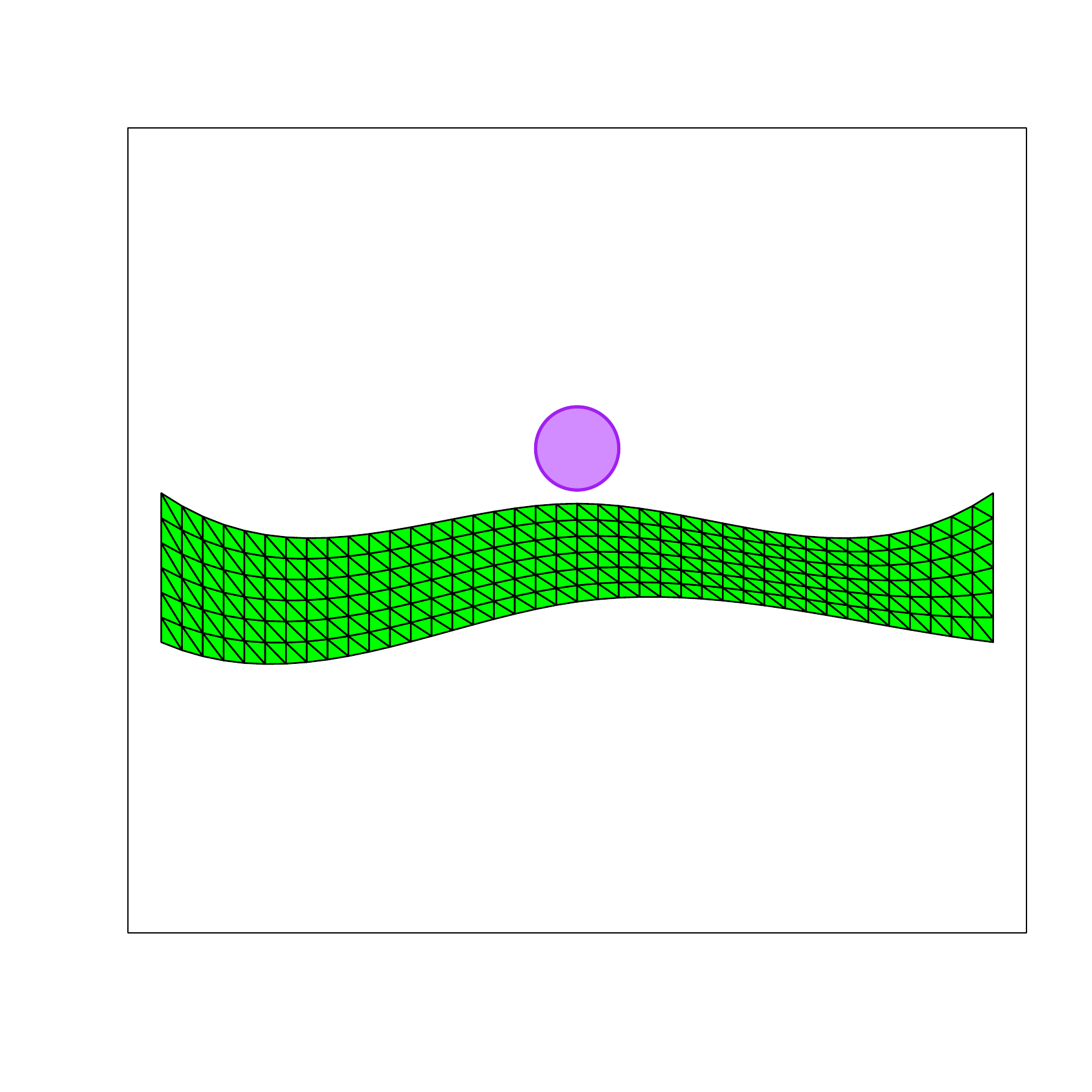}}
		\subfloat{\includegraphics[width=0.2\textwidth]{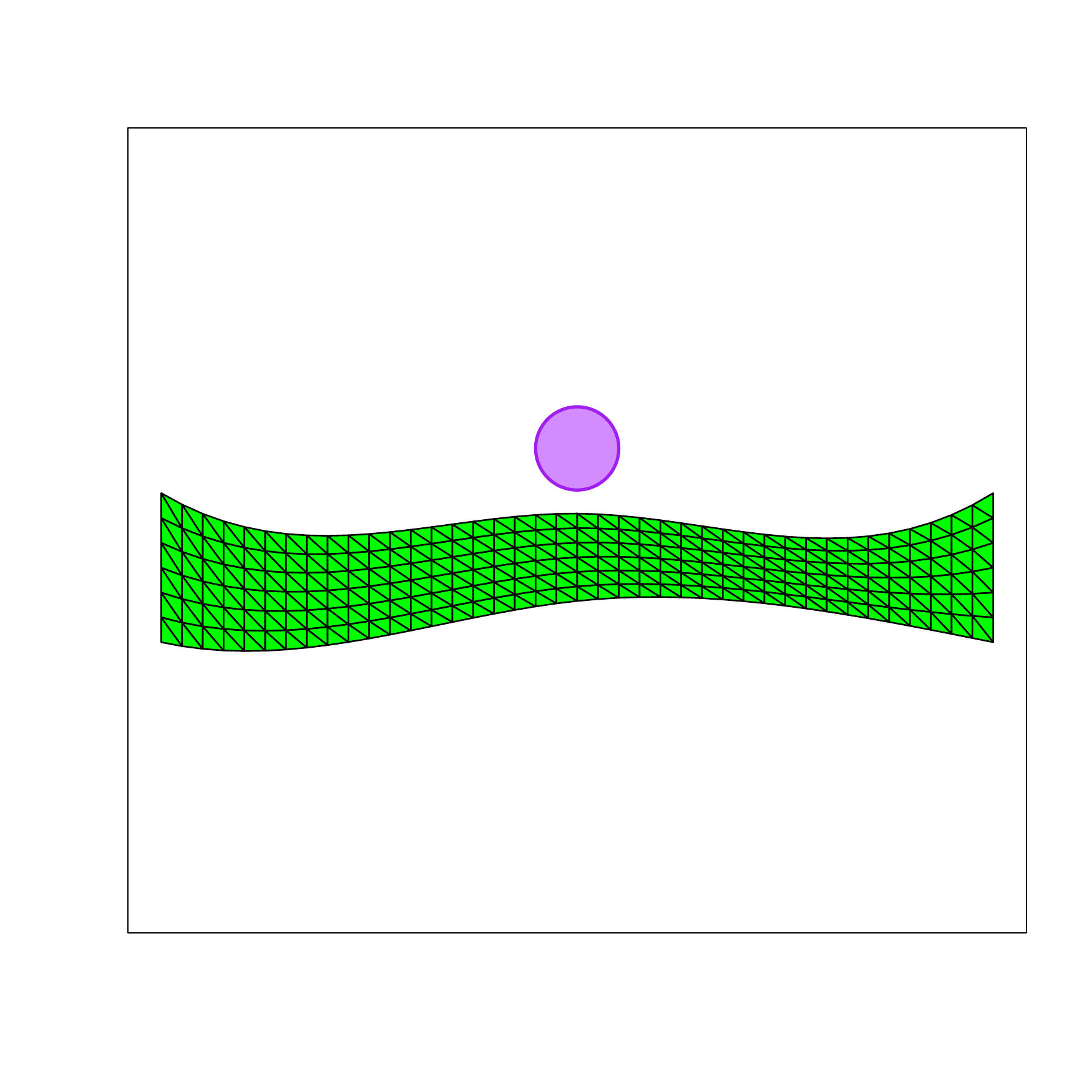}}
		\subfloat{\includegraphics[width=0.2\textwidth]{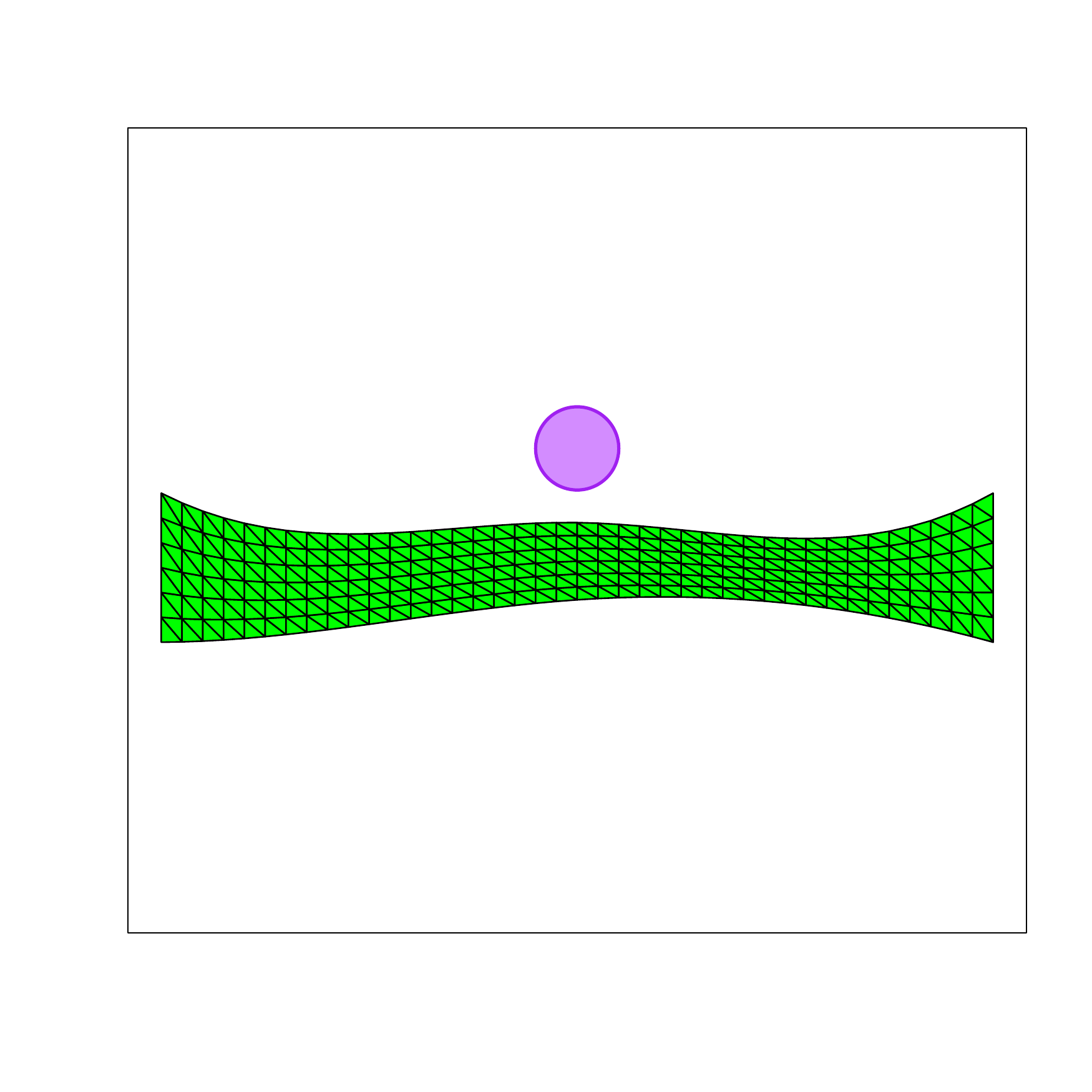}}
		\subfloat{\includegraphics[width=0.2\textwidth]{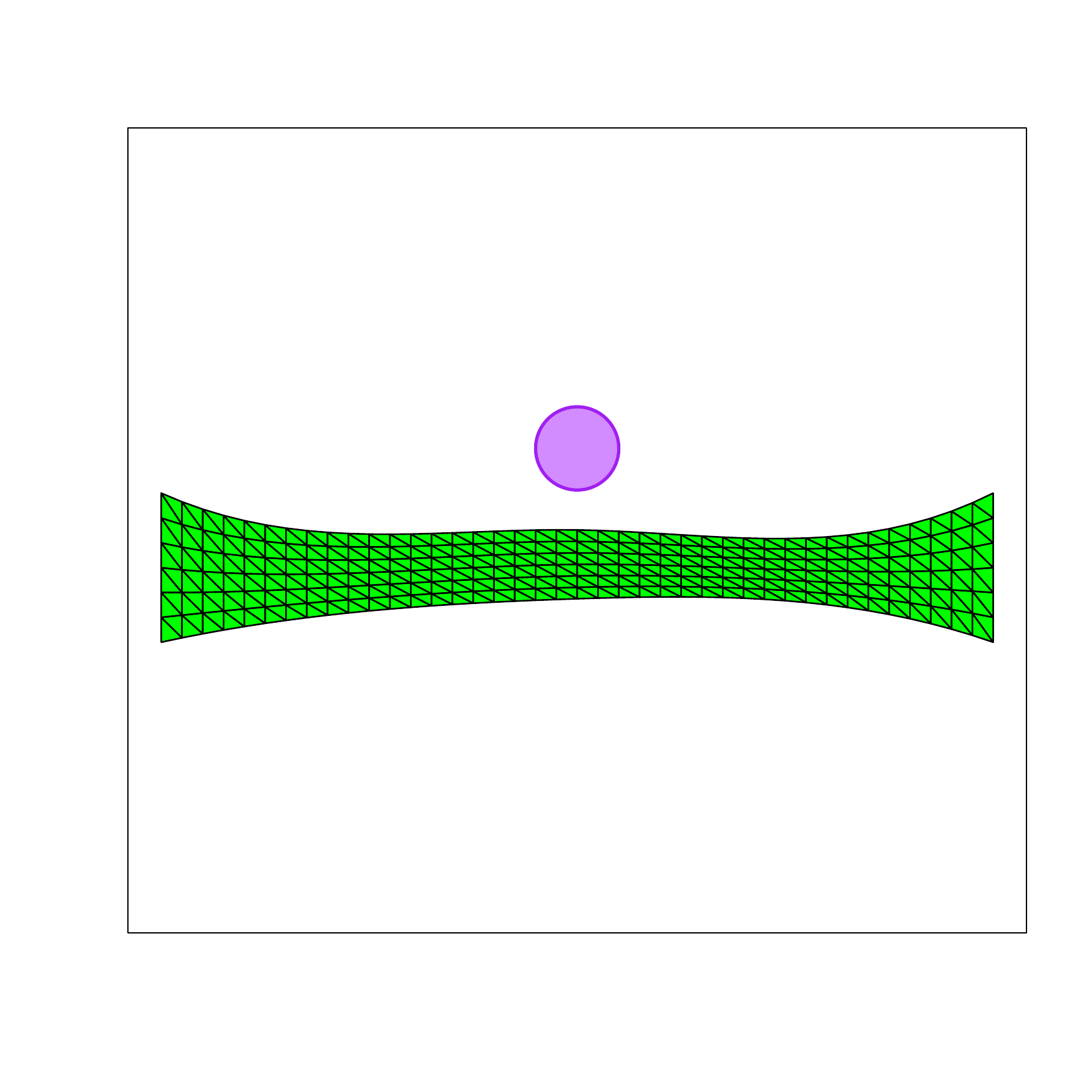}}

		\subfloat{\includegraphics[width=0.2\textwidth]{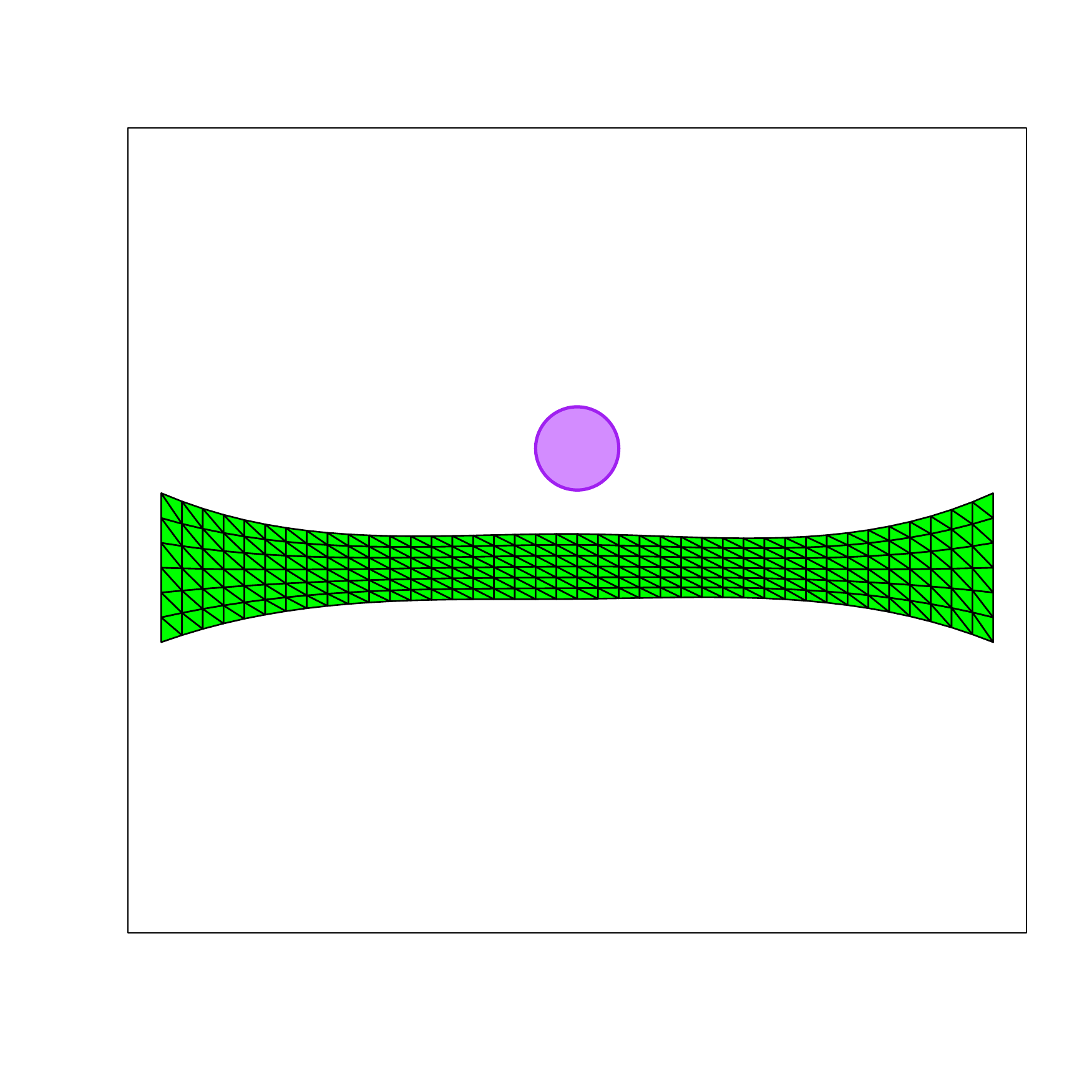}}
		\subfloat{\includegraphics[width=0.2\textwidth]{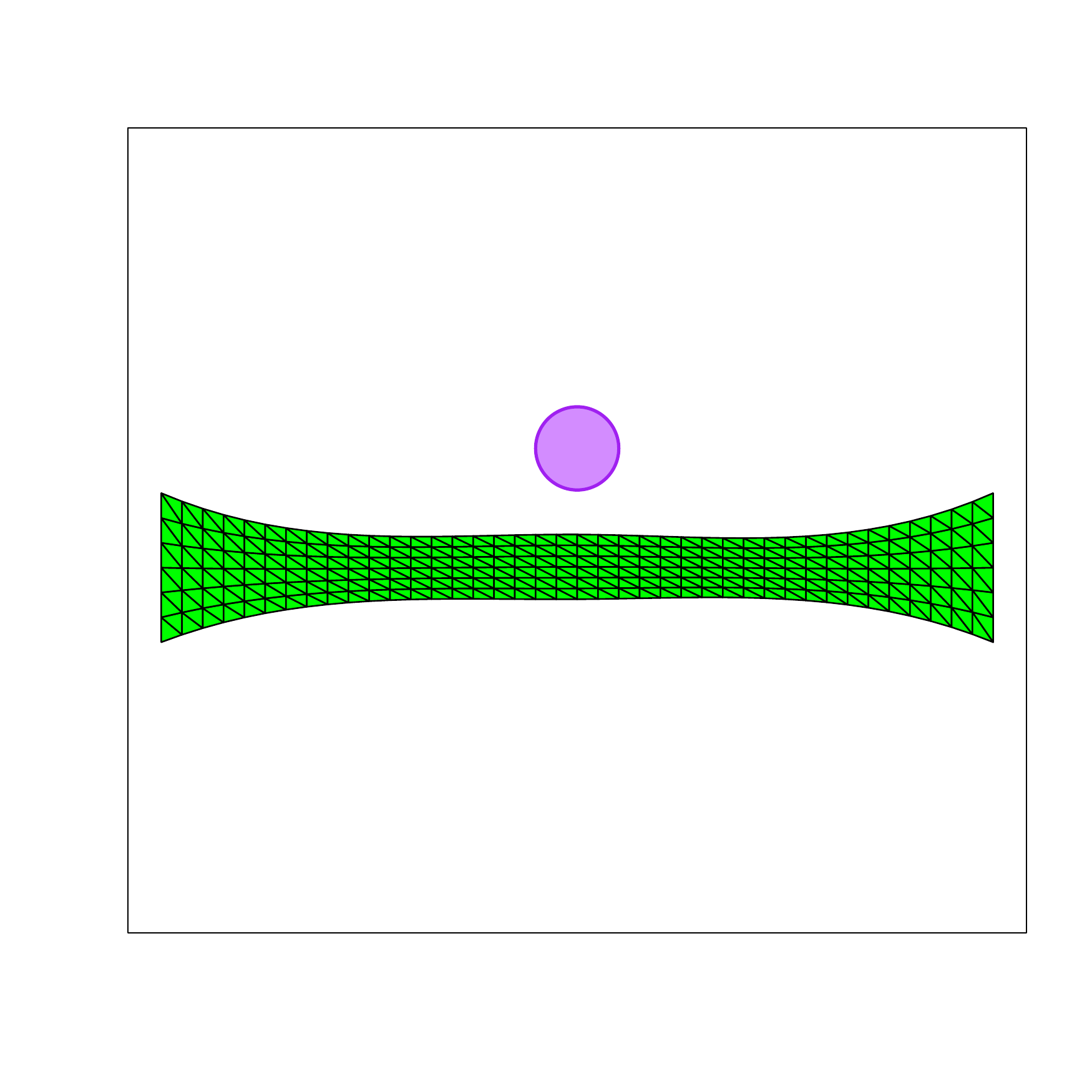}}
		\subfloat{\includegraphics[width=0.2\textwidth]{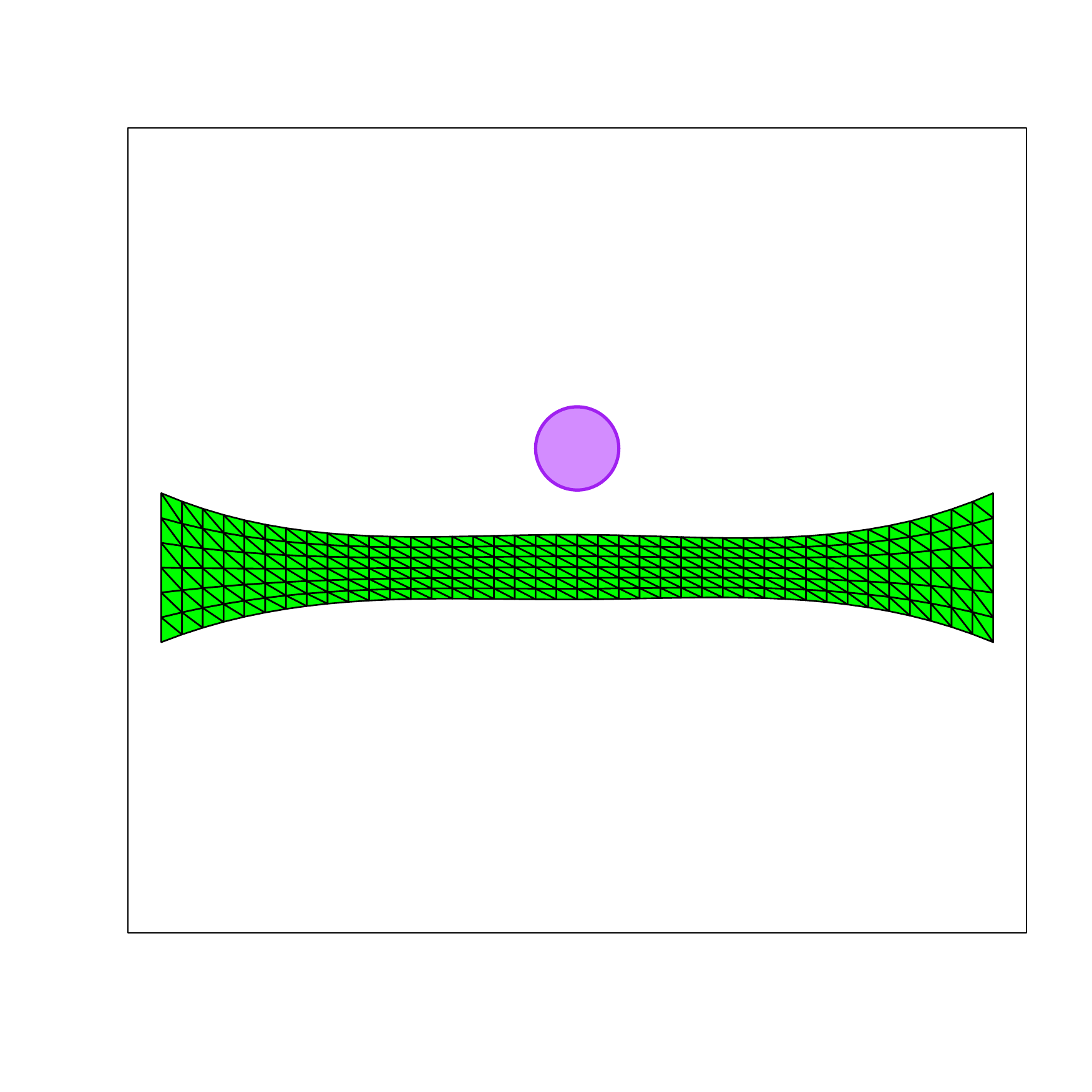}}
		\subfloat{\includegraphics[width=0.2\textwidth]{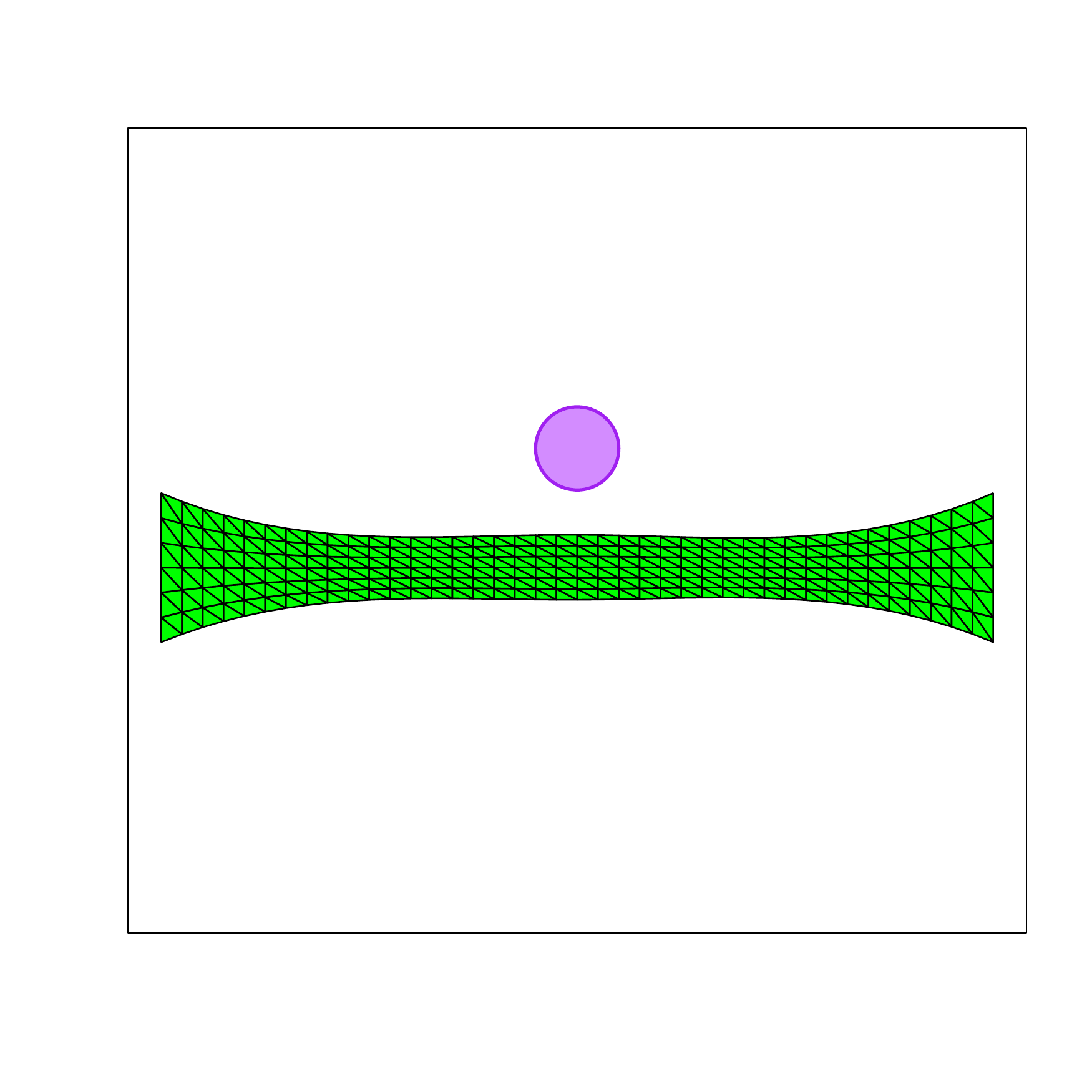}}
		\subfloat{\includegraphics[width=0.2\textwidth]{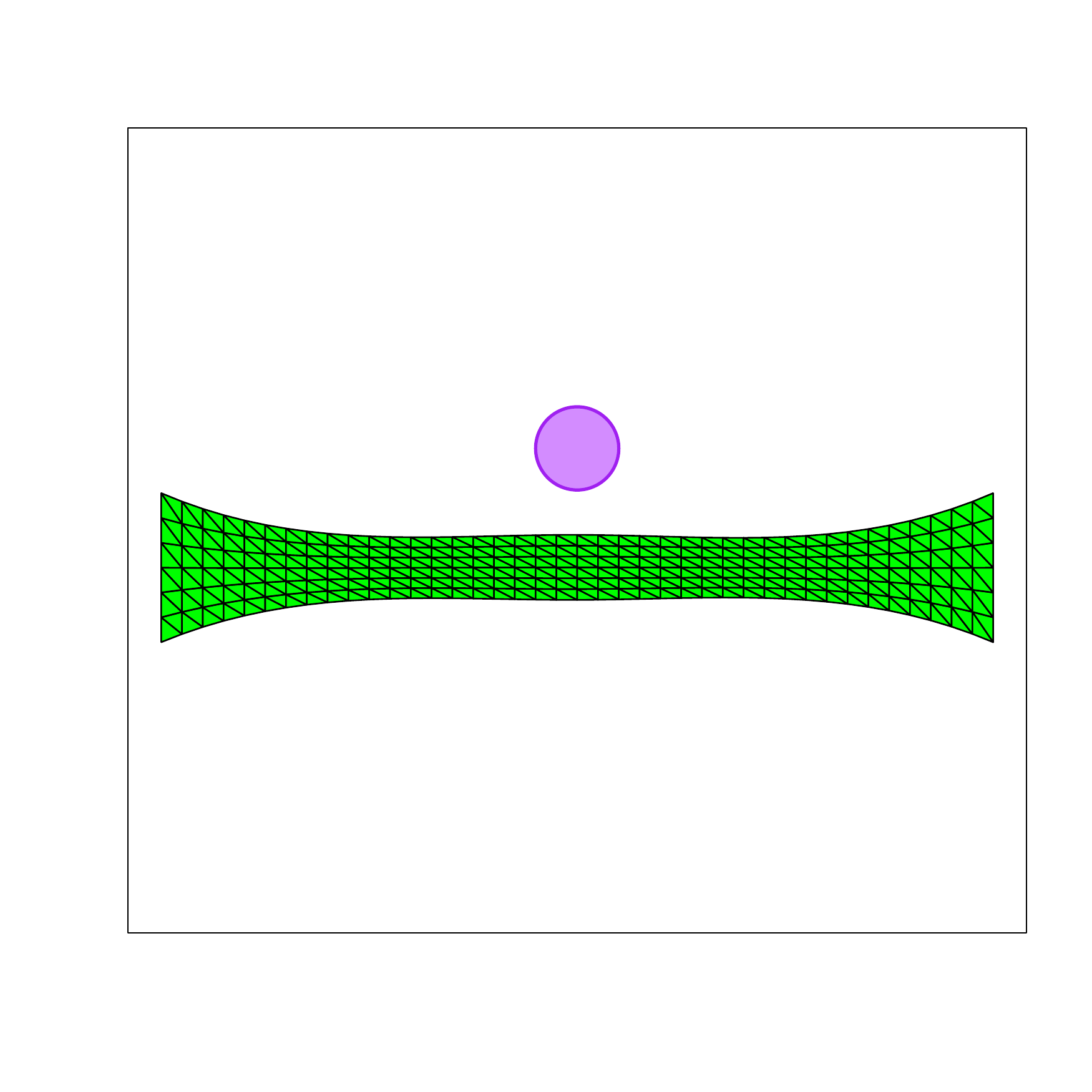}}

 	\end{center}		
	\caption{Test case 1 -- selected solutions \kk{of the discretized Hamiltonian flow} computed by the \kk{momentum method} 
    starting in $q^0$ and ending in $q^{\text{HF}}$. 
    Recall that the plotted circle indicates the correct position but does not always reflect its true radius, \kk{particularly for larger shapes when the $y$-axis is rescaled. \href{https://youtu.be/7vsPvzD7trs}{See video.}\label{fig:TC1_shapeflow}}}

\end{figure}

 \kk{Several exemplary} shapes that were computed \kk{as part of the discretized} the Hamiltonian flow are shown in Figure~\ref{fig:TC1_shapeflow}, capturing the course of the \kk{momentum method.} 
In Figure~\ref{fig:TC1_shapeflow}, the first \otd{shape} corresponds to the starting shape $q^0$ and the \otd{last shape to the solution $q^{\text{HF}}$, respectively}. The \kk{remaining} shapes are chosen to illustrate the \kk{Hamitonian} flow \kk{that is approximated} by the \kk{momentum method.} 

\subsubsection*{Test Case 2: An S-Shaped Joint} 
\label{subsubsec:S-Joint}
\begin{figure}[]
	\begin{center}
		\subfloat[Exemplary \kk{locally Pareto} optimal solution of \cite{bolten2021tracing}, resembling the contour of a whale, with added circular obstacle $\varpi$ (purple) placed above its fin.\label{subfig:TC2_PT_circ}]{\includegraphics[width=0.45\textwidth]{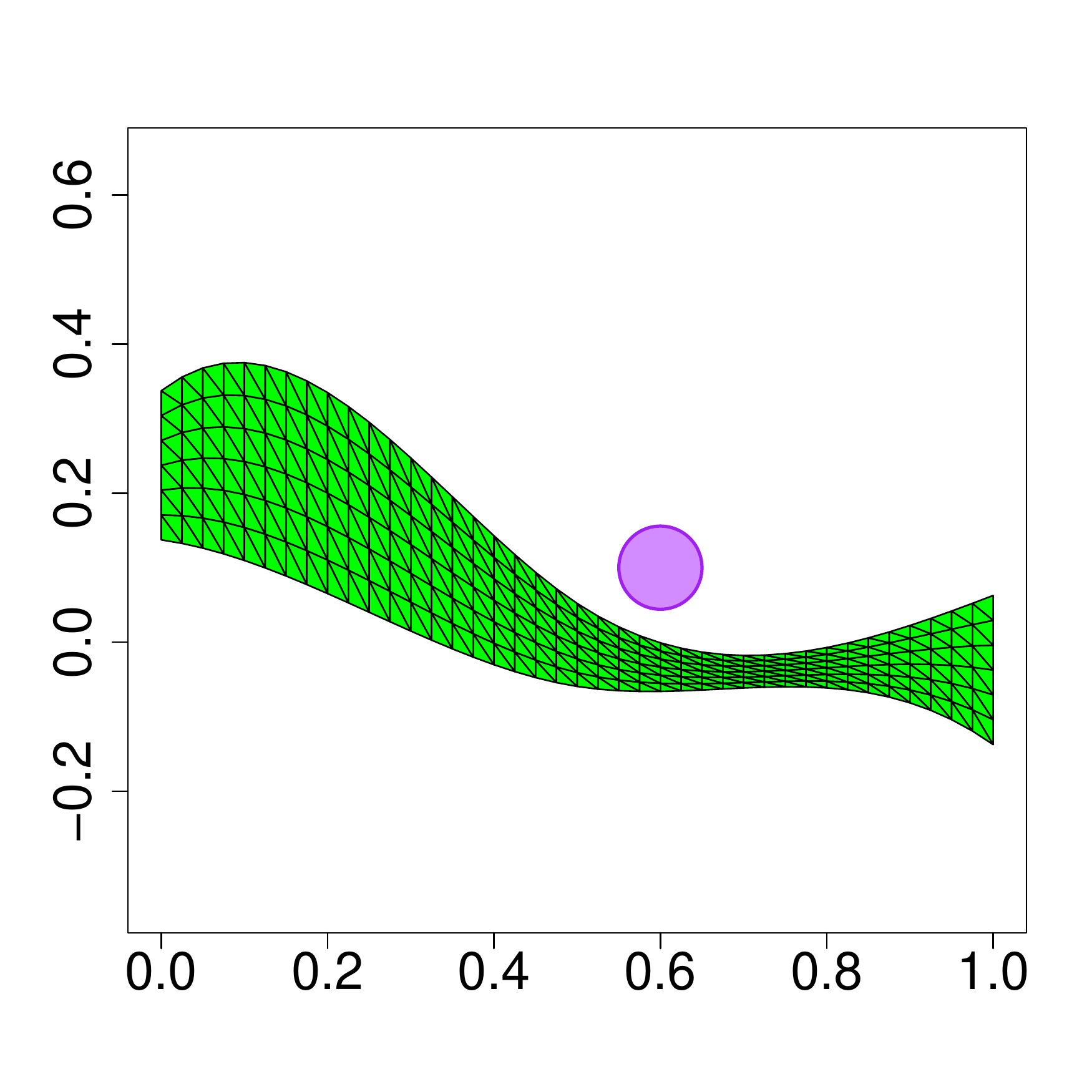}}
		\hspace{\fill}
		\subfloat[Initial shape $q^0$ for the \kk{gradient descent and momentum method} located above the obstacle $\varpi$. \label{subfig:TC2_initShape}
		]{\includegraphics[width=0.45\textwidth]{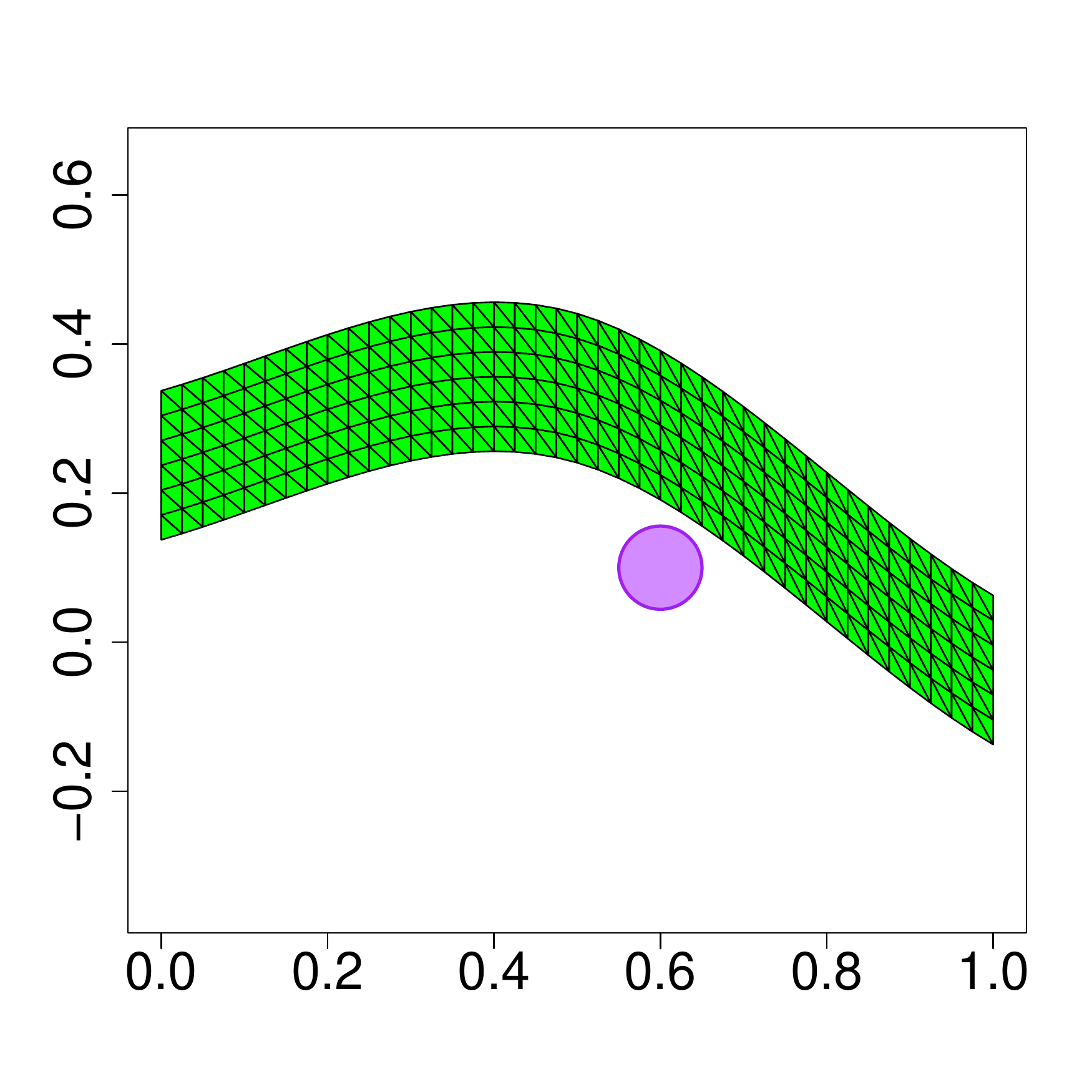}}
 	\end{center}		
	\caption{Test case 2 -- \kk{placement of circular obstacle $\varpi$ and initial shape.}\label{fig:TC2_start}}
\end{figure}

Following \cite{bolten2021tracing,Doganay2019}, for the second test case \kk{we consider an S-shaped joint} the right boundary of which is located $0.27\,\text{m}$ and hence beneath its left boundary. 
The locally \kk{Pareto} optimal solutions of the unpenalized biobjective shape optimization problem investigated in \cite{bolten2021tracing,Doganay2019} resemble the profiles of whales with varying volume. We now place a circular obstacle with midpoint $x_{\text{mp}}=(0.6, 0.1)$ and radius $r=0.05$, i.e., $\varpi=\varpi((0.6, 0.1), 0.05)$, above the right part of \kk{one exemplary locally Pareto optimal shape}, see Figure~\ref{subfig:TC2_PT_circ}, and choose the initial shape such that it is located above the circular obstacle $\varpi$, see Figure~\ref{subfig:TC2_initShape}.

For the numerical experiments we set the penalizing parameter $c_\text{P}$ of $J_3$ to $c_\text{P}=30$ \otd{and choose the dissipation parameter as $\gamma=10$}. The remaining parameters are identical to the ones used in the first test case, i.e., \otd{the weight vector $\lambda$ is $\lambda=(0.4, 0.3, 0.3)$}, the gradient descent method implements Armijo step lengths, starts in $q^0$ and has the same stopping condition and maximum number of iterations as before, and the parameters for the \kk{momentum method are} $\kappa=10^{-3}$, $m=10$, \otd{$T=1$ and $250$ time steps}, with the initial momentum set to $p^0=0\in\R^6$.

\begin{figure}[h]
	\begin{center}
		\subfloat[\kk{Final} solution $q^{\text{GD}}$ of the gradient descent method \kk{after $200$ iterations,} starting in $q^0$. Note that the actual intersection are with $\varpi$ is much smaller than what the plot suggests.
  \label{subfig:TC2_solGD}]{\includegraphics[width=0.45\textwidth]{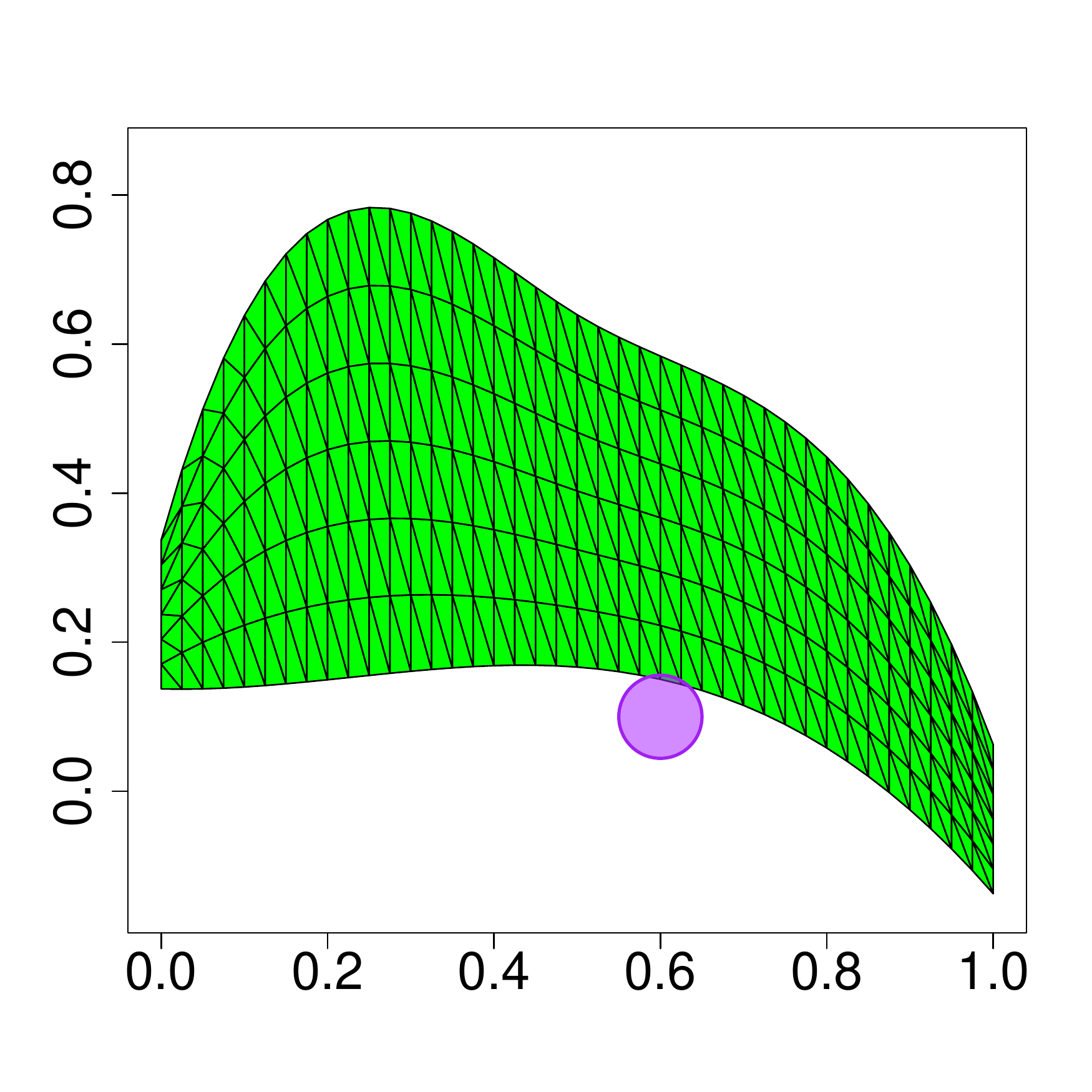}}
		\hspace{\fill}
		\subfloat[\kk{Final} solution $q^{\text{HF}}$ of the \kk{momentum method after \otd{$250$} iterations} starting in $q^0$. 
        \label{subfig:TC2_soldHF}
		]{\includegraphics[width=0.45\textwidth]{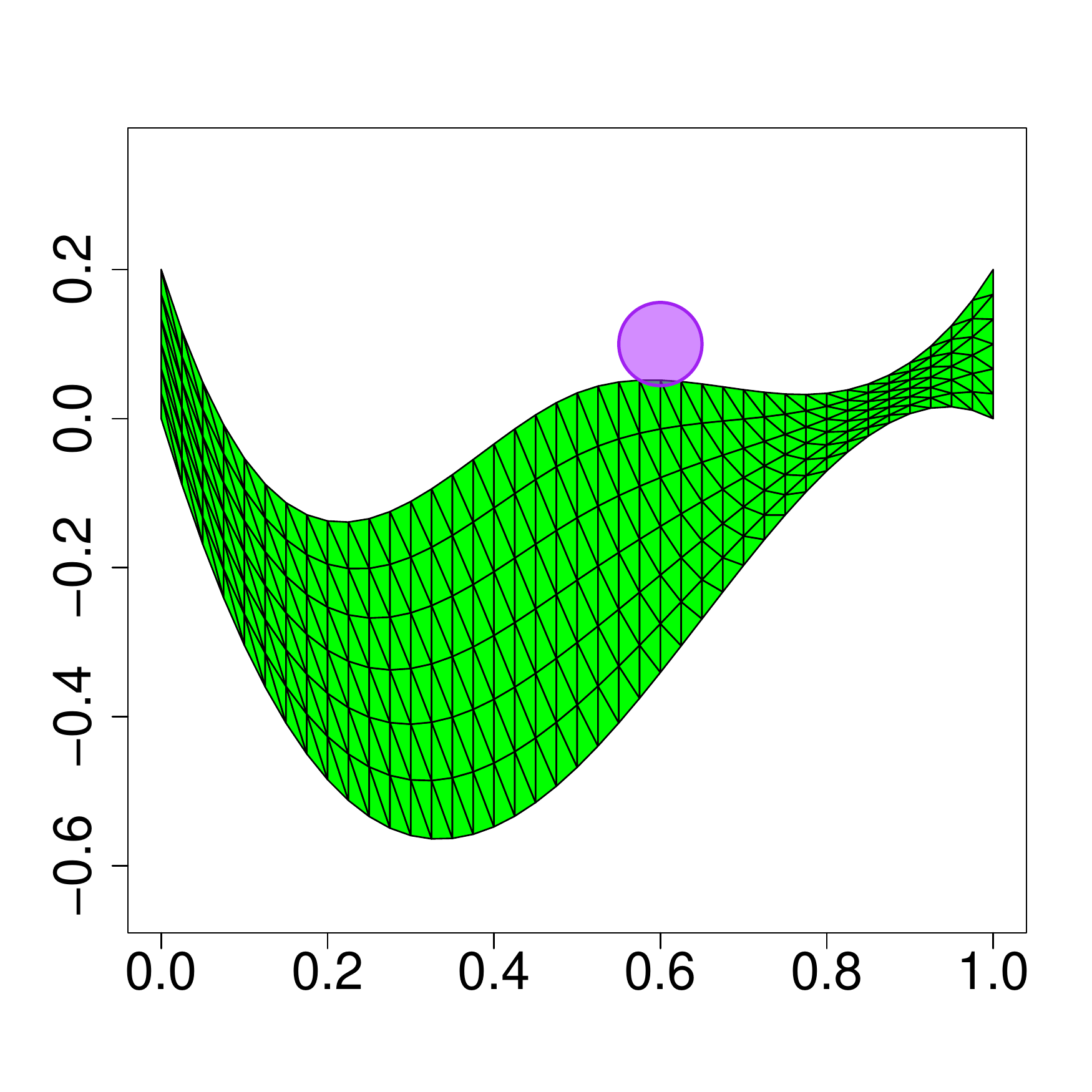}}
 	\end{center}		
	\caption{Test case 2 -- comparison of the results of gradient descent and \kk{momentum method.}\label{fig:TC2_comp}}
\end{figure}


\begin{figure}[h]
	\begin{center}
		\subfloat[Test case 1 \label{fig:TC1_Ehist}]{\includegraphics[width=0.45\textwidth]{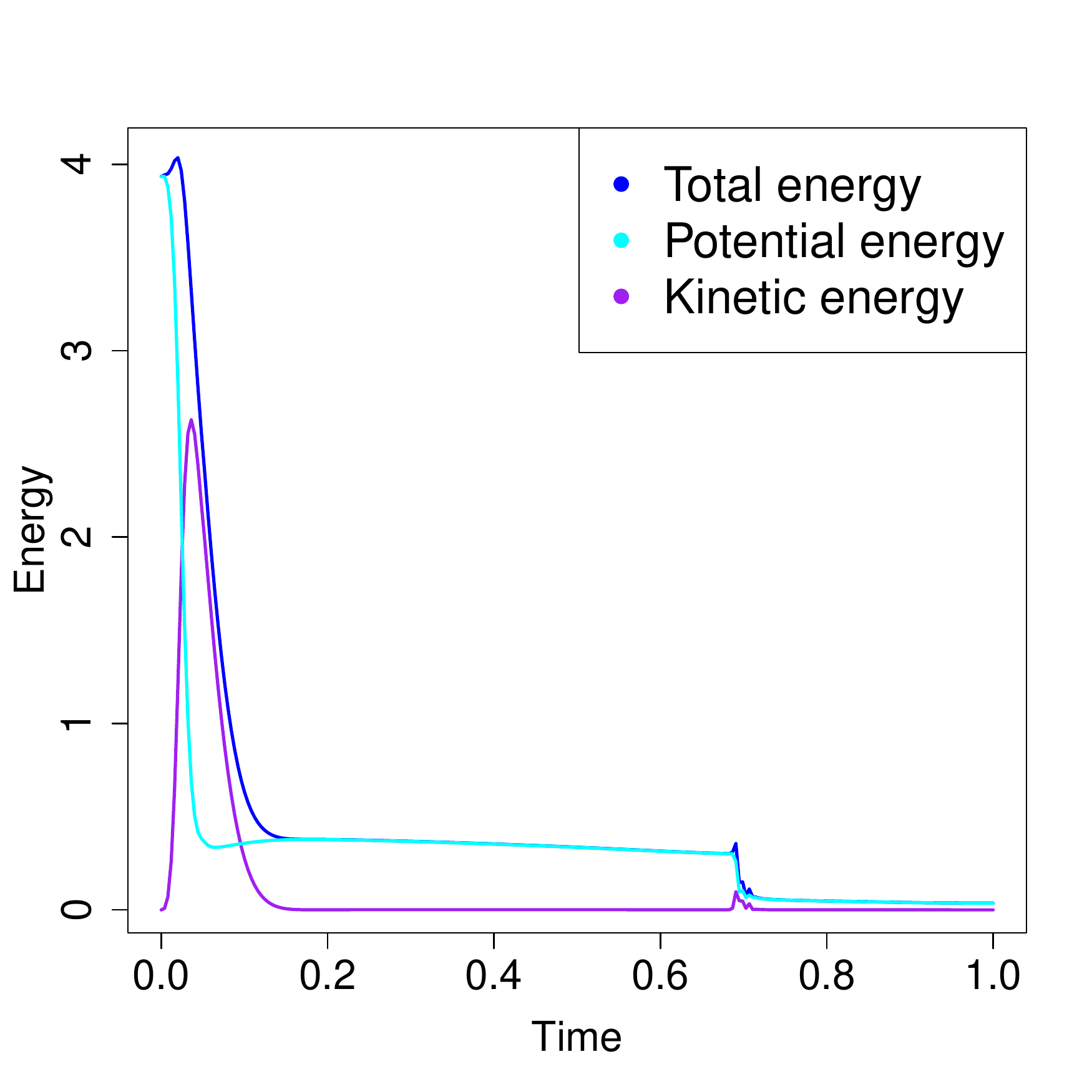}}
		\hspace{\fill}
		\subfloat[Test case 2 \label{fig:TC2_Ehist}
		]{\includegraphics[width=0.45\textwidth]{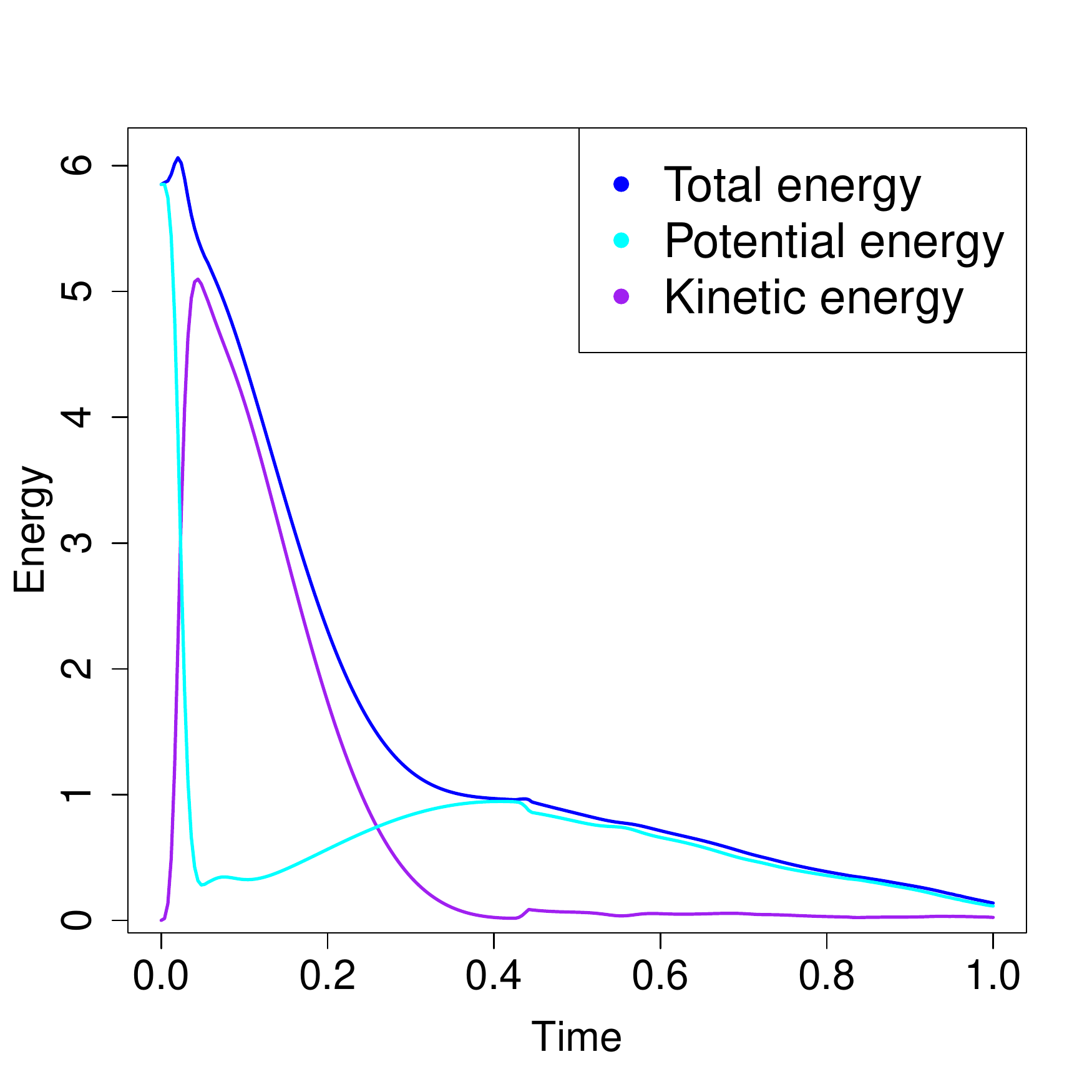}}
 	\end{center}		
	\caption{Test cases 1 and 2 -- histories of the total energy (blue), potential energy (cyan), and the kinetic energy (purple) during the momentum method starting in $q^0$}
\end{figure}

 As in the first test case,  the final solution $q^{\text{GD}}$ of the gradient descent method did not overcome the circular obstacle $\varpi$ and terminated in a solution above $\varpi$, see Figure~\ref{subfig:TC2_solGD}. \kk{It terminated with an objective value of $J_\lambda(q^{\text{GD}})\approx 0.1795$ after} the maximum number of \kk{$200$} iterations was reached. Hence, there is not a guarantee that this solution is a local minimum of $J_\lambda$. 
 \kk{A significantly better result was obtained with the momentum method that terminated in the} solution $q^{\text{HF}}$ with objective value \otd{$J_\lambda(q^{\text{HF}})\approx0.1141$}. It is again located beneath $\varpi$ and resembles the shape of a spoon, see Figure~\ref{subfig:TC2_soldHF}. \kk{We observe that the indicator for the probability of failure for the final spoon shape} is \otd{$J_1(q^{\text{HF}})\approx 0.0413$} and hence a magnitude smaller than that for the whale shape depicted in Figure~\ref{subfig:TC1_PT_circ}, which was $J_1=0.4819$. 

\begin{figure}
    \centering
    \includegraphics[width=0.6\textwidth]{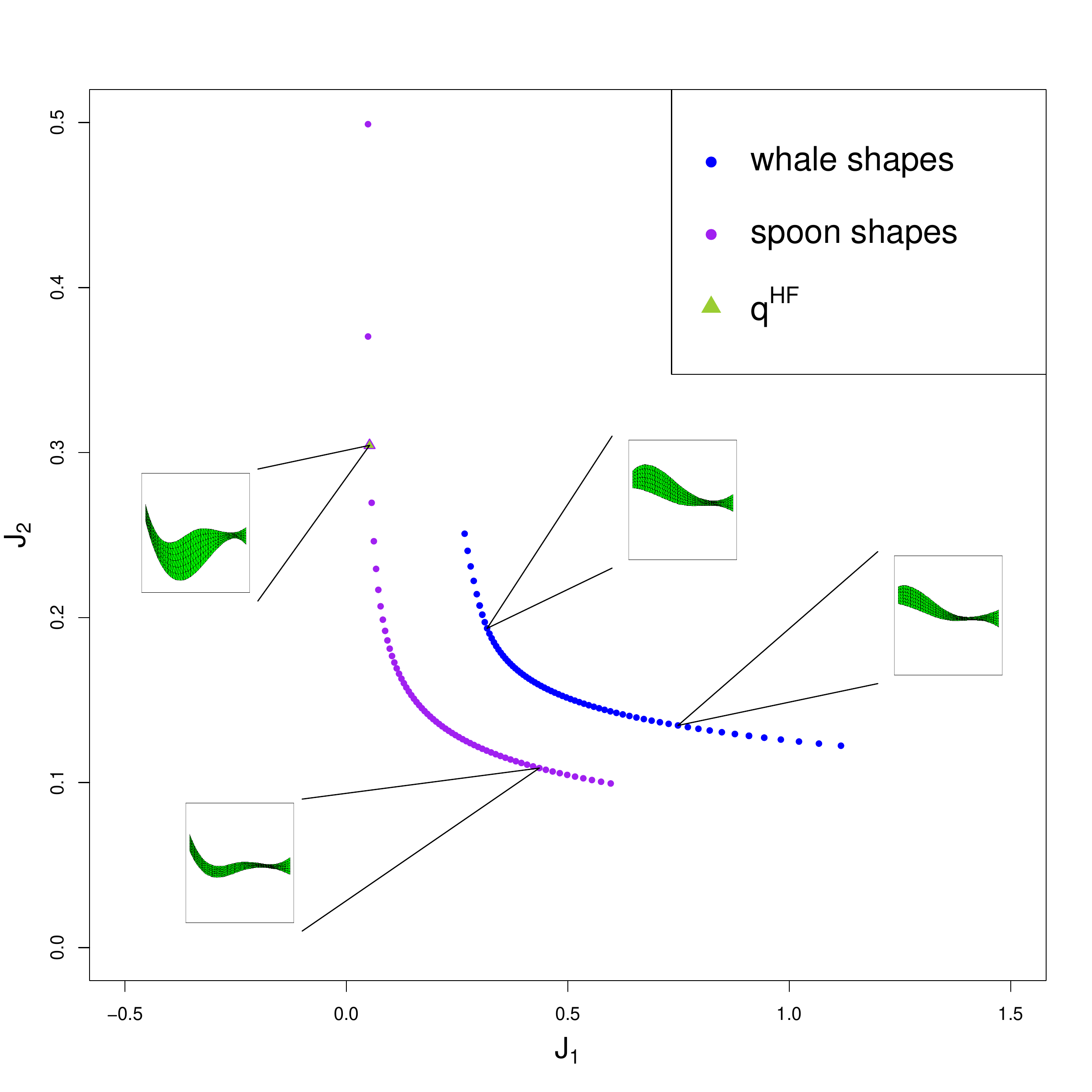}
    \caption{\otd{Test case 2 -- comparison of the (local) biobjective Pareto fronts, i.e., only $J_1$ and $J_2$ objective values are considered, w.r.t.\ whale shapes (blue, see \cite{bolten2021tracing}) and spoon shapes (purple, with initial shape $q^{\text{HF}}$ as a green triangle) computed with the tracing method of \cite{bolten2021tracing}.}\label{fig:PTspoon}}
    
\end{figure}

 \otd{While the whale shape from Figure~\ref{subfig:TC1_PT_circ} has a smaller volume than the spoon shape $q^{\text{HF}}$ (we obtained $J_2(q^{\text{HF}})\approx 0.3251$), 
a comparison of a spoon with volume of $J_2 = 0.2068$ (and probability of failure of $J_1 =0.0771$) computed via the tracing approach of \cite{bolten2021tracing} starting in $q^{\text{HF}}$ as the initial shape with another whale shape computed in \cite{bolten2021tracing} that has a comparable volume of $J_2=0.2073$ reveals that the whale shape still has a considerably larger probability of failure of $J_1=0.3009$. Furthermore, comparing the (local) Pareto fronts w.r.t.\ $J_1$ and $J_2$ of whale shapes computed in \cite{bolten2021tracing} and spoon shapes yield that the spoon shapes dominate the whale shapes, see Figure~\ref{fig:PTspoon}.}

Thus, in comparison with the gradient descent method employed in \cite{Doganay2019} the momentum method yields preferable solutions also for the unpenalized, biobjective shape optimization problem of \cite{bolten2021tracing,Doganay2019}.

 \otd{Figure~\ref{fig:TC2_Ehist} shows the histories of the total energy, potential energy and the kinetic energy during the discretized Hamiltonian flow starting in $q^0$. We start with a potential energy of $5.8513$ which decreases during the first steps of the optimization to $0.2862$ and then slightly increases again in the next iterations. On the other hand, the kinetic energy reaches a value of $5.0986$ after 12 iterations 
 and then decreases in the next iterations while the potential energy increases again. Except for a small increase during the initial iterations due to the limited precision of the symplectic Euler method, the total energy decreases during the approach.}

\begin{figure}[]
	\begin{center}
		\subfloat{\includegraphics[width=0.2\textwidth]{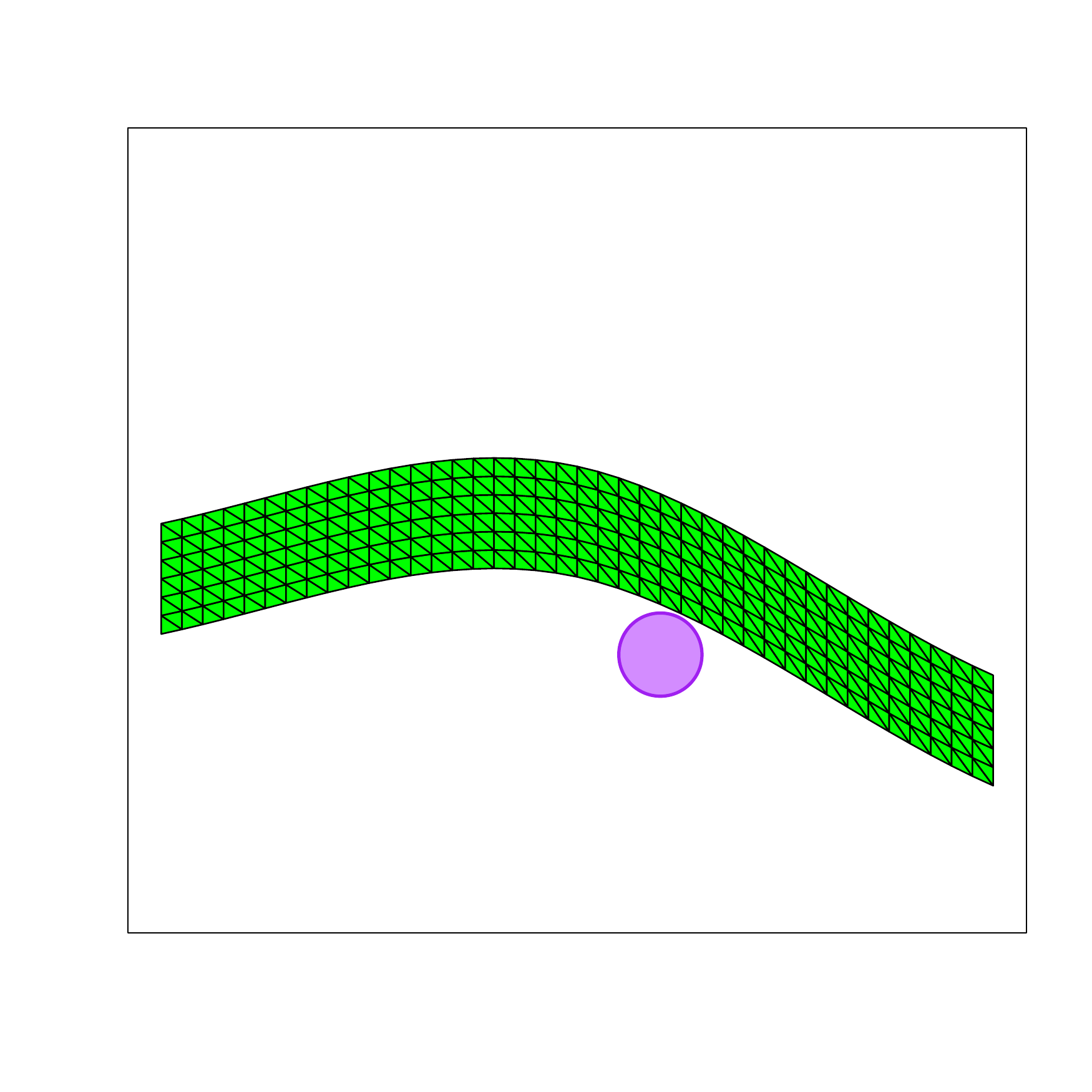}}
		\subfloat{\includegraphics[width=0.2\textwidth]{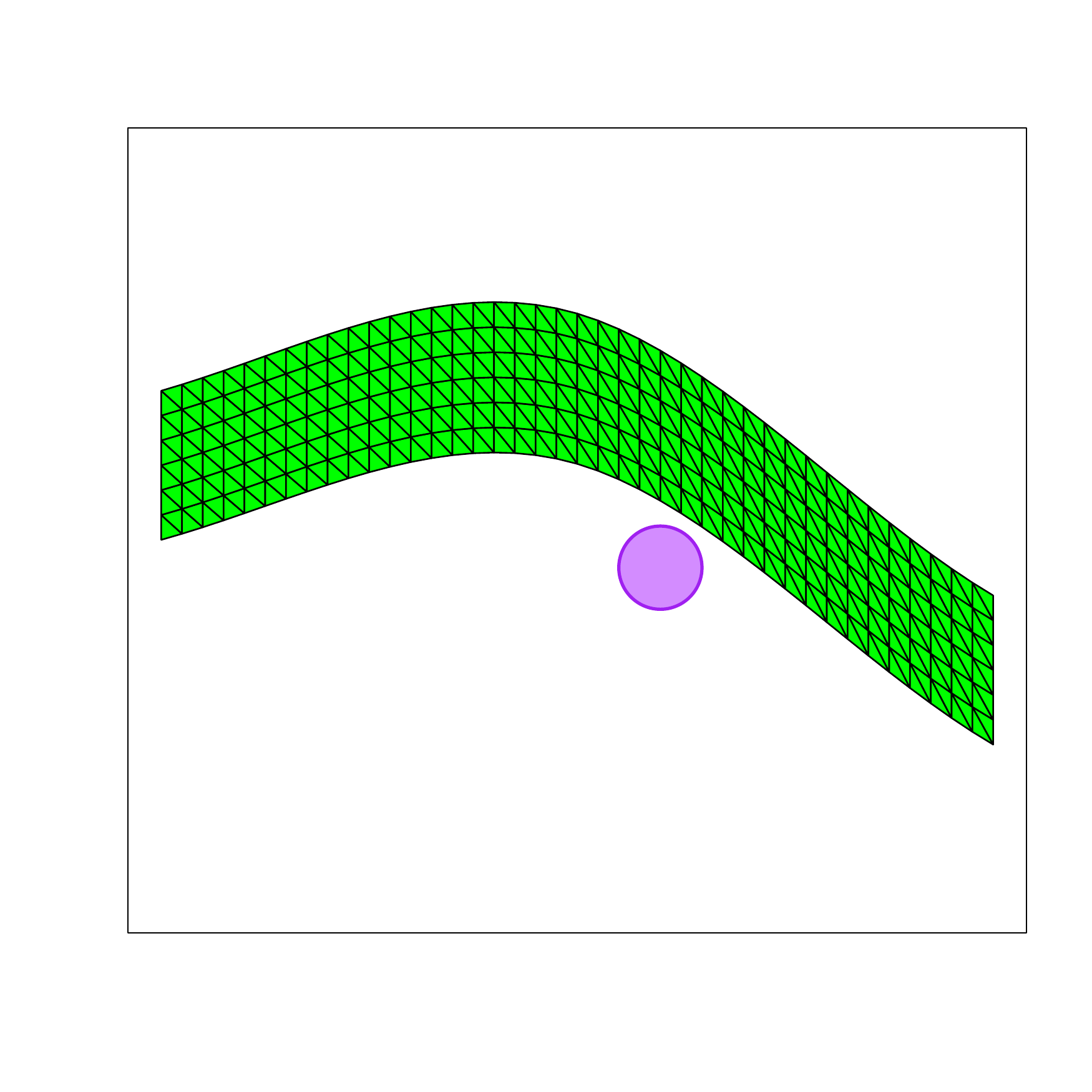}}
		\subfloat{\includegraphics[width=0.2\textwidth]{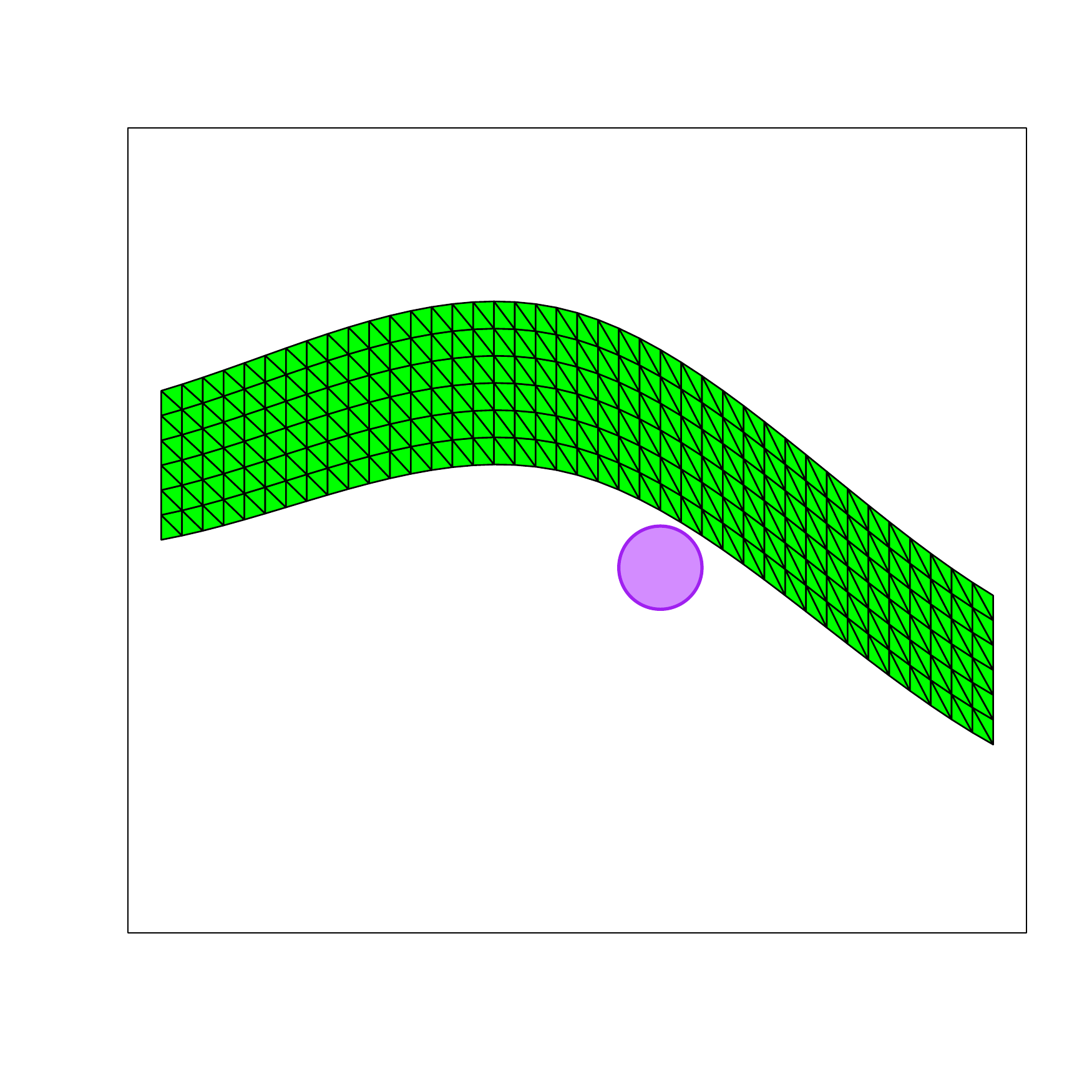}}
		\subfloat{\includegraphics[width=0.2\textwidth]{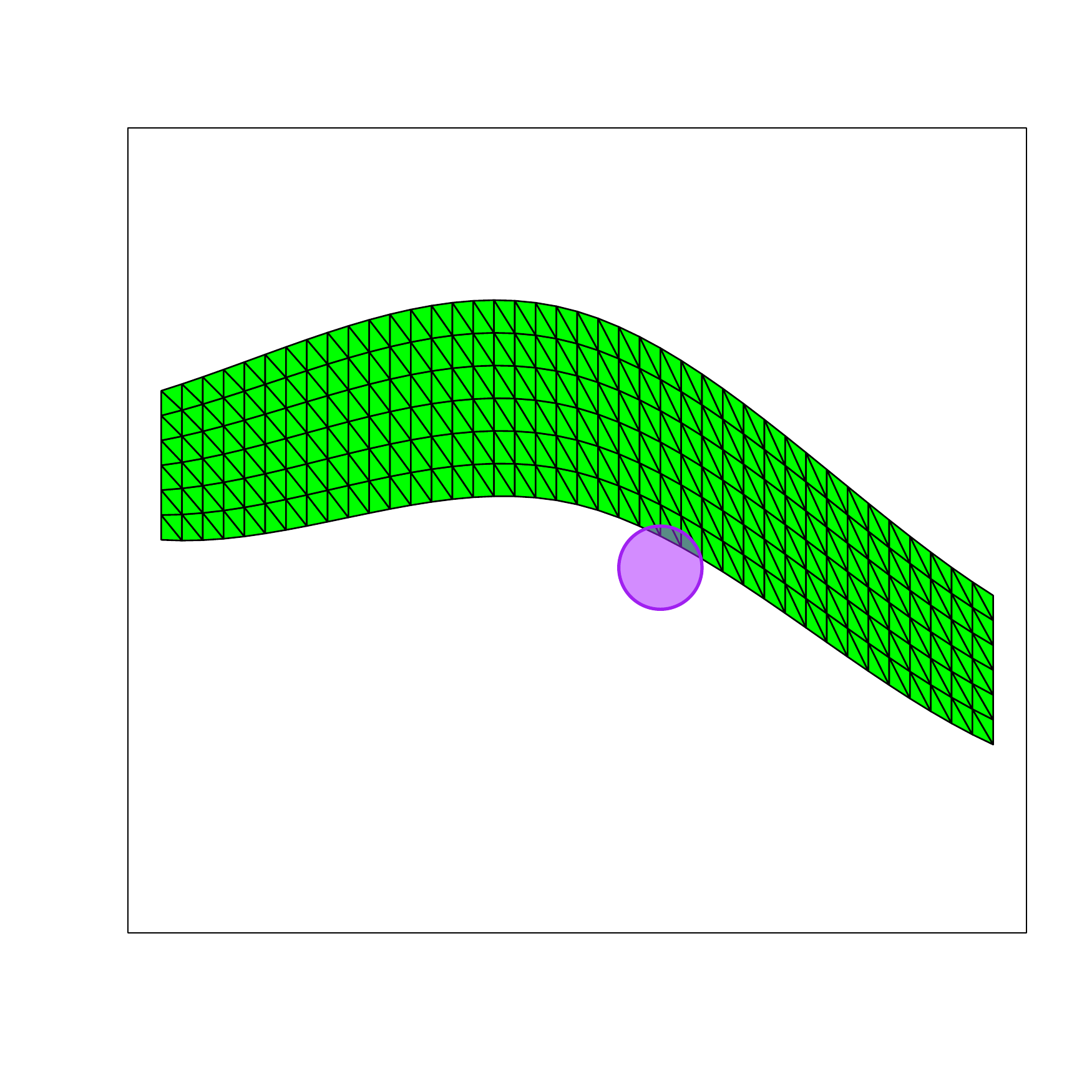}}
		\subfloat{\includegraphics[width=0.2\textwidth]{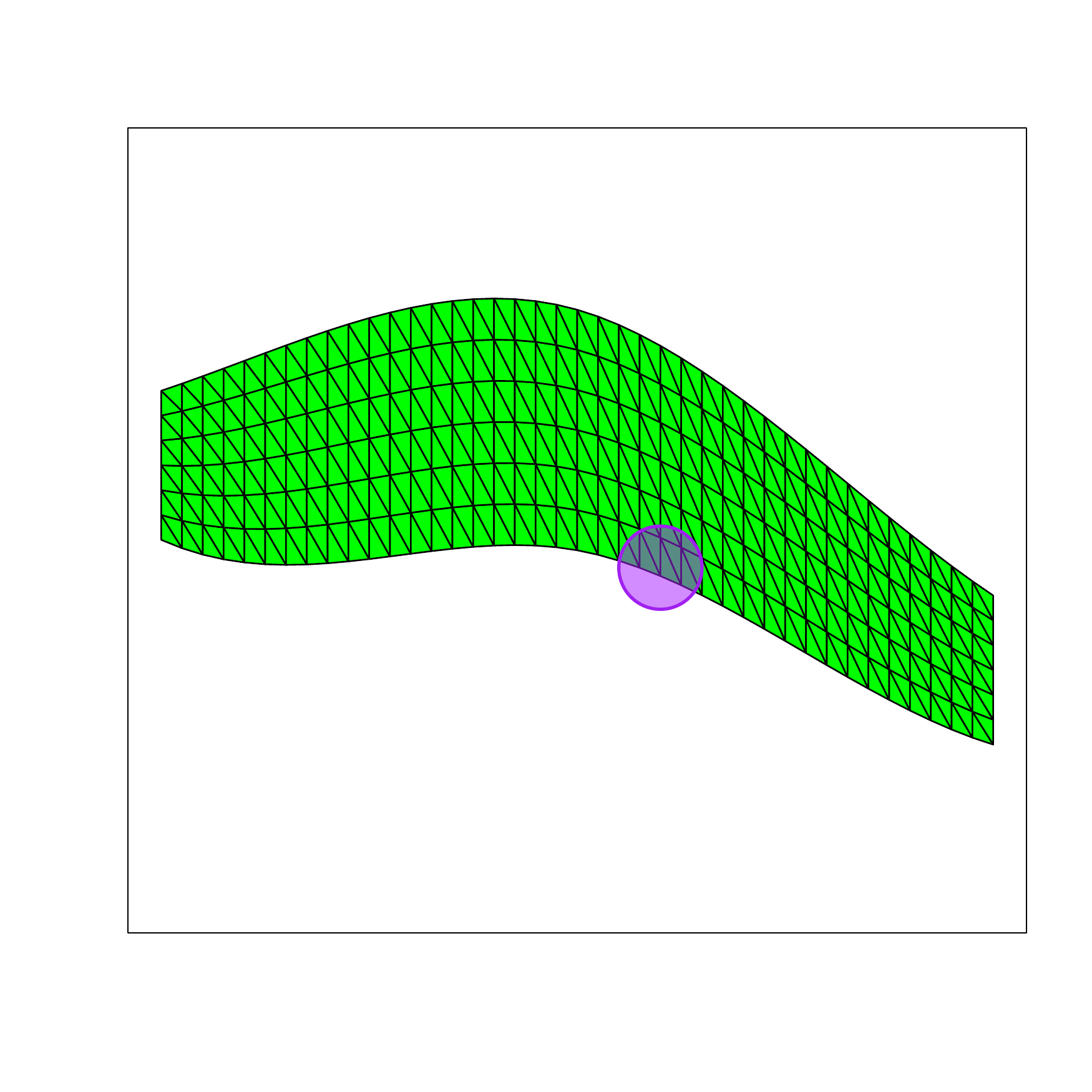}}

		\subfloat{\includegraphics[width=0.2\textwidth]{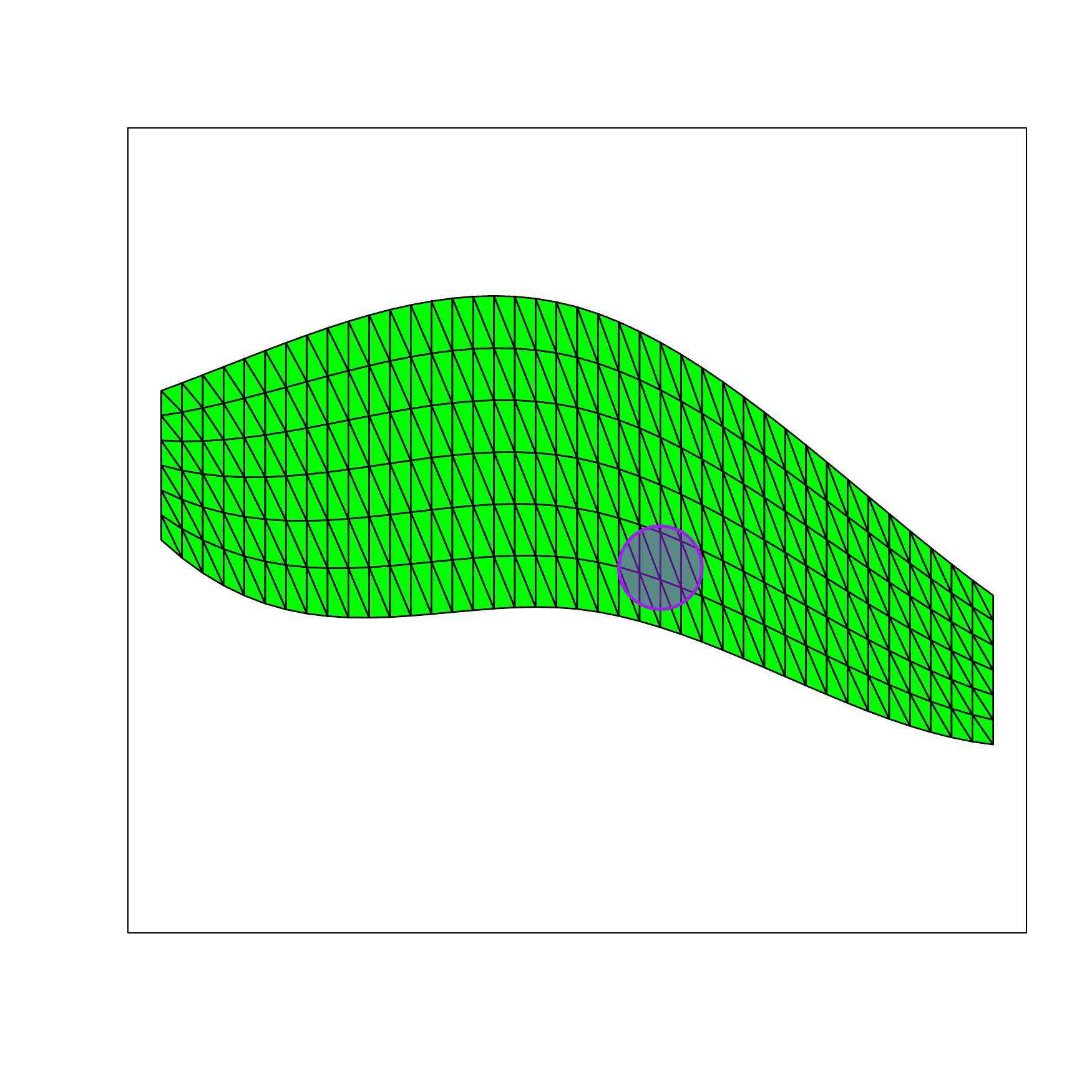}}
		\subfloat{\includegraphics[width=0.2\textwidth]{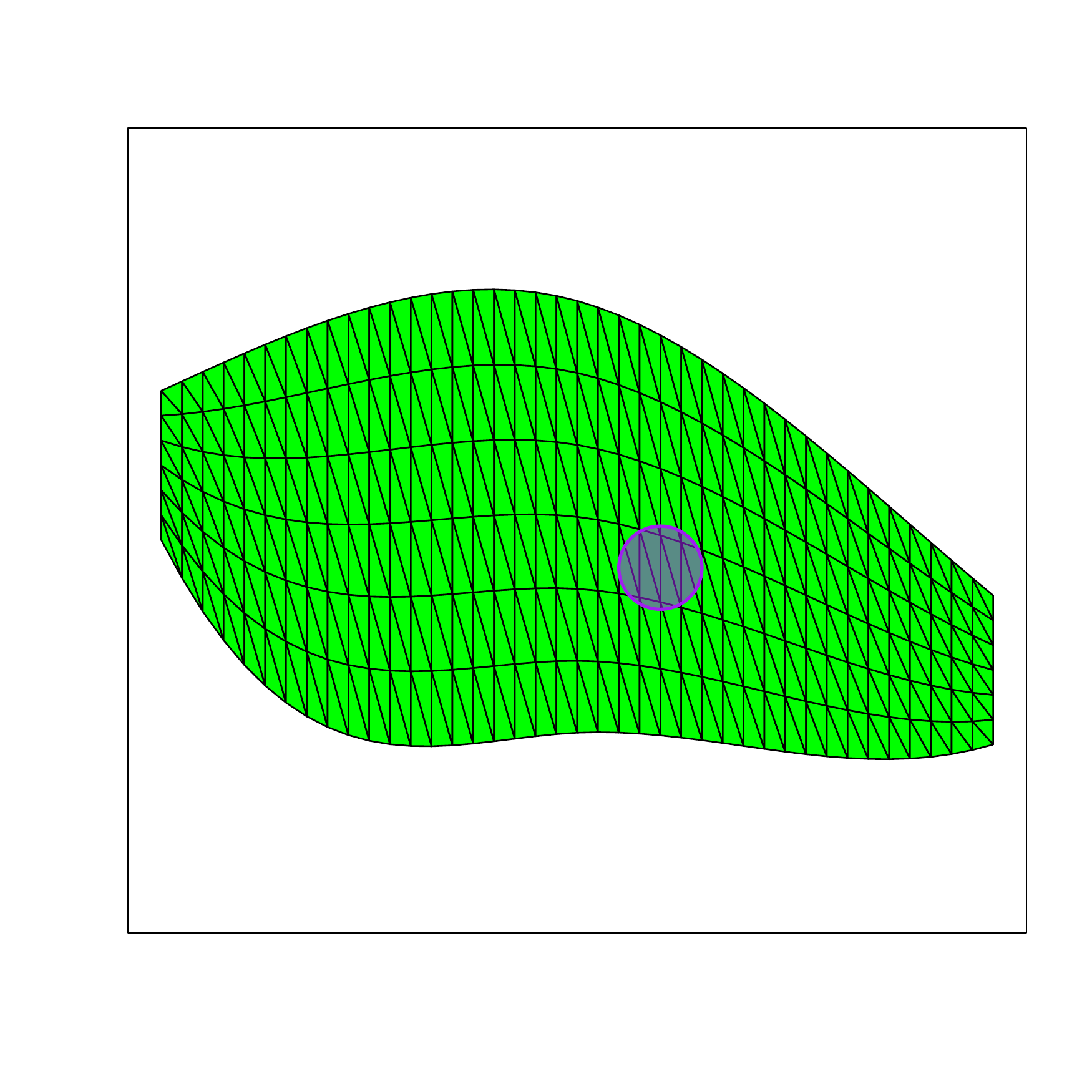}}
		\subfloat{\includegraphics[width=0.2\textwidth]{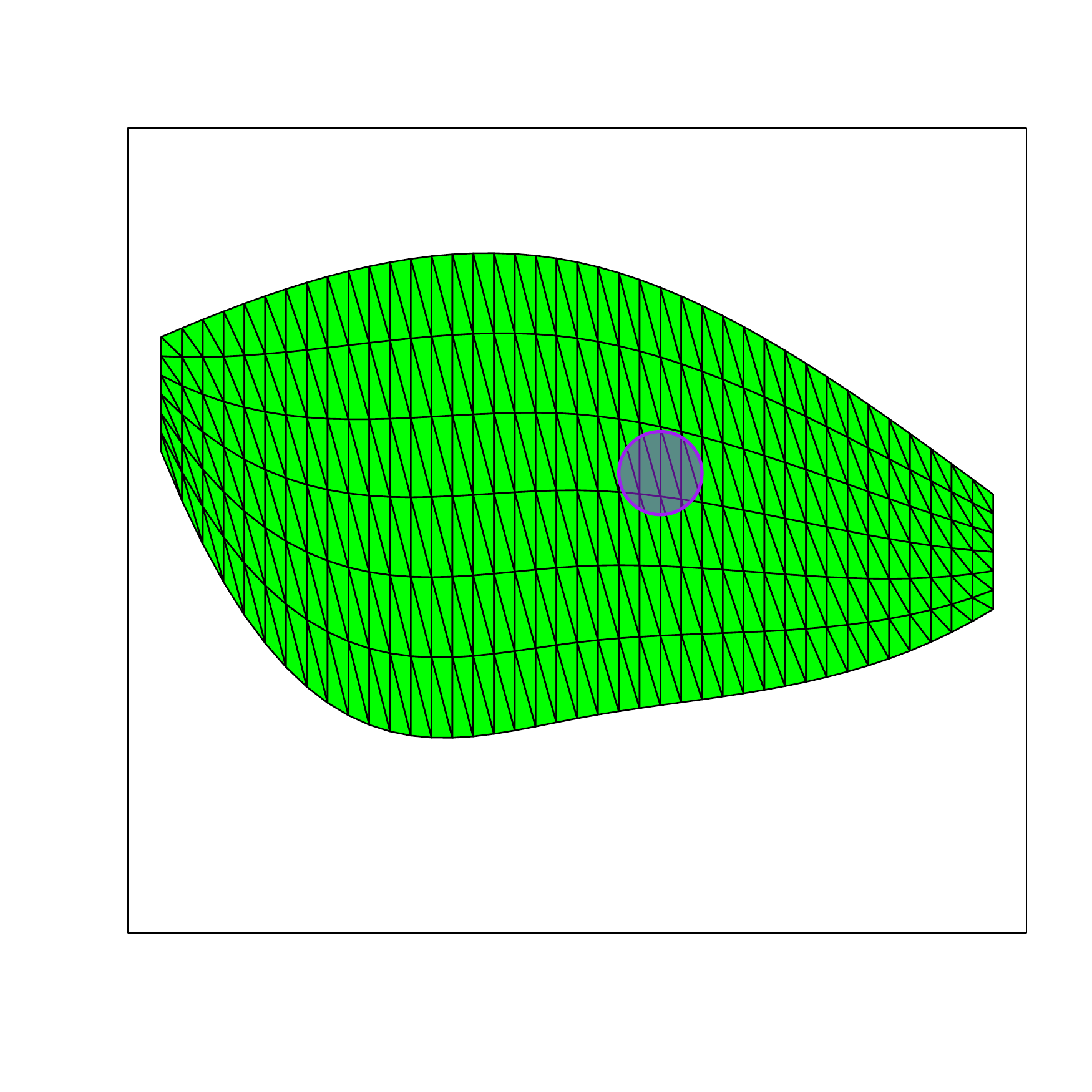}}
		\subfloat{\includegraphics[width=0.2\textwidth]{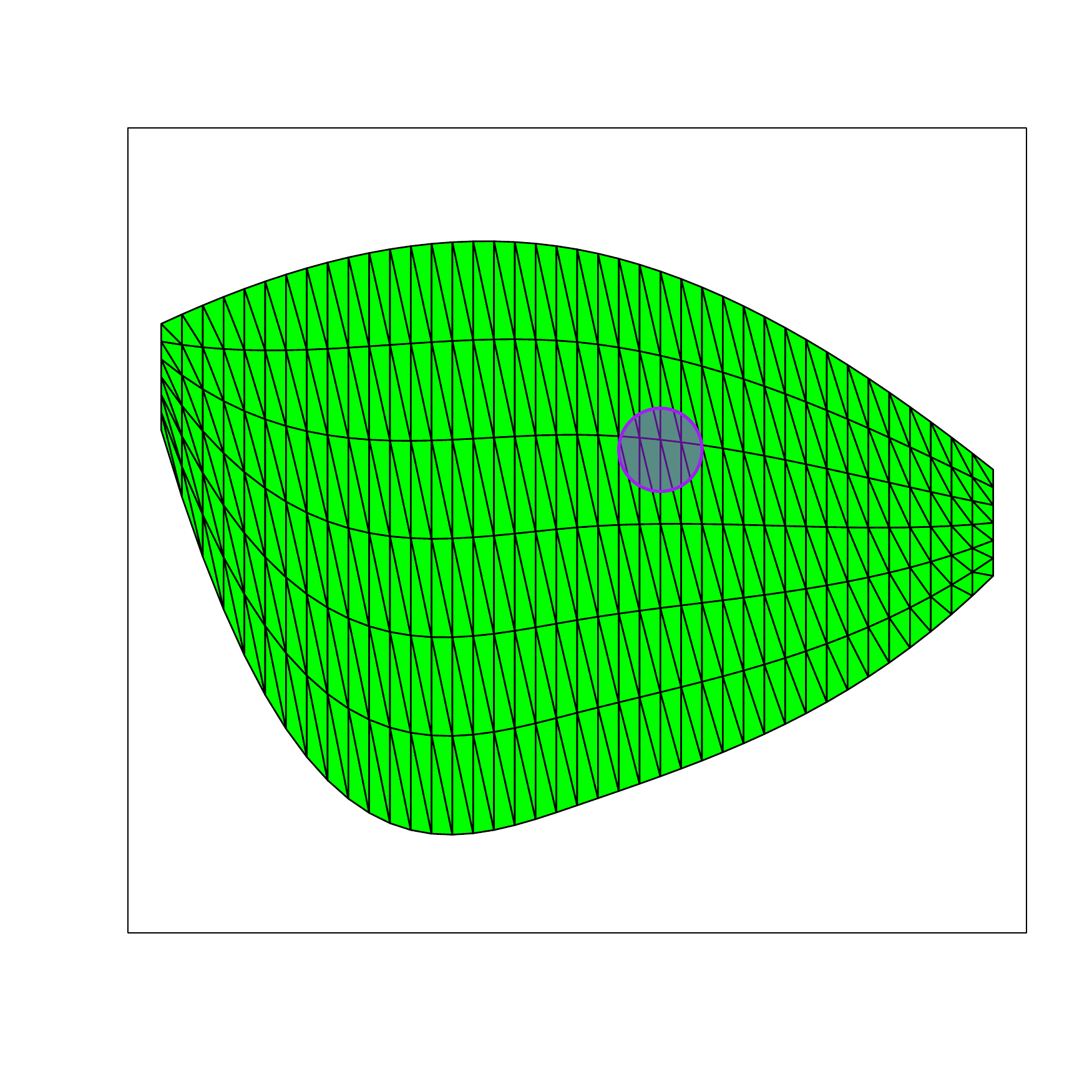}}
		\subfloat{\includegraphics[width=0.2\textwidth]{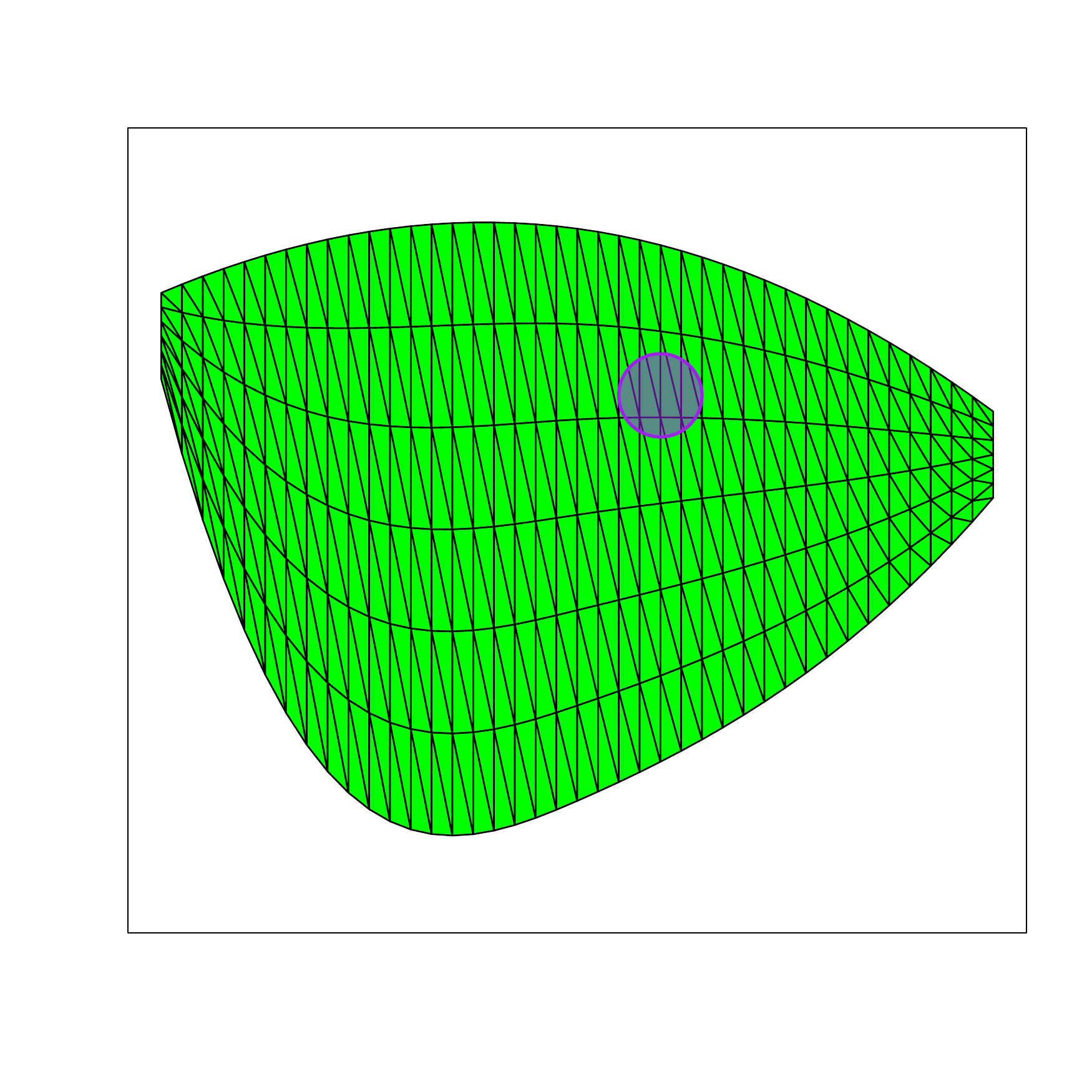}}

		\subfloat{\includegraphics[width=0.2\textwidth]{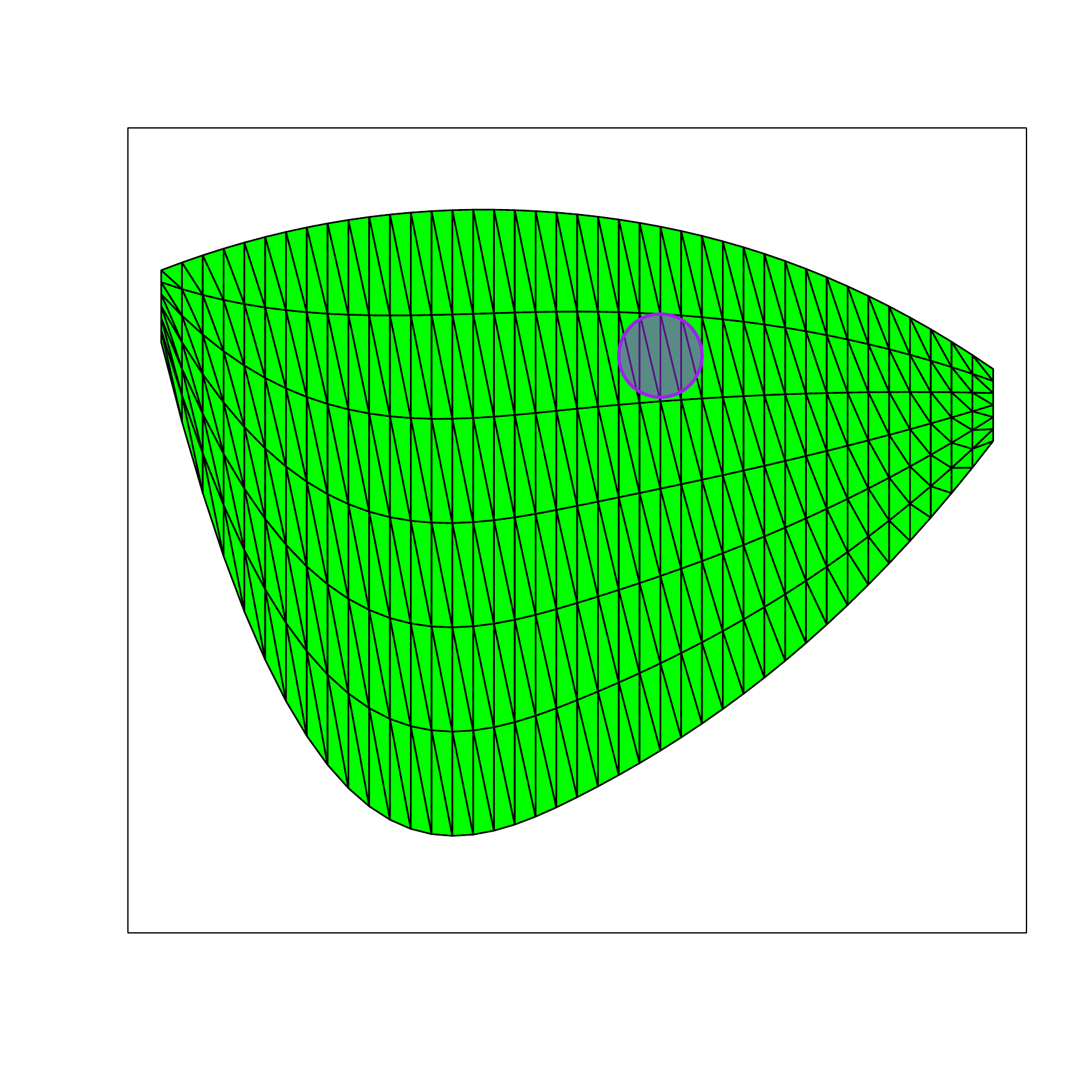}}
		\subfloat{\includegraphics[width=0.2\textwidth]{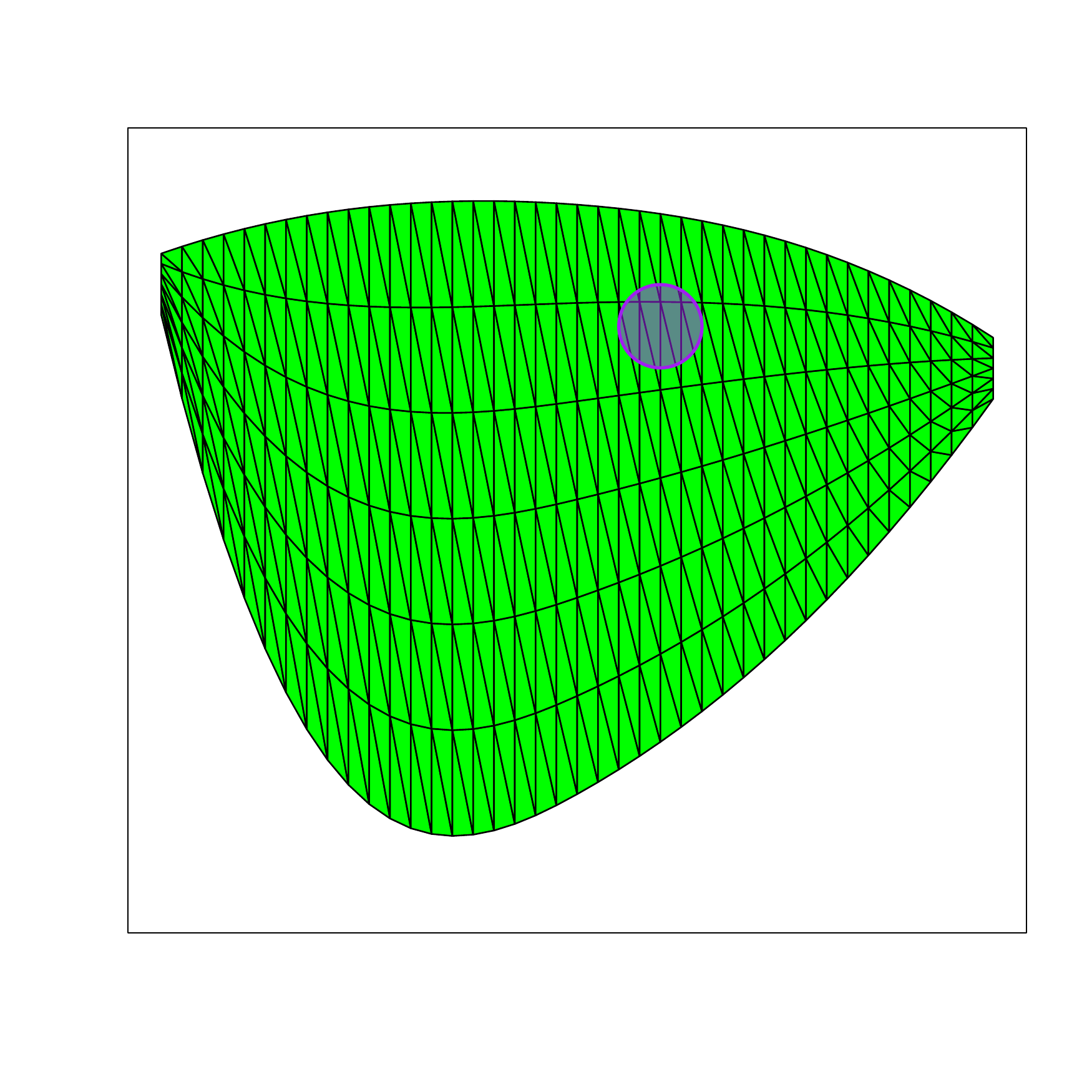}}
		\subfloat{\includegraphics[width=0.2\textwidth]{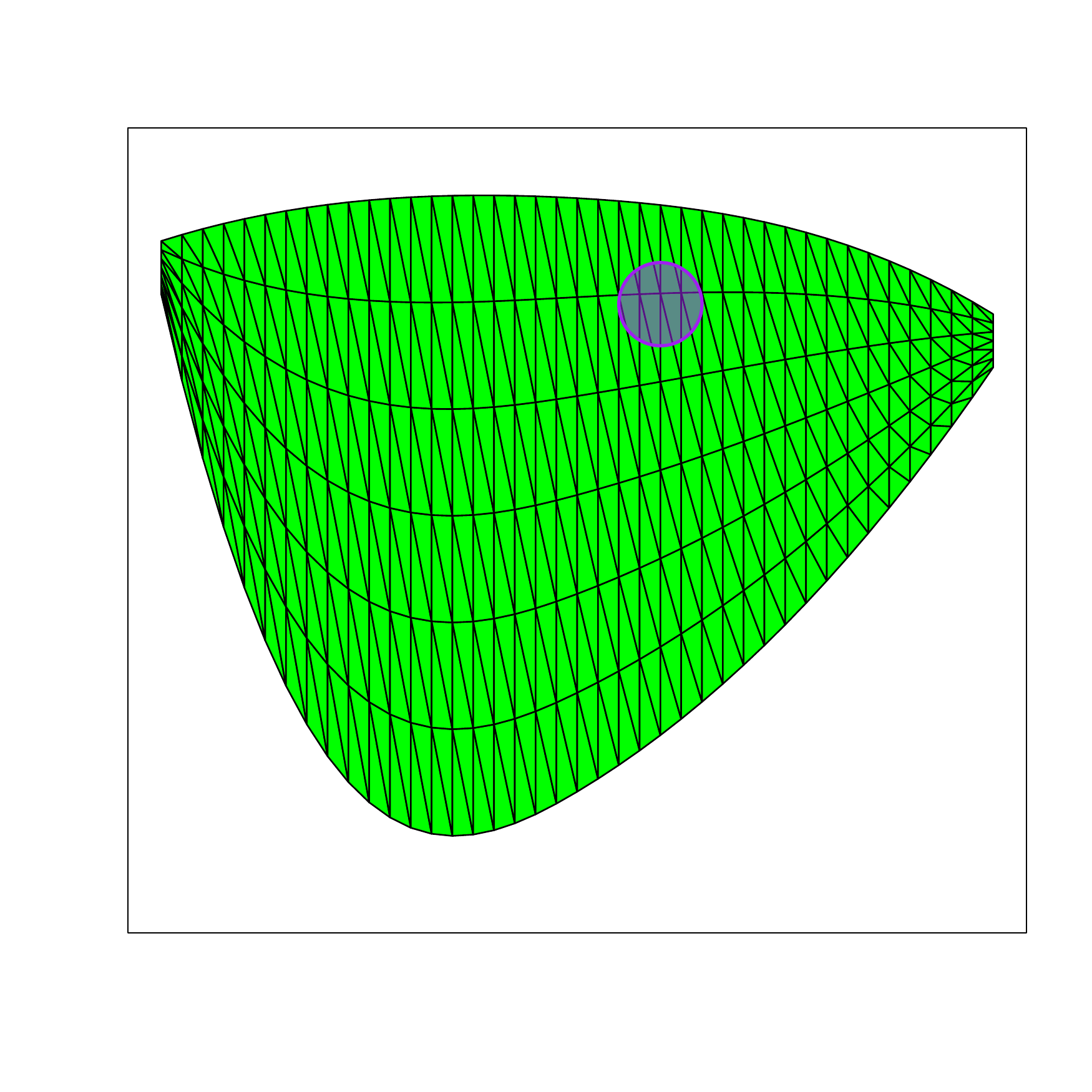}}
		\subfloat{\includegraphics[width=0.2\textwidth]{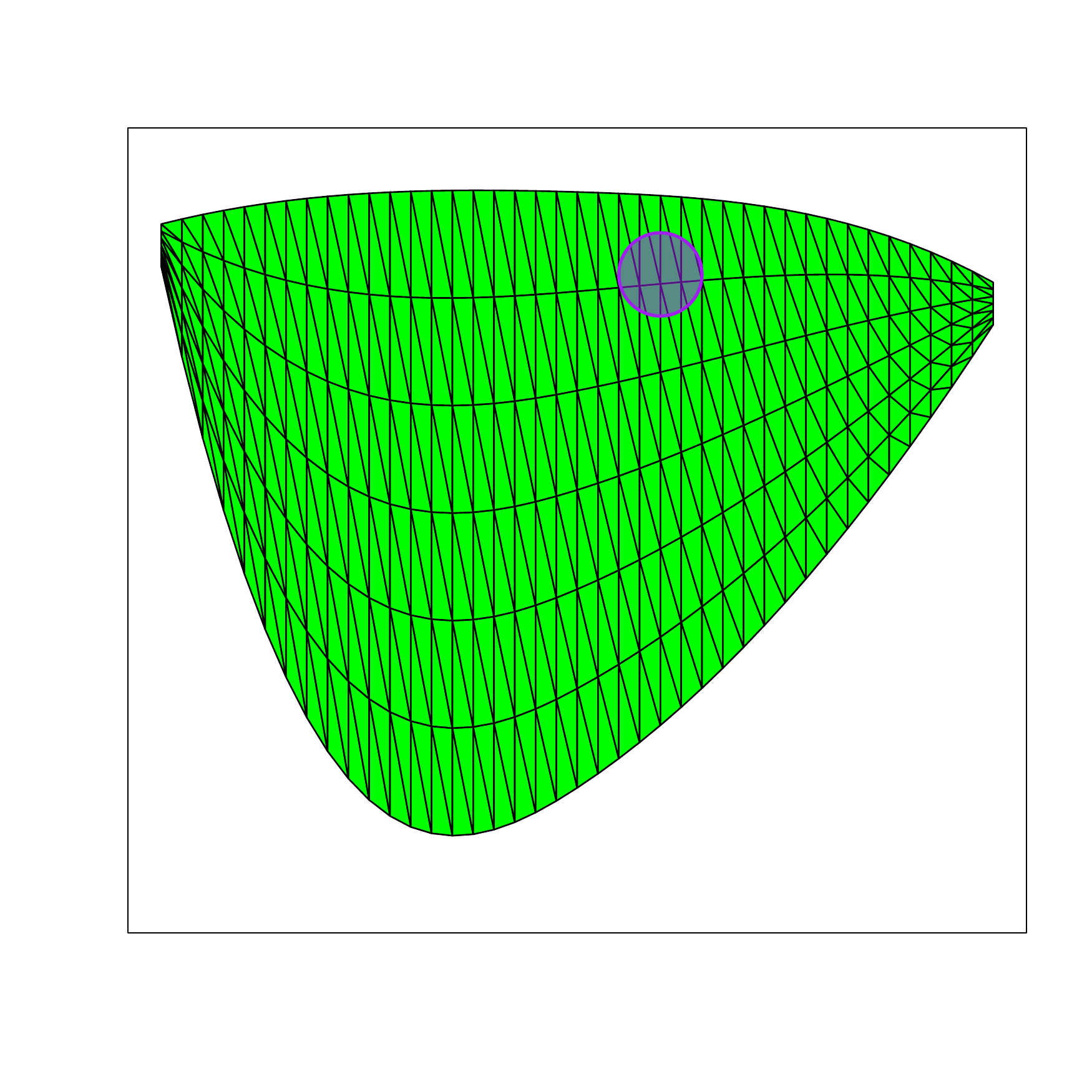}}
		\subfloat{\includegraphics[width=0.2\textwidth]{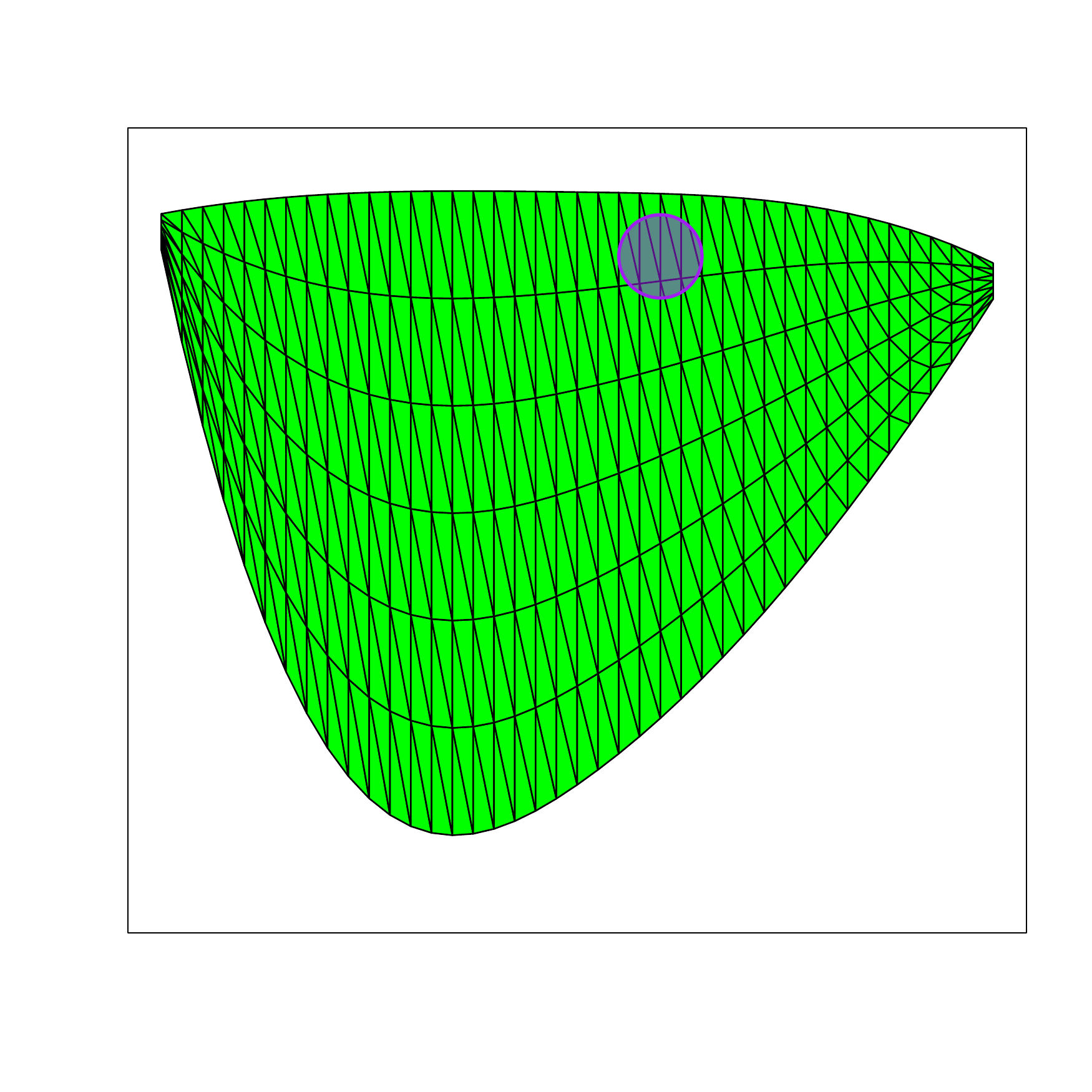}}

		\subfloat{\includegraphics[width=0.2\textwidth]{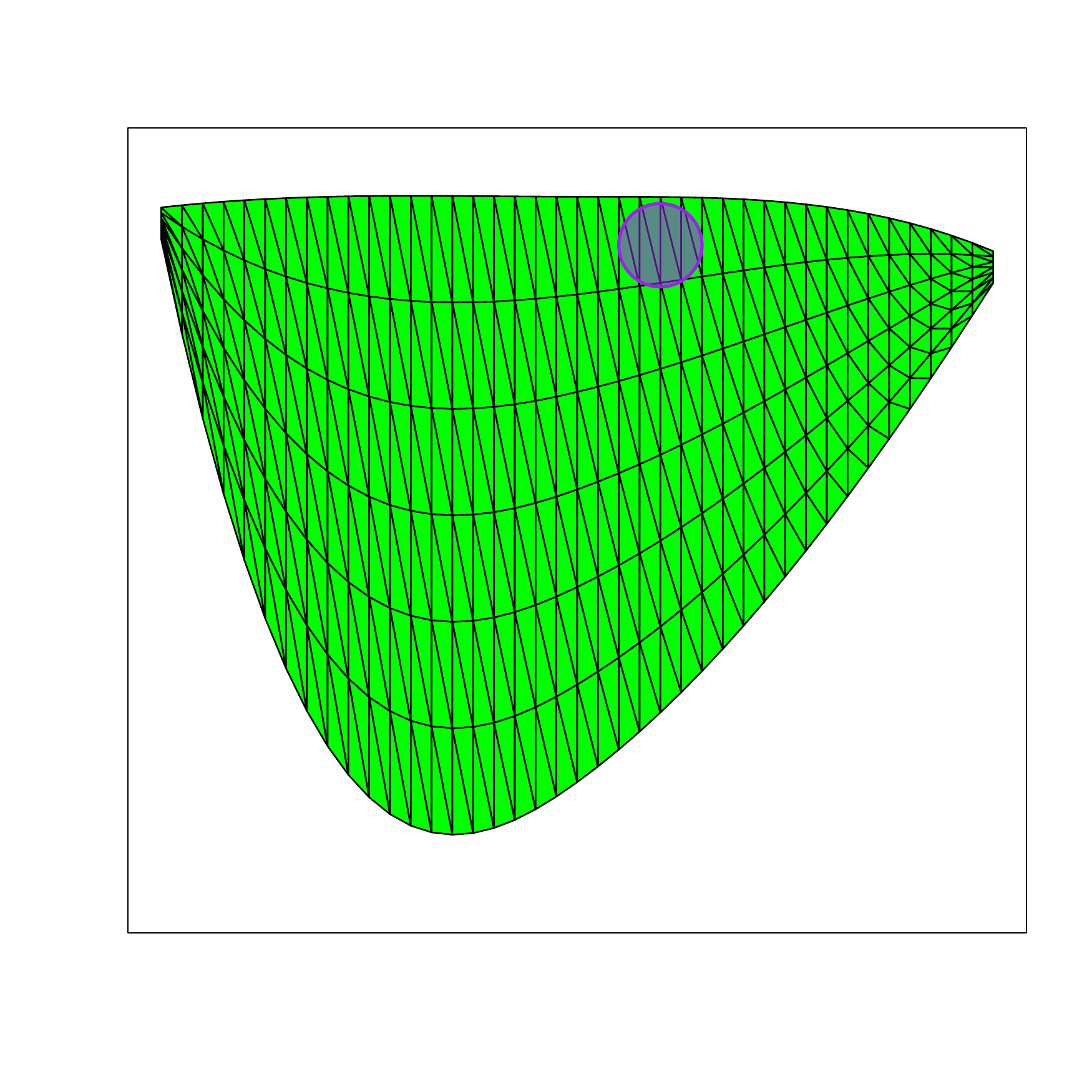}}
		\subfloat{\includegraphics[width=0.2\textwidth]{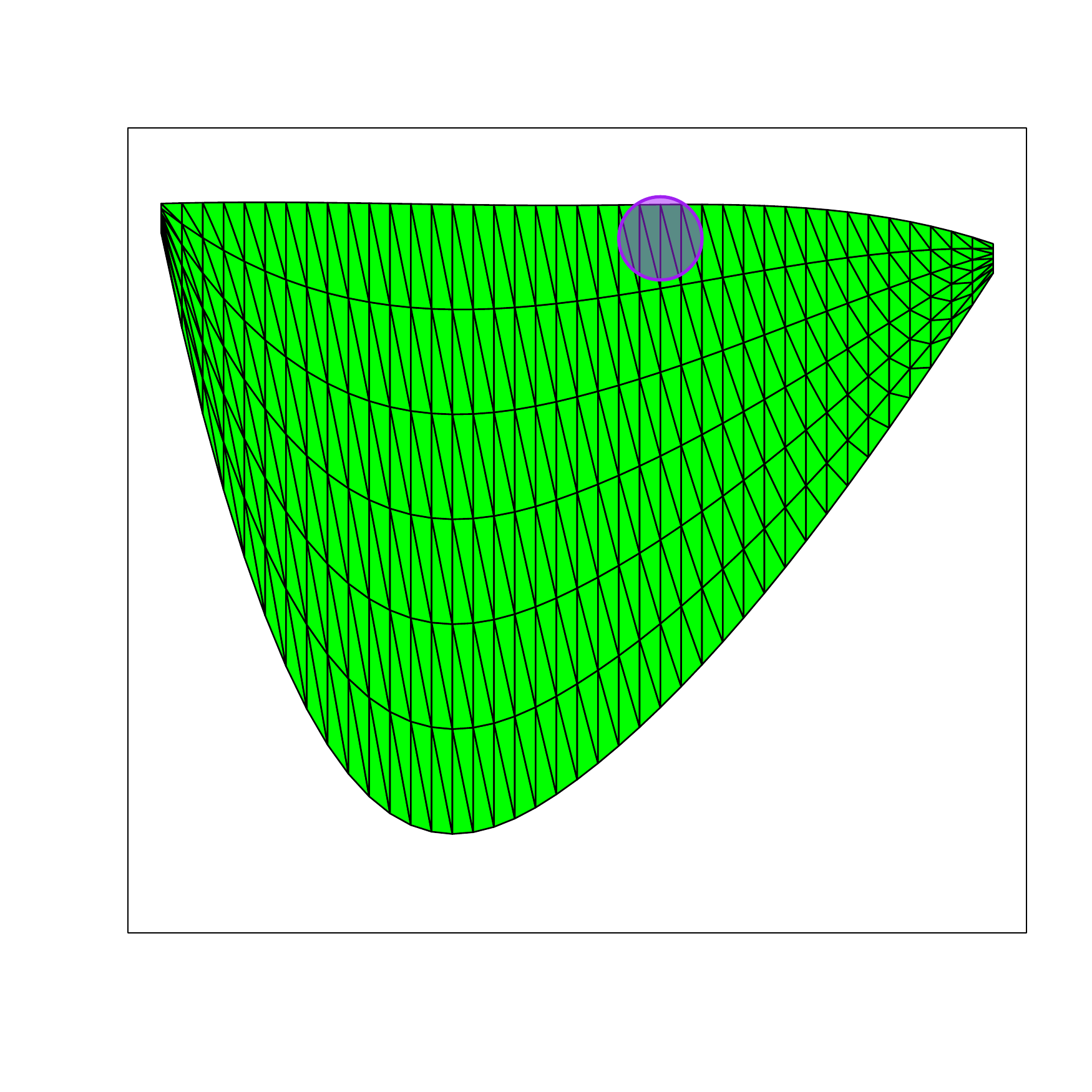}}
		\subfloat{\includegraphics[width=0.2\textwidth]{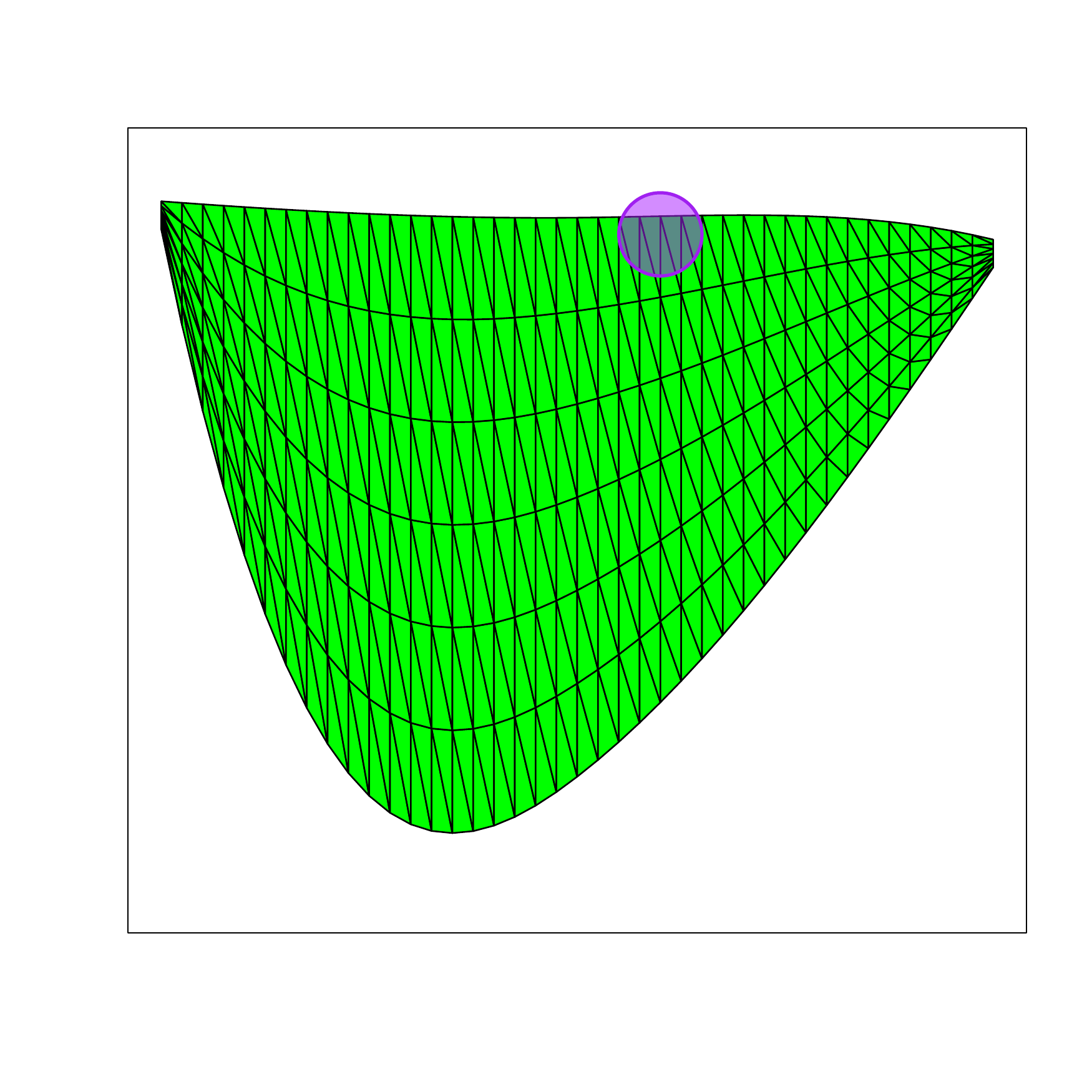}}
		\subfloat{\includegraphics[width=0.2\textwidth]{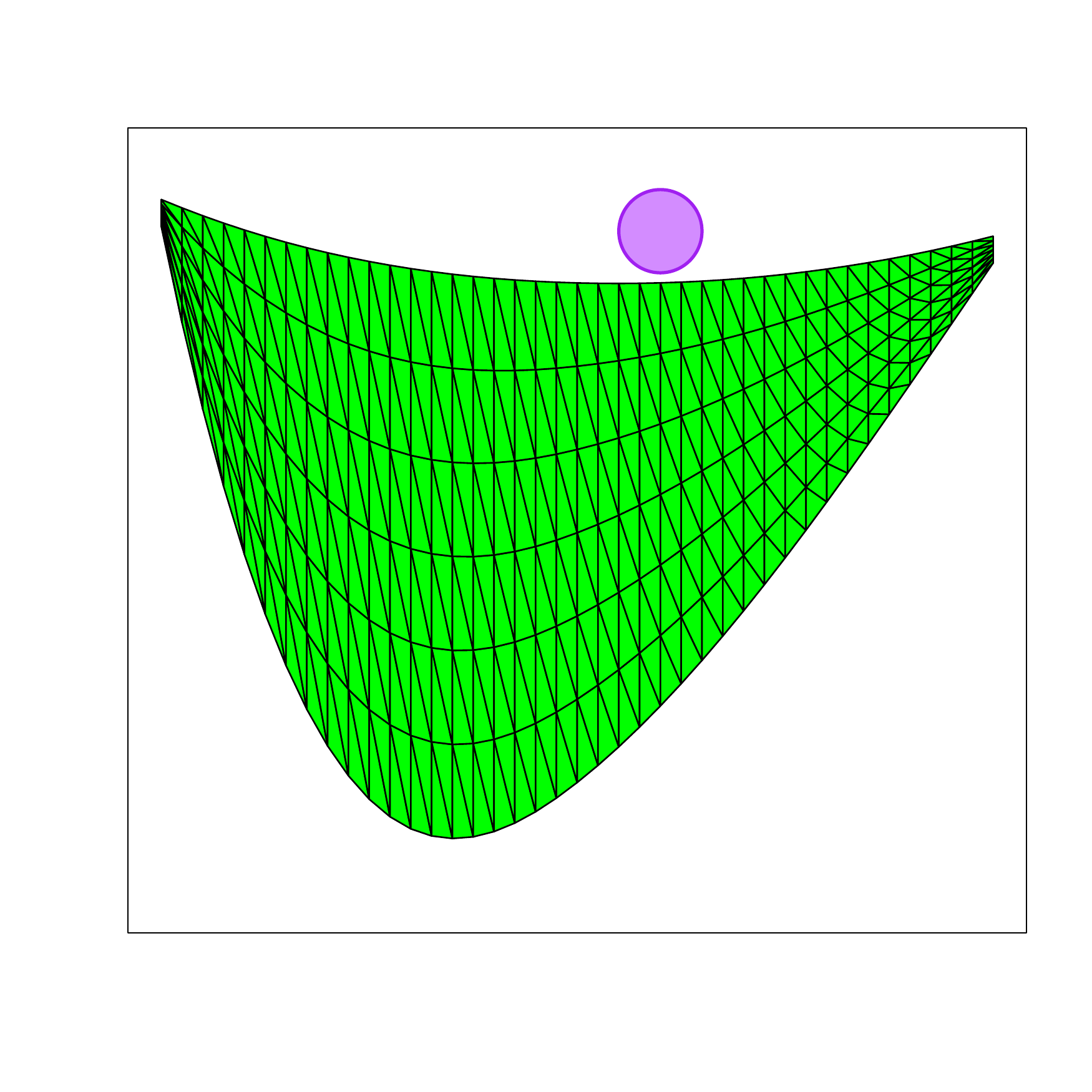}}
		\subfloat{\includegraphics[width=0.2\textwidth]{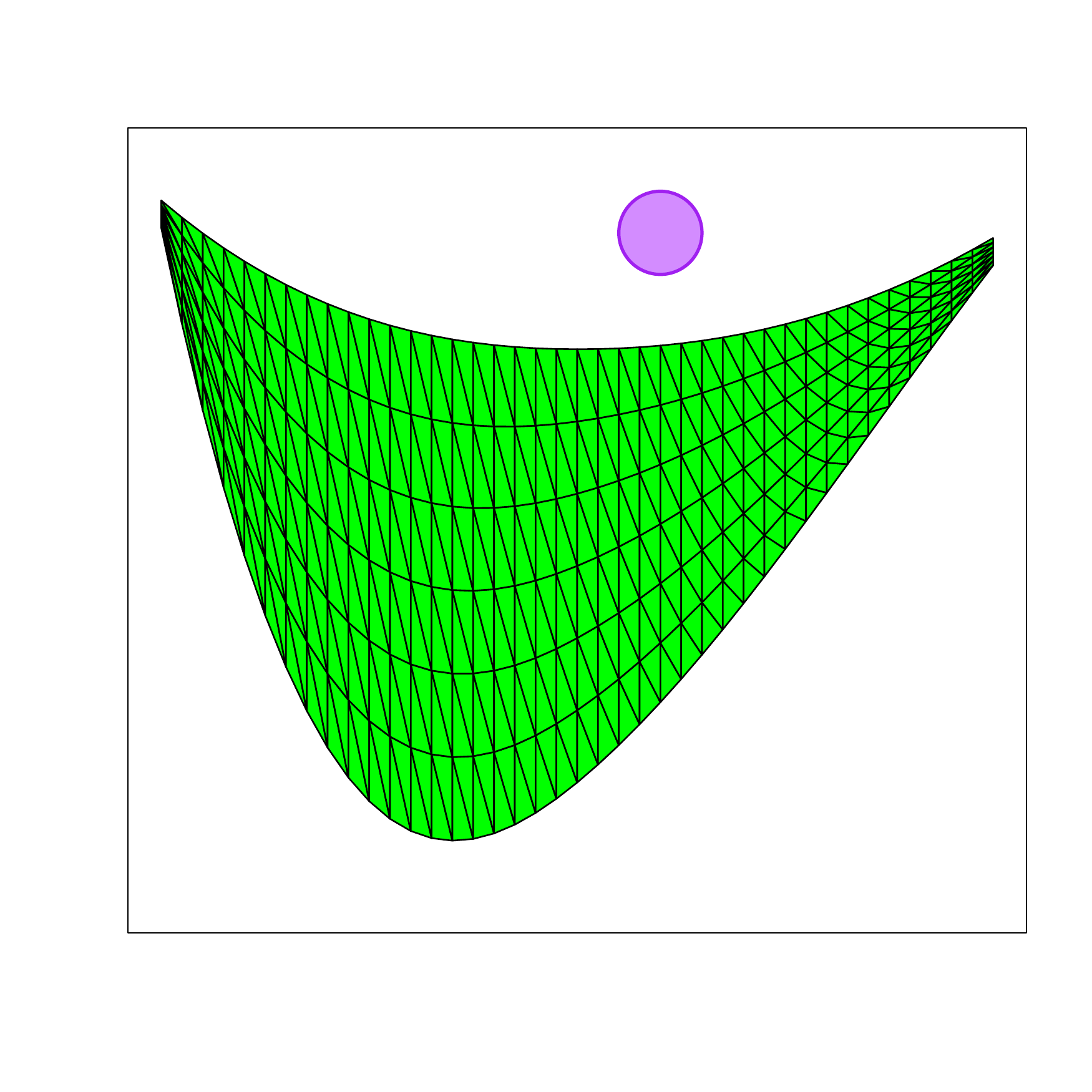}}

		\subfloat{\includegraphics[width=0.2\textwidth]{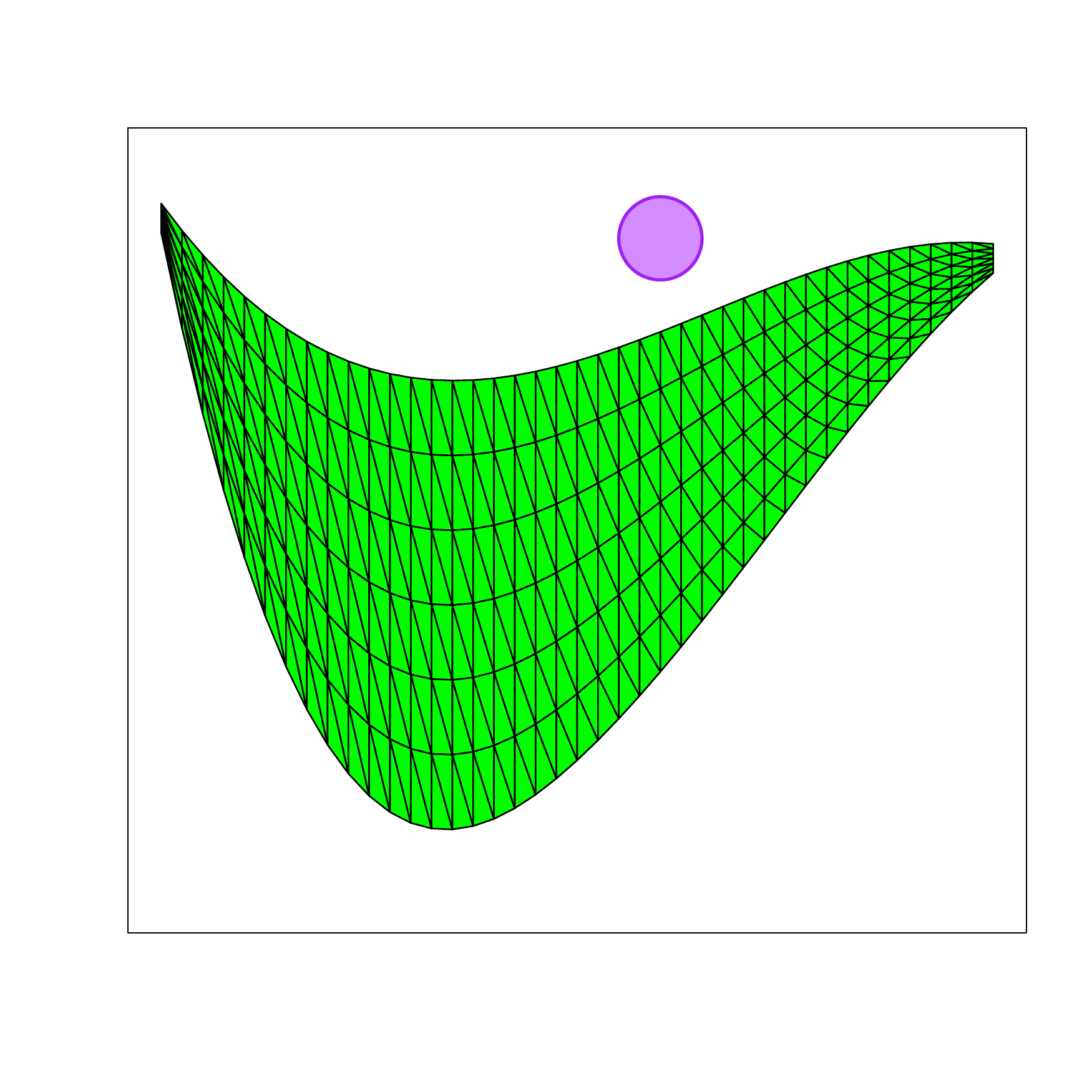}}
		\subfloat{\includegraphics[width=0.2\textwidth]{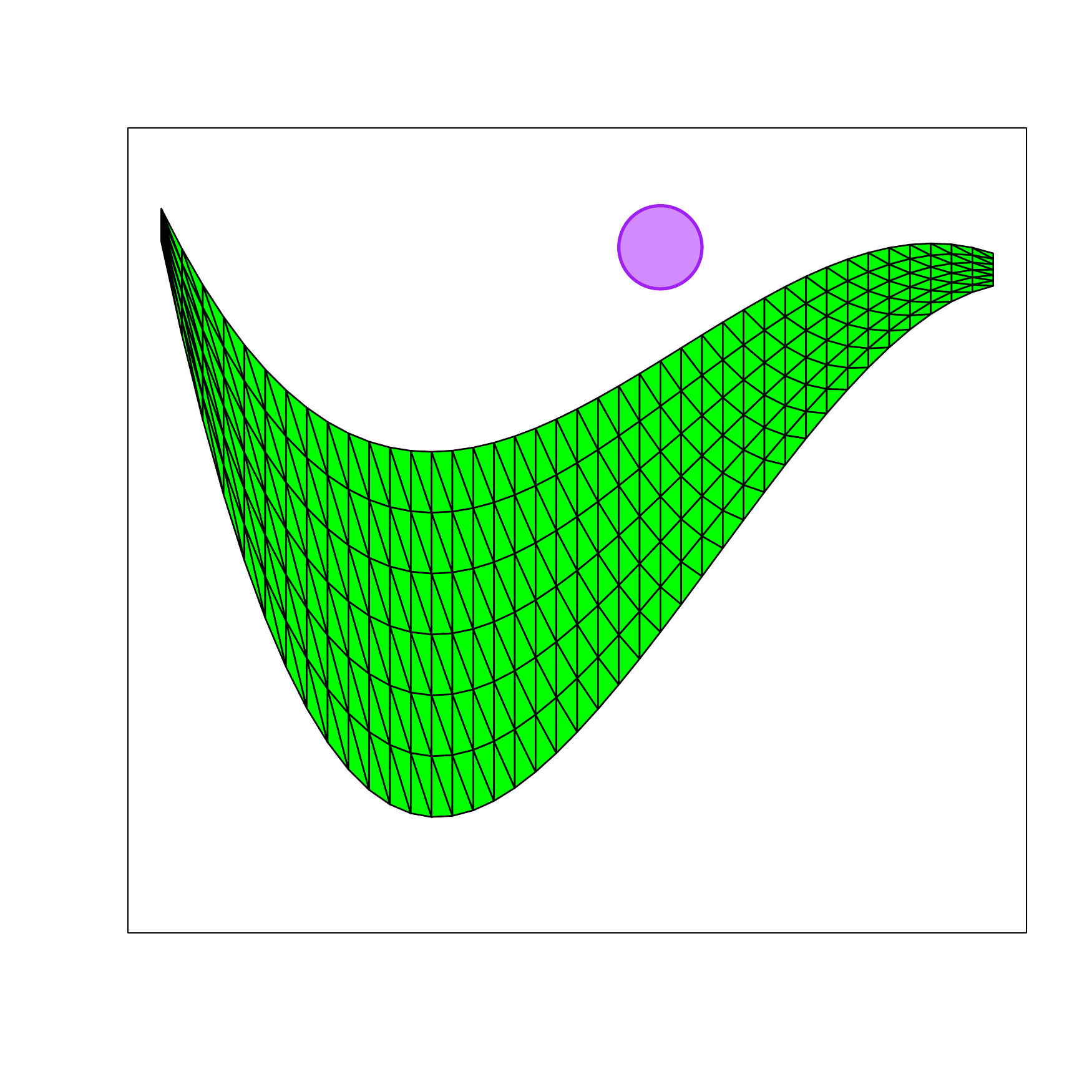}}
		\subfloat{\includegraphics[width=0.2\textwidth]{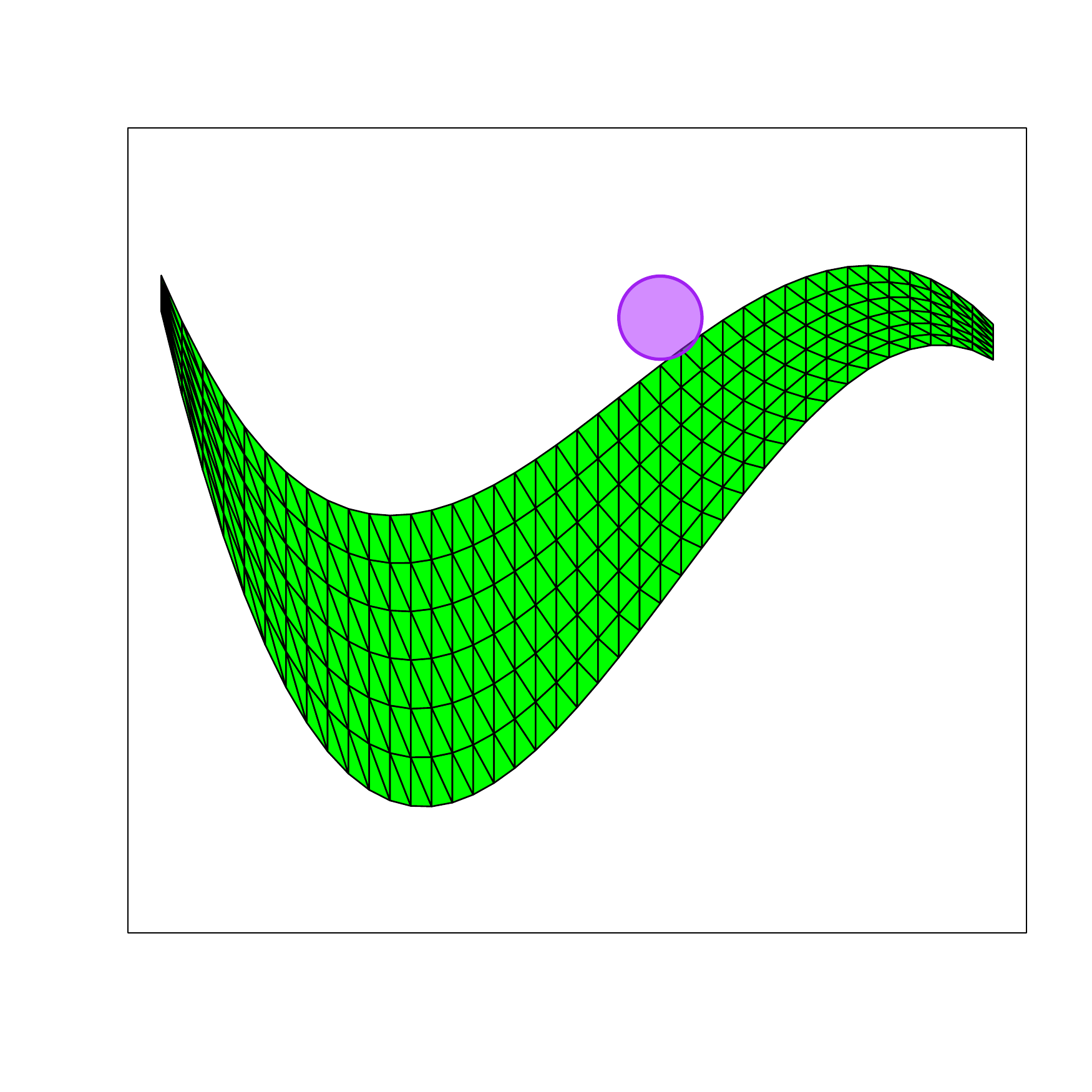}}
		\subfloat{\includegraphics[width=0.2\textwidth]{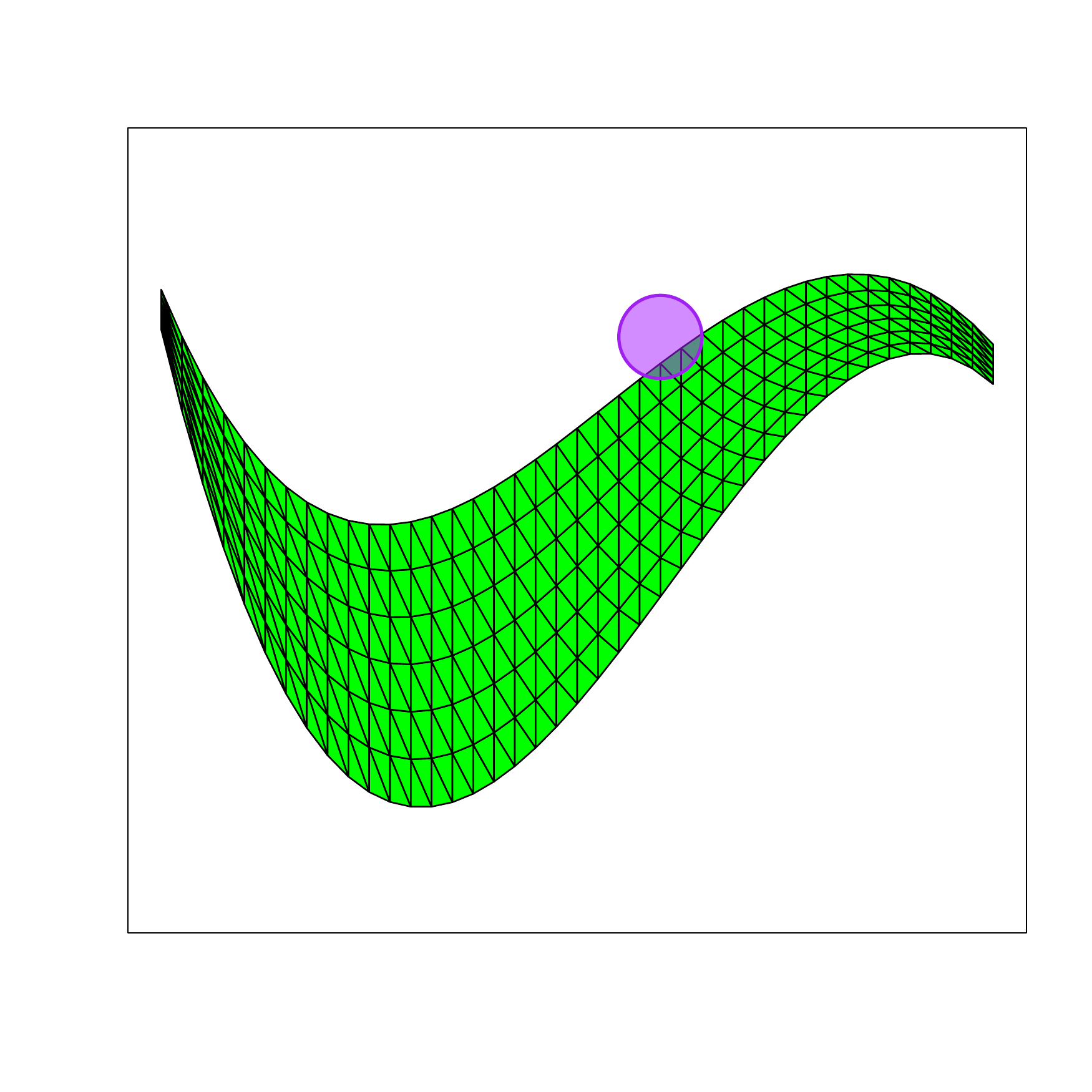}}
		\subfloat{\includegraphics[width=0.2\textwidth]{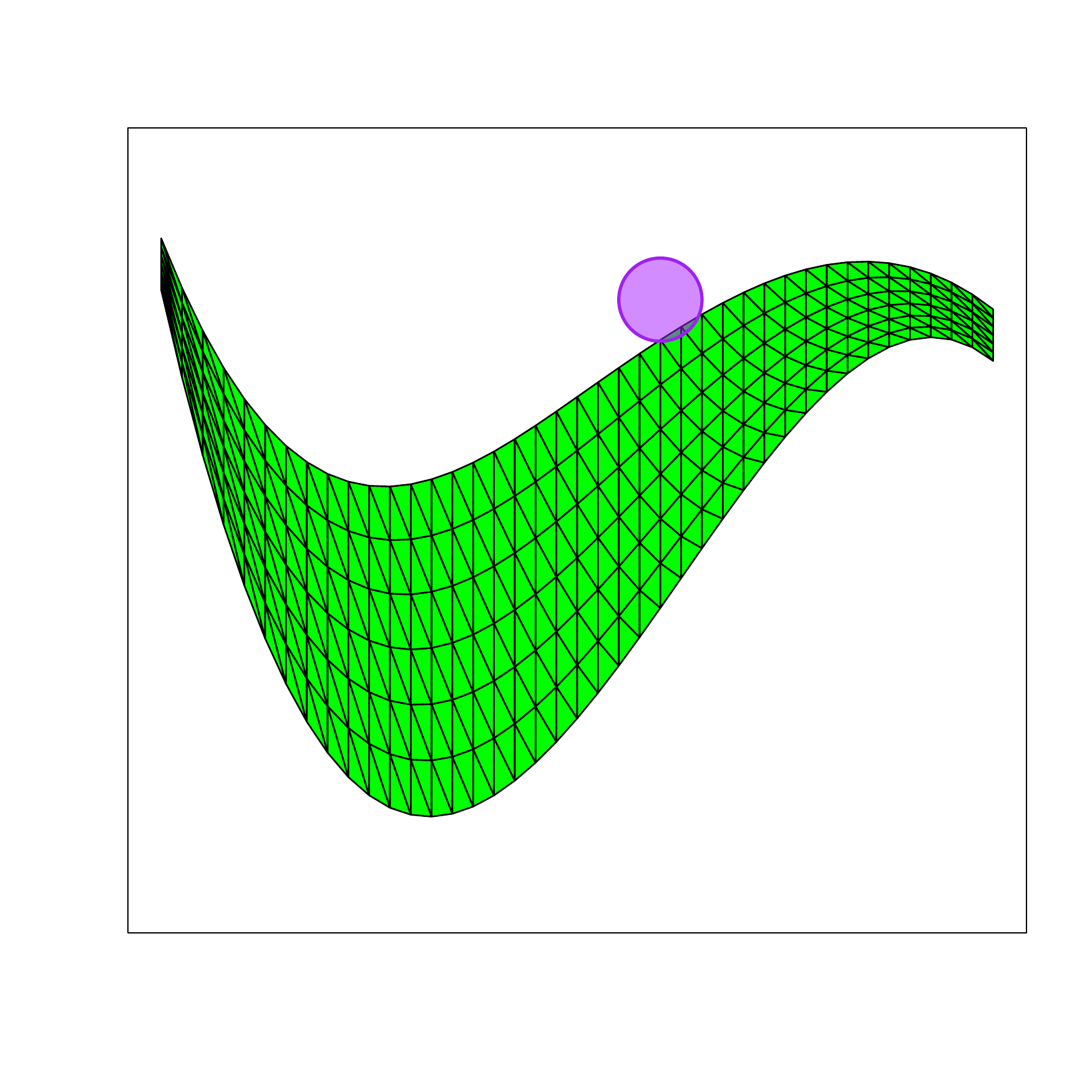}}

		\subfloat{\includegraphics[width=0.2\textwidth]{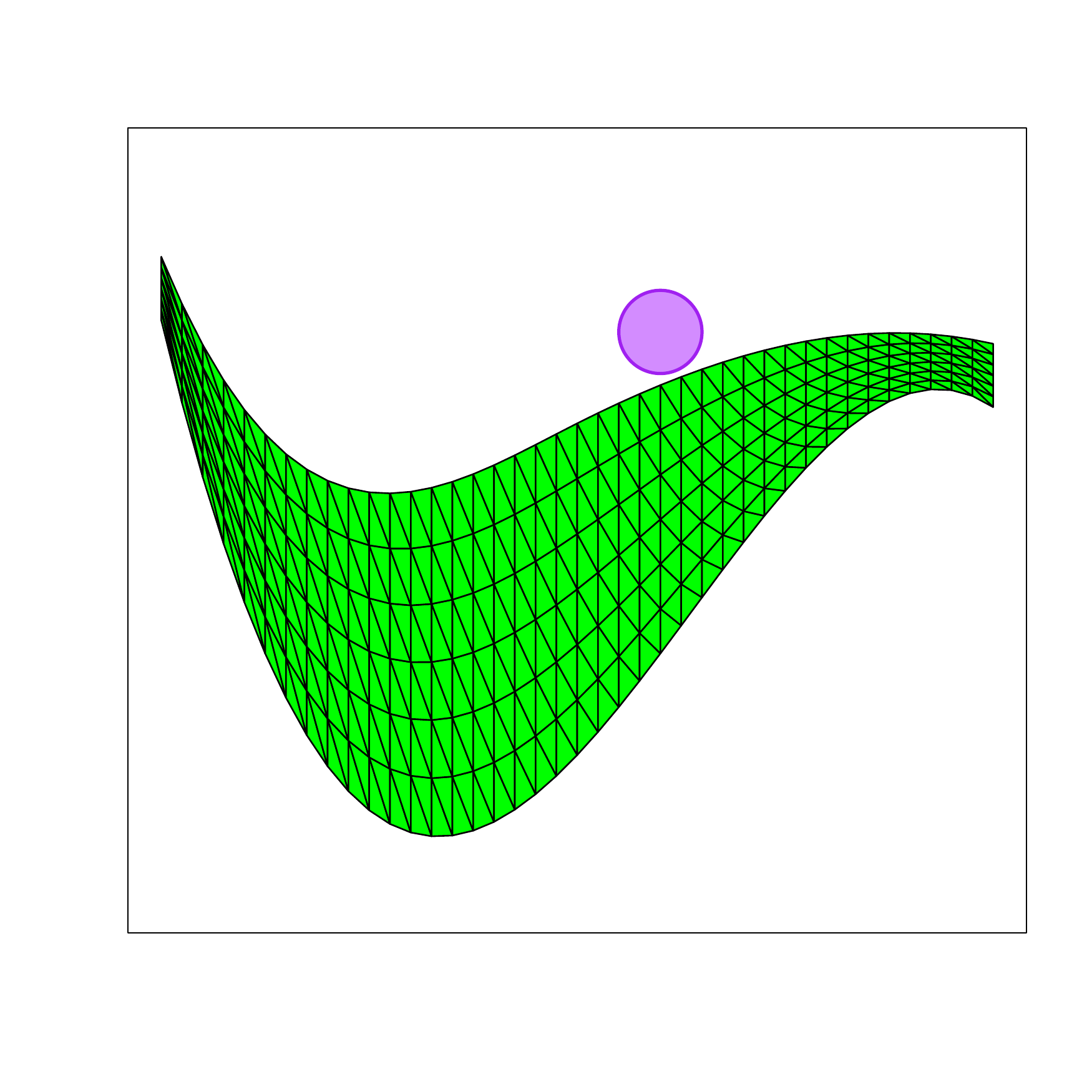}}
		\subfloat{\includegraphics[width=0.2\textwidth]{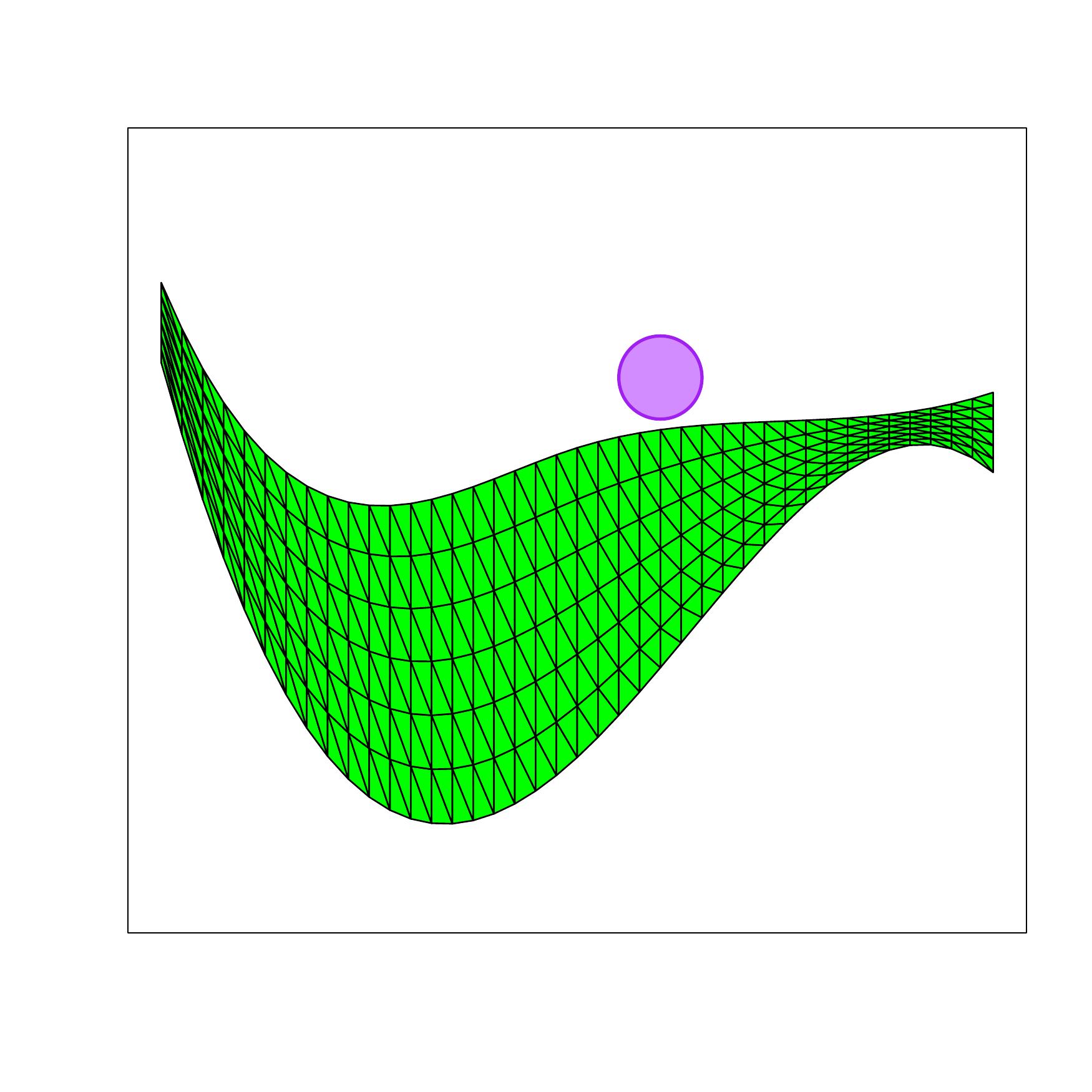}}
		\subfloat{\includegraphics[width=0.2\textwidth]{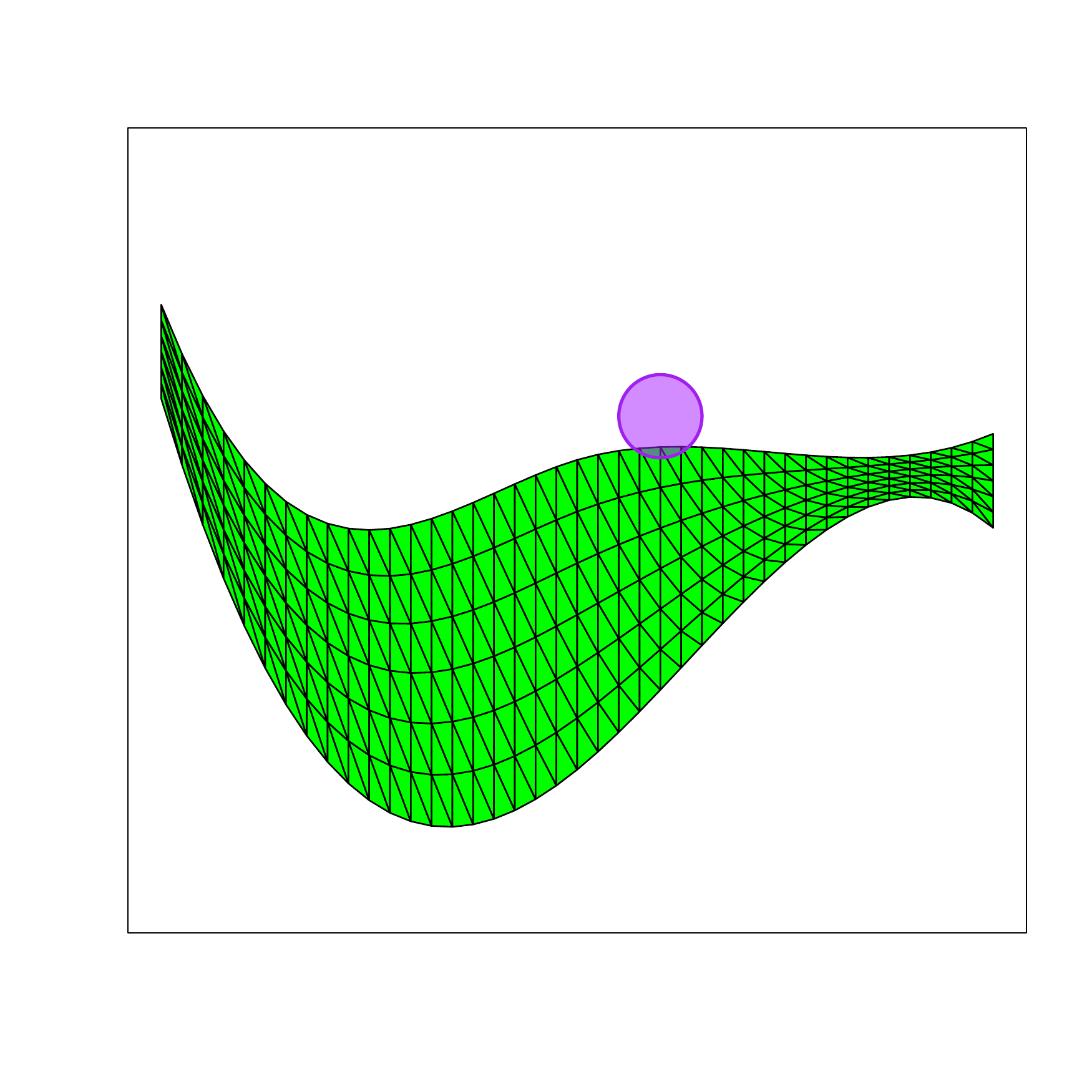}}
		\subfloat{\includegraphics[width=0.2\textwidth]{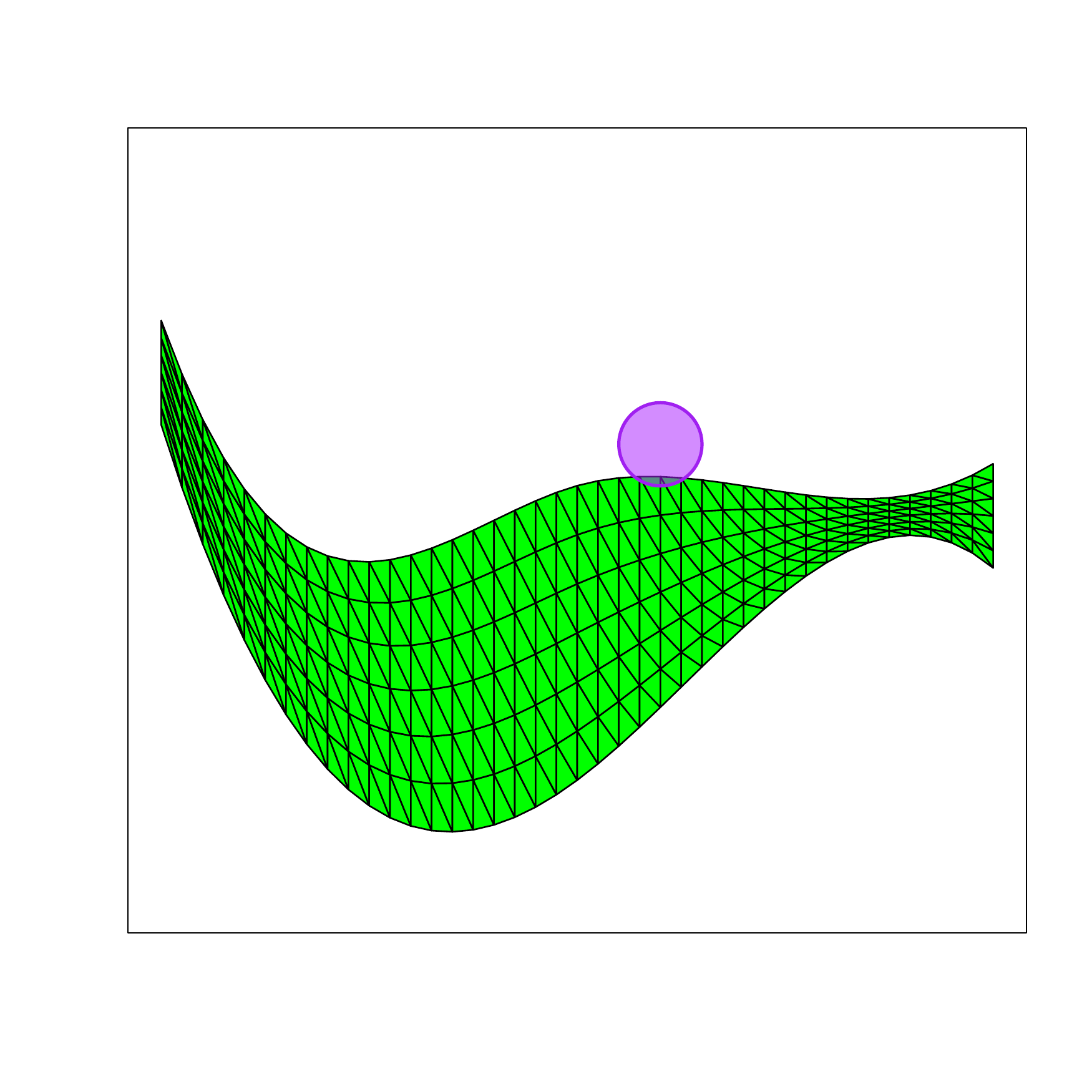}}
		\subfloat{\includegraphics[width=0.2\textwidth]{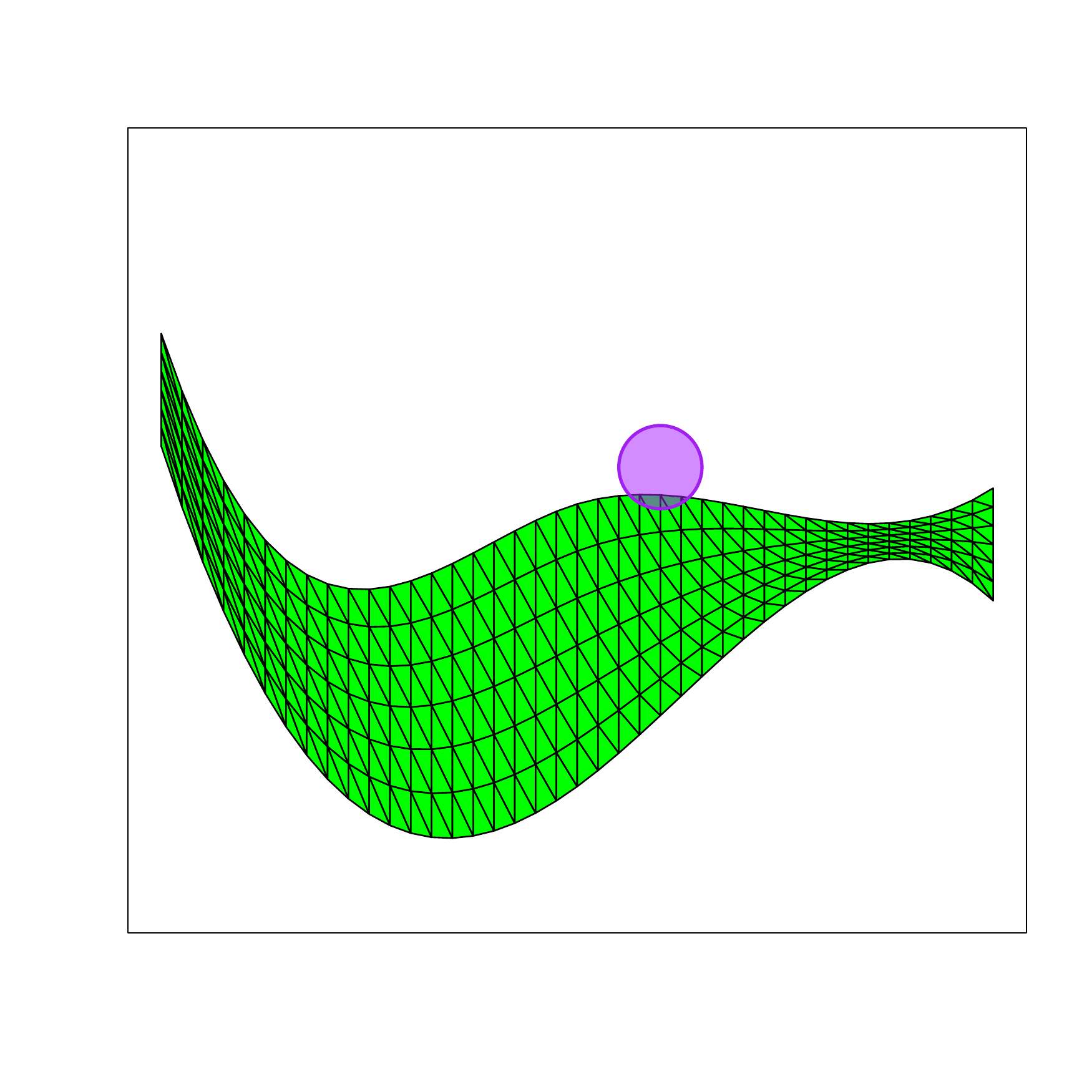}}

 	\end{center}		
	\caption{Test case 2 -- selected shapes computed during the \kk{momentum method} 
    starting in $q^0$ and ending in $q^{\text{HF}}$. 
    Recall that the plotted circle does not always coincide with the actual circle. \href{https://youtu.be/IhtRV9hX8AA}{See video.}\label{fig:TC2_shapeflow}}
\end{figure}

 \kk{Selected exemplary} shapes capturing the course of the \kk{discretized Hamiltonian flow} is illustrated in Figure~\ref{fig:TC2_shapeflow}. Here, as for the first test case, the first \otd{shape} depicts the initial shape $q^0$ and the \otd{last shape the solution $q^{\text{HF}}$, respectively}, while the remaining shapes are chosen in a way to best illustrate the Hamiltonian shape flow. 

\section{\label{sec:CO}Conclusions and outlook}
In this paper we consider explicitly non-con\-vex  problems in shape optimization, which previously attained little attention. In particular, we demonstrate that avoidance of occupied areas in an installation space creates additional local minima which render the shape optimization problem more difficult. To overcome such difficulties at least partially, we adapt momentum methods from \kk{the fields of non-convex optimization and} machine learning to shape optimization. Using a (discretized) dissipative Hamiltonian flow instead of a (discrtized) gradient flow, we numerically integrate shape flows in a way that the obstacles in the installation space are overcome. We provide numerical experiments for 2D mechanical shape optimization problems where the objective functions are given by a scalarization of material consumption, reliability, \kk{and constraint violation.}

We consider this work as a starting point for further developments in non-convex shape optimization. First, a better understanding of convergence to stationary points beyond local asymptotic Lyapunov stability is desirable. This can  be achieved either by a detailed analysis of global convergence of dissipative Hamiltonian flows to stationary points or by extension to port Hamiltonian flows where ports may be used to impose guarantees on the dissipated energy, which should be useful for global convergence to stationary points. More generally, the physics based intuition connected to Hamiltonian flows will contribute to the design of control strategies for shape optimization algorithms beyond the gradient descent paradigm.  

Also, the Hamiltonian perspective in connection with muti-objective optimization offers most interesting connections between topological properties of dynamical systems ('bifurcations', see e.g. \cite{williams1997chaos}) and the \kk{choice} of weighting parameters for scalarizations in multi-objective optimization. The notable stability of topological properties of dynamical systems away from bifurcating parameter settings might also be considered as a theoretical foundation for tracing methods as proposed in \cite{bolten2021tracing}. 

In this work, we introduce a generic definition of shape geometry using splines and thus achieve a finite dimensional parametrization once the spline basis is fixed. It would also be of interest to consider the infinite dimensional shape optimization setting, where we have to introduce momentum in the cotangent bundle of the infinite dimensional manifold of shapes. Clever choices of the skew-symmetric matrix $J$ and the dissipation $R$ could also be useful to provide the required smoothing for update steps in shape optimization in order to maintain the regularity class of the boundary.

We intend to come back to the indicated problems in future research.

\vspace{.3cm}

\paragraph{Acknowledgement} We thank Camilla Hahn, Marco Reese, Johanna Schultes, Volker Schulz and Michael Stiglmayr for interesting discussions.

\bibliographystyle{siamplain}
\bibliography{nonconvexSO}

\end{document}